\documentclass[11pt]{article}
\usepackage[margin=1in]{geometry}  
\usepackage{graphicx}
\usepackage{amssymb}             
\usepackage{amsmath}
\usepackage{mathtools}
\usepackage{slashed}            
\usepackage{amsfonts}              
\usepackage{amsthm}
\usepackage{breqn}
\usepackage{verbatim}              
\usepackage{enumitem}
\usepackage[numbers]{natbib}
\usepackage{mathtools}

\theoremstyle{definition}
\newtheorem{thm}{Theorem}[section]
\newtheorem{nota}[thm]{Notation}

\newtheorem{mydef}{Definition}[section]
\newtheorem{lem}[thm]{Lemma}

\newtheorem{cor}[thm]{Corollary}

\newtheorem*{rmk}{Remark}

\DeclareMathOperator{\id}{id}
\newcommand{\RR}{\mathbb{R}}      

\newcommand{\q}{\quad}
\newcommand{\p}{\partial}
\newcommand{\dd}{\mathfrak{D}}
\newcommand{\DD}{\mathcal{D}}
\newcommand{\curl}{\text{curl}\,}

\newcommand{\nab}{\nabla}
\newcommand{\n}{\hat{n}}
\newcommand{\lap}{\Delta}
\newcommand{\no}{\nonumber}
\newcommand{\di}{\text{div}\,}

\newcommand{\lleq}{\stackrel{L}{=}}

\newcommand{\cp}{\overline{\partial}{}}

\newcommand{\vol}{\text{vol}\,}

\newcommand{\kk}{\kappa}

\newcommand{\PP}{\mathcal{P}}
\newcommand{\uu}{\textbf{u}_0}
\newcommand{\ww}{\textbf{w}_0}
\newcommand{\vv}{\textbf{v}_0}
\newcommand{\pp}{\textbf{p}_0}
\newcommand{\qq}{\textbf{q}_0}
\newcommand{\wf}{\widetilde{\mathcal{F}}}

\newcommand{\rr}{\mathfrak{R}_\kk}
\numberwithin{equation}{section}

\usepackage[bookmarksnumbered,pdfpagelabels=true,plainpages=false,colorlinks=true,
            linkcolor=black,citecolor=black,urlcolor=black]{hyperref}
\usepackage{xcolor}
\newcommand{\bw}{\mathfrak{R}_\kappa}
\newcommand{\norm}[1]{||#1||}
\newcommand{\bou}{\Gamma}
\newcommand{\poly}{\mathcal{P}}

\begin{document}
\mathtoolsset{showonlyrefs=true} 

\title{On the Incompressible Limit for the Compressible Free-Boundary  Euler Equations with Surface Tension in the Case of a Liquid}

\author{Marcelo M. Disconzi\thanks{Vanderbilt University, Nashville, TN, USA. \texttt{marcelo.disconzi@vanderbilt.edu}}\,\,\thanks{MMD gratefully acknowledges support from NSF grant \# 1812826,
from a Sloan Research Fellowship provided by the Alfred P. Sloan foundation,
from a Discovery grant administered by Vanderbilt University, and from 
a Dean's Faculty Fellowship.} \and Chenyun Luo\thanks{Chinese University of Hong Kong, Shatin, NT, Hong Kong. \texttt{cluo@math.cuhk.edu.hk}}\,\,\,\thanks{Part of this work is done when CL was a faculty at Vanderbilt University.}}

\maketitle

\begin{abstract}
In this paper we establish the incompressible limit for the compressible free-boundary Euler equations with surface tension in the case of a liquid. Compared to the case without
surface tension treated recently in \cite{lindblad2018priori, luo2017motion}, the presence
of surface tension introduces severe new technical challenges, in that several
boundary terms that automatically vanish when surface tension is absent now
contribute at top order. Combined with the necessity of producing estimates
uniform in the sound speed in order to pass to the limit, such difficulties imply that
neither the techniques employed for the case without surface tension, nor estimates
previously derived for a liquid with surface tension and fixed sound speed, are applicable here.
In order to obtain our result, we devise a suitable sound-speed-weighted energy that takes into account
the coupling of the fluid motion with the boundary geometry. Estimates are closed by exploiting
the full non-linear structure of the Euler equations and invoking several geometric properties
of the boundary in order to produce some remarkable cancellations. We stress that we do
not assume the fluid to be irrotational.
\end{abstract}
\tableofcontents

\section{Introduction}
We consider the motion of a compressible liquid with free surface boundary in $\RR^3$. We use the notation $\DD_t$ to represent the bounded domain occupied by the fluid at each time $t$, whose boundary is advected by the fluid. The motion of the fluid is described by the compressible Euler equations
\begin{align}
\begin{cases}
\rho (\p_t u+\nab_u u) = -\nab p, &\text{in}\,\,\DD,\\
 \p_t \rho+\nab_u \rho + \rho\di u=0, & \text{in}\,\,\DD,\\
 p=p(\rho), &\text{in}\,\,\DD.
\end{cases} \label{EE}
\end{align}
Here, $\DD = \cup_{0\leq t\leq T} \{t\}\times \DD_t$, $u=u(t,x)$ is the velocity of the fluid, whereas $p=p(t,x)$ and $\rho=\rho(t,x)$ are the pressure and density, respectively. The density is bounded from 
below away from zero, i.e., $\rho\geq \text{constant} >0$. This condition on the density is what characterizes the fluid
as a liquid.
The initial and boundary conditions are
\begin{align}
\begin{cases}
\{x:(0,x)\in \DD\} =\DD_0,\\
u=u_0,\rho=\rho_0 \q\text{in}\,\{0\}\times\DD_0,
\end{cases}
\begin{cases}
(\p_t +\nab_u )|_{\p\DD}\in T(\p\DD),\\
p|_{\p\DD}= \sigma\mathcal{H}, \label{IBE}
\end{cases}
\end{align}
where $\mathcal{H}$ is the mean curvature of $\p\DD_t$,  $\sigma\geq 0$ is a constant,
and $T(\p\DD)$ is the tangent bundle of $\p\DD$ (the condition $(\p_t +\nab_u )|_{\p\DD}\in T(\p\DD)$
expresses the fact that the boundary moves with speed equal to the normal component of the
velocity). Finally, the equation of state is assumed to be a strictly increasing function of the density, i.e., 
\begin{align}
p=p(\rho),\q p'(\rho)>0.
\nonumber
\end{align}
We shall consider the specific equation of state given by \eqref{equation of state} in this manuscript. The unknowns in \eqref{EE}-\eqref{IBE} are $u,\rho$ and $\DD_t$, and hence, $\mathcal{H}$ and $p$ are function of the unknowns, and therefore, are not known a priori. 

Problem \eqref{EE}-\eqref{IBE} behaves significantly different depending on whether $\sigma =0$
or $\sigma > 0$. The former is known as the case without surface tension whereas the latter
is the case with surface tension, which is the situation treated in this manuscript. 
Our goal is to show that, for $\sigma >0$, the motion of a free-boundary incompressible fluid
with surface tension (corresponding to the idealized situation of a constant density fluid) is well-approximated by \eqref{EE}-\eqref{IBE} when an appropriate notion of compressibility is
very small. It is well-known that solutions to the incompressible equations,
written in section \ref{sec_background} below, cannot be obtained by simply
setting $\rho$ to a constant in \eqref{EE}-\eqref{IBE} (see, e.g., \cite{lindblad2018priori}).
The correct way
of setting the incompressible limit is via the 
fluid's sound speed introduced in section
\ref{section 1.3}.

The study of the incompressible limit has a long history in fluid dynamics, see section 
\ref{sec_background}. For the case of a motion with free-boundary, the only 
results we are aware are the recent works \cite{lindblad2018priori, luo2017motion} 
by Lindblad and the second author, both treating the case $\sigma=0$. In particular, to
the best of our knowledge this is the first proof of the incompressible limit
for the free-boundary compressible Euler equations with surface tension, i.e., $\sigma >0$.
Despite many new difficulties introduced by the presence of surface tension, 
which are discussed in section \ref{section 1.5}, it is important to consider the case $\sigma >0$
because real fluids have surface tension. Thus, this feature has to be incorporated
in the construction of more realistic models. We remark that we do \emph{not} assume that the fluid
is irrotational.

\subsection{Lagrangian coordinate and the reference domain}\label{section 1.1}
We introduce Lagrangian coordinates, under which the moving domain becomes fixed. Let $\Omega$ be a bounded domain
 in $\RR^3$. Denoting coordinates on $\Omega$ by $y=(y_1,y_2,y_3)$, we define $\eta:[0,T]\times \Omega\to \DD$ to be the flow of the velocity $u$, i.e., 
 \begin{align}
 \p_t \eta(t,y) & =u(t,\eta(t,y)),
 \nonumber
\\
\eta(0,y) &= y.
\nonumber
 \end{align}
We introduce the Lagrangian velocity, density and pressure, respectively, by $v(t,y):=u(t,\eta(t,y))$, $R(t,y):=\rho(t,\eta(t,y))$ and $q(t,y):=p (t,\eta(t,y))$. Therefore, 
\begin{align}
\p_t \eta = v. \label{flow}
\end{align}
For the sake of simplicity and clean notation, here we consider the model case when $\DD_0=\Omega = \mathbb{T}^2\times (0,1)$. We set 
$$
\Gamma_0 := \mathbb{T}^2\times \{x_3=0\},\q
\Gamma_1 := \mathbb{T}^2\times \{x_3=1\},
$$
so that $\Gamma:=\p\Omega = \Gamma_0\cup\Gamma_1$. Using a partition of unity, as in, e.g.,
\cite{coutand2013LWP0.5,DKT}, a general domain can be treated with the same tools
we shall present. Choosing $\Omega$ as above, however, allows us to focus on the 
real issues of the problem without being distracted by the cumbersomeness of the 
partition of the unity. We also note that one might  want to consider a situation more
akin the finite-depth water waves problem, where the bottom boundary, $\Gamma_0$, remains
fixed. This case requires only minor modifications from our presentation but, again, we believe
that this would be a distraction from the main problem.

Let $\p$ be the spatial derivative with respect to the spatial variable $y$. We introduce the matrix $a=(\p \eta)^{-1}$. This is well-defined since $\eta(t,\cdot)$ is almost $\text{id}$ (i.e., the identity diffeomorphism on $\Omega$) whenever $t$ is sufficiently small. Define the cofactor matrix
\begin{align}
A = Ja,
\nonumber
\end{align}
where $J=\det(\p \eta)$. Then, $A$ satisfies the Piola identity:
\begin{align}
\p_\mu A^{\mu \alpha} =0.
\nonumber
\end{align}
Here, the summation convention is used for repeated upper and lower indices, and in above and throughout, we adopt the convention that the Greek indices range over $1,2,3$, while the Latin indices
 range over $1$ and $2$. 

In terms of $v,R,q$ and $a$, the system \eqref{EE}-\eqref{IBE} becomes
\begin{align}
\begin{cases}
R\p_t v^\alpha + a^{\mu\alpha}\p_\mu q=0, &\text{in}\,\, [0,T]\times \Omega\\
\p_t R + R a^{\mu\alpha}\p_\mu v_\alpha =0, &\text{in}\,\, [0,T]\times \Omega\\
q=q(R), &\text{in}\,\, [0,T]\times \Omega\\
A^{\mu\alpha}N_\mu q +\sigma\sqrt{g}\lap_g\eta^\alpha=0, & \text{on}\,\, [0,T]\times\Gamma,\\
\eta(0,\cdot)=\id, \q R(0,\cdot)=R_0(=\rho_0), \q v(0,\cdot)=v_0,
\end{cases}
\label{E}
\end{align}
where $N$ is the unit outward normal to $\Gamma$, and $\lap_g$ is the Laplacian of the metric $g_{ij}$ induced on $\Gamma(t)=\eta(t,\Gamma)$ by the embedding $\eta$, i.e., 
\begin{align}
g_{ij}=\p_i \eta^\mu \p_j \eta_\mu,\q 
\lap_g(\cdot) = \frac{1}{\sqrt{g}}\p_i(\sqrt{g}g^{ij}\p_j(\cdot)),
\end{align}
where $g=\det g$. Since $\eta(0,\cdot)=\id$, the initial Eulerian and Lagrangian velocities 
(i.e., $u_0$ and $v_0$) agree. In addition, we also have $a(0,\cdot)=I$, where $I$ is the identity matrix. Finally, $J=\text{det}(\p \eta)$ satisfies
\begin{align}
\p_t J = Ja^{\mu\nu}\p_\mu v_\nu,\q [0,T]\times \Omega. \label{p_t J}
\end{align}
This, together with the second equation of \eqref{E} imply
\begin{align}
RJ=\rho_0, \q [0,T]\times \Omega, \label{RJ=rho_0}
\end{align}
and hence the first equation in \eqref{E} is equivalent to
\begin{align}
\rho_0 \p_t v^\alpha + A^{\mu\alpha}\p_\mu q=0,\q \text{in}\,\, [0,T]\times \Omega. \label{E mod}
\end{align}

\subsection{Background\label{sec_background}}
The study of the motion of a fluid has a long history in mathematics. In particular,  the study of free-boundary fluid problems has blossomed over the past decade or so.  However, much of this activity has focused on the study of the incompressible free-boundary Euler equations, i.e., 
\begin{align}
\begin{cases}
\beta\mathfrak{v}_t^\alpha+\mathfrak{a}^{\mu\alpha}\p_{\mu} \mathfrak{q}=0,
& \text{in}\,\,[0,T]\times \Omega
\\
\di \mathfrak{v}=0, & \text{in}\,\,[0,T]\times \Omega \\
\mathfrak{A}^{\mu\alpha}N_\mu \mathfrak{q}+\sigma\sqrt{\mathfrak{g}}\lap_{\mathfrak{g}}\tilde{\eta}^\alpha=0, & \text{on}\,\,[0,T]\times\Gamma,
\end{cases}
\label{incompressible}
\end{align}
where $\beta$ is a positive constant corresponding to the fluid's constant density, $\mathfrak{v}$ and $\mathfrak{q}$ are the incompressible
Lagrangian velocity and pressure, 
 $\mathfrak{a}=(\p \tilde{\eta})^{-1}$, $\mathfrak{A}=\det (\p\tilde{\eta})\mathfrak{a}$, where 
$\tilde{\eta}$ is the Lagrangian map associated with $\mathfrak{v}$. 

It is well-known that for the incompressible equations, $\mathfrak{q}$ is not determined by an equation of state. Rather, it is a Lagrange multiplier enforcing the constraint $\di \mathfrak{v} = 0$.  
The local well-posedness for the incompressible free-boundary Euler equations
has been studied by many authors, see \cite{bieri2017motion, christodoulou2000motion, coutand2007LWP, coutand2010LWP, disconzi2014limit,disconzi2016free, disconzi2017prioriI, ignatova2016local, kukavica2017local, lindblad2002, lindblad2005well,  lindblad2009priori,
nalimov1974cauchy,  SchweizerFreeEuler, shatah2008geometry, shatah2008priori, shatah2011local, 
 wu1997LWPww, wu1999LWPww, zhang2008free} and references therein. It is worth mentioning here that when $\DD_0$ is unbounded (with finite or infinite depth) and the velocity $\mathfrak{v}_0$ is irrotational (i.e., $\curl \mathfrak{v}_0=0$, a condition that is preserved by the evolution), this problems is called the water-waves problem, which has received a great deal of attention \cite{alazard2013global, alazard2013sobolev, alazard2018morawetz, deng2016global, germain2012GWP, germain2015GWP, harrop2017finite, hunter2016two, ifrim2014two, ifrim2015two, ifrim2017lifespan,  ionescu2018global, ionescu2014global, ionescu2015global, ionescu2016global,totz2012rigorous, wang2015global, wang20173, wang2018global, wu2009GWP, wu2011GWP}. 

However, the theory of the free-boundary compressible Euler equations is far less developed. It is known that for suitable initial data, the system \eqref{EE} modeling a liquid admits a local (in time) solution, e.g., \cite{coutand2013LWP0.5, disconzi2017prioriC, disconzi2018local, lindblad2003well, lindblad2005well', trakhinin2009local}, and for the gas model, the existence of a local solution was obtained in \cite{ coutand2010priori, coutand2011LWP1D, coutand2012LWP, jang2009LWPgas, jang2015LWPgas, luo2014well}.

In this paper we study how the solutions to \eqref{E} and \eqref{incompressible} are related. Intuitively,  one expects that the solution of \eqref{E} should converge to that of \eqref{incompressible} when the ``compressibility vanishes".   The proper way to define this problem is via the fluid's sound speed 
(see \eqref{sound} below), which corresponds to the speed of propagation of sound waves inside the fluid and captures the fluid's compressibility in that stiffer fluids have larger sound 
speed\footnote{This is an experimental fact, see, e.g., \cite{white1999fluid}.}. 

The incompressible limit problem consists in proving that if a sequence $(v_{0,\kappa}, R_{0, \kappa})$ of \textit{well-prepared} initial data for \eqref{E} converges to $(\mathfrak{v}_0, \beta)$,
 where $\mathfrak{v}_0$ is the initial data for the incompressible problem \eqref{incompressible}, and the sound speed at time zero diverges to infinity,  then the respective solution $(v, R)$ of \eqref{EE} converges to $(\mathfrak{v}, \beta )$, 
 where $\mathfrak{v}$ solves \eqref{incompressible}. Here, well-prepared initial data means that, in addition to satisfying the compatibility conditions, the initial data has to be tailored  
 to the above limit (see Theorem \ref{data}).  

The incompressible limit for the compressible Euler equations in a fixed domain (i.e., $\DD_t=\DD_0$ or the whole space) was established by
several authors under different assumptions, see \cite{alazard2006low, alazard2008minicourse, disconzi2017inL, ebin1977motion, ebin1982motion, klainerman1981singular, klainerman1982compressible, Metivier2011, Schochet1986The} and references therein. In addition, the incompressible limit for the compressible free-boundary Euler equations was solved by Lindblad and the second author in \cite{lindblad2018priori} with $\sigma=0$ in a bounded domain, and by the second author \cite{luo2017motion} in the same case but with unbounded domain.  To our best knowledge, the aforementioned works \cite{lindblad2018priori, luo2017motion} are the only known results in the study of the incompressible limit for equations \eqref{E}. In particular,  no result is available for the case with $\sigma>0$.  We will establish a priori estimates for \eqref{E} that are uniform in the sound speed (see sections \ref{section 3}-\ref{section 4}). In addition, we will construct a sequence of well-prepared data for \eqref{E} which converges to that of \eqref{incompressible} when the sound speed tends to infinity (see section \ref{section 5}). As a consequence,  we conclude the convergence of the compressible solution to the incompressible one by an Arzel\`a-Ascoil-type theorem.

\subsection{The sound speed}\label{section 1.3}
Physically, the sound speed is defined as $c=\sqrt{q^\prime \circ R}$. 
To set up the incompressible limit, it is conveninet to view the sound speed as a paramenter. 
As in \cite{disconzi2017inL, ebin1982motion}, we consider a family
$\{q_\kk (R)\}$  parametrized by $\kk \in [0,\infty)$, where 
\begin{equation}
\kk:= q_{\kk}'(R)|_{R=\beta}. \label{sound}
\end{equation}
Here, $'=\frac{d}{dR}$, and
\begin{equation}
q_\kk(R)=c_\gamma\kk (R^\gamma-\beta^\gamma),\q c_\gamma=\gamma^{-1}>0, \,\beta>0,\,\,\gamma\geq 1. \label{equation of state}
\end{equation} 
We slightly abuse terminology and call $\kk$ the sound speed.
In order to consider the incompressible limit, we view the density as a function of the pressure, i.e., $R_\kk=R_{\kk}(q)=[(c_\gamma \kk)^{-1} q+\beta^\gamma]^{1/\gamma}$, 
and we see that $R_\kk' (q)$ satisfies
\begin{equation}
\frac{1}{c_0}\mathfrak{R}_\kk \leq R'_\kk (q) \leq c_0 \mathfrak{R}_\kk, \label{equivlent R_kk}
\end{equation}
for some fixed constant $c_0>0$, where $\mathfrak{R}_\kk=(c_\gamma \kk)^{-\frac{1}{\gamma}}$.
Also,  for $0\leq k\leq 4$, we have that:
\begin{align}
|R_{\kk}^{(k)}(q)|\leq c_0,\q |R_{\kk}^{(k)}(q)|\leq c_0|R_\kk'(q)|^k\leq c_0|R_\kk'(q)| \label{R_kk assumption},\\
|q_{\kk}^{(k)}(R)| \leq c_0 |q_{\kk}'(R)|,
\label{q_kk assumption}
\end{align}
hold uniformly in $\kappa$.


\subsection{The main results}
\textbf{Notations.} All notations will be defined as they are introduced. In addition, a list of symbols is given at the end of this section for a quick reference. 

\mydef The $L^2$-based Sobolev spaces are denoted by $H^s(\Omega)$, with the corresponding norm denoted by $||\cdot||_s$; note that $||\cdot||_0=||\cdot||_{L^2(\Omega)}$. We denote by $H^s(\Gamma)$ the Sobolev space of functions defined on $\Gamma$, with norm $||\cdot||_{s,\Gamma}$.

\thm \label{main theorem 2} Let $\Omega=\mathbb{T}^2\times(0,1)$ and $v_{0,\kk}$ be a smooth\footnote{By ``smooth" we mean ``as smooth as necessary for the qualitative arguments (such as integration by parts) to go through." However, all of our quantitative estimates depend only on the Sobolev norms mentioned in Theorem \ref{main theorem 2}.} vector field. Let $\rho_{0,\kk}$ be a smooth function satisfying $\rho_{0,\kk}\geq c>0$ and $q_{0,\kk}$ be the associated pressure given by \eqref{equation of state}. Suppose that for some $\mathfrak{m}\in \mathbb{R}$ such that
\begin{align}
||v_{0,\kk}||_4,||v_{0,\kk}||_{4,\Gamma},||q_{0,\kk}||_4, ||q_{0,\kk}||_{4,\Gamma} \leq \mathfrak{m},\q \text{for all}\,\, \kk>0
\label{initial norm}
\end{align} 
holds.
Then there exist a $T>0$ and a constant $\mathfrak{M}$ such that any smooth solution $(v_{\kk}, R_{\kk})$ to \eqref{E} defined on the time interval $[0,T]$ satisfies 
 \begin{equation}
 \mathcal{N}(t)\leq \mathfrak{M},
 \end{equation}
 where
\begin{align}
\mathcal{N}=||v_\kk||_4^2+||\rr\p_tv_\kk||_3^2+||\rr\p_t^2v_\kk||_2^2+||\rr^{\frac{3}{2}}\p_t^3v_\kk||_1^2\no\\
+||R_\kk||_4^2+||\p_tR_\kk||_3^2+||\sqrt{\rr}\p_t^2R_\kk||_2^2+||\rr\p_t^3R_\kk||_1^2\no\\
+||\p_t v_{\kk}||_2^2+||\sqrt{\rr}\p_t^2 v_\kk||_1^2+||\p_t^2 R_\kk||_1^2+E,
\label{Nkk}
\end{align}
where $E$ is defined as Definition \ref{def E}.

The next theorem is a direct consequence of Theorem \ref{main theorem 2} together with
the  Arzel\`a-Ascoli theorem.

\thm \label{main theorem 1}Let $\mathfrak{v}_0\in H^{6.5}(\Omega)$ be a divergence free vector field and let $\mathfrak{v}$ be the solution to the incompressible free-boundary Euler equations \eqref{incompressible} with data $\mathfrak{v}_0$ defined on a small 
time interval $[0,T]$.
Let $(v_{0,\kk}, R_{0,\kk})\in H^4(\Omega)\times H^4(\Omega)$ be a sequence of initial data for the compressible free-boundary Euler equations \eqref{E} satisfying the compatibility conditions up to order $3$ (see section \ref{section 5.1} for a statement of the compatibility conditions). Furthermore, assume that $(v_{0,\kk}, R_{0,\kk})\to (\mathfrak{v}_0, \beta)$ in $C^2(\Omega)$ as $\kk\rightarrow\infty$ and that \eqref{initial norm} holds. Let $(v_\kk, R_\kk)$ be the solution for \eqref{E} with the equation of state \eqref{equation of state}. Then:
\begin{enumerate}
\item For $\kk$ sufficiently large, $(v_\kk, R_\kk)$ is defined on $[0,T]$.
\item $
(v_\kk, R_\kk)\to (\mathfrak{v},\beta)$ in $C^0([0,T],C^2(\Omega))$
after possibly passing to a subsequence.
\end{enumerate}
\rmk 
$\mathfrak{v}_0\in H^{6.5}(\Omega)$ is required so that the initial norms are uniformly bounded. We refer the proof of Theorem \ref{data construction} for details. 

 Finally, we need the following theorem to show that the data required in Theorem \ref{main theorem 2} and Theorem \ref{main theorem 1} exists.

\thm \label{data}
Let $\mathfrak{v}_0\in H^{6.5}(\Omega)$ be a divergence free vector field in $\Omega$. Then there exists initial data $(v_{0,\kk}, R_{0,\kk})\in H^4(\Omega)\times H^4(\Omega)$ satisfying the compatibility conditions up to order $3$ (see section \ref{section 5.1} for a statement of the compatibility conditions) such that $(v_{0,\kk}, R_{0,\kk})\rightarrow (\mathfrak{v}_0, \beta)$ in $C^2(\Omega)$ as $\kk\rightarrow \infty$, and \eqref{initial norm} holds. 

\nota For the sake of clean notations, we will drop the $\kk$-indices on $v_\kk, R_\kk, q_\kk$, i.e., we will denote $(v_\kk, R_\kk, q_\kk)=(v, R, q)$ when no confusion can arise.

\subsection{On existence of solutions\label{S:On_existence}}

In Theorems \ref{main theorem 2} and \ref{main theorem 1} we have assumed
that a solution is given in the stated function spaces, whereas in Theorem
\ref{data} we showed how to construct initial data for solutions in the corresponding
spaces without, however, establishing the existence of solutions.
In this section we show that existence of solutions in the spaces we use follow
from the existence result of \cite{coutand2013LWP0.5}, although such an existence
result and its corresponding estimates do not suffice to obtain the incompressibe limit, as we
also discuss further below.

We begin noticing that given a solution with regularity as in \cite{coutand2013LWP0.5},
for each fixed $\kappa$, the norms appearing 
in $\mathcal{N}$ are well-defined,
where $\mathcal{N}$ is given in equation \eqref{Nkk}. The issue is that the time of existence 
of the solutions obtained in \cite{coutand2013LWP0.5}, as well as the a priori bounds in \cite{coutand2013LWP0.5}, depend on
$\kappa$, whereas one needs bounds and a time interval that is uniform on $\kappa$
in order to pass to the limit $\kappa \rightarrow \infty$.

The crucial point is that while our estimates hold on a small time interval $[0,T]$,
the \emph{smallness of $T$} does not depend on $\kappa$ provided that
$\kappa$ is sufficiently large. 
In a nutshell, the logic to obtain solutions in the spaces where we take the incompressible
limit is the following: (i) \cite{coutand2013LWP0.5} is used to obtain, for each $\kappa$, a solution
defined on a time interval $[0,T_\kappa)$; (ii) We apply our estimates to show
that the solution from \cite{coutand2013LWP0.5} can be controlled on a time interval $[0,T)$ that is uniform on
the sound speed $\kappa$. This uniform control follows from the use of our weighted-in-$\kappa$
estimates (i.e., the estimates with $\mathfrak{R}_\kappa$-weights) and, in fact, cannot be obtained from the energy used in \cite{coutand2013LWP0.5}, as we also show below; 
(iii) A more or less standard continuation argument is then used to obtain that 
$T_\kappa \geq T$ for all $\kappa$ sufficiently large. In this way we obtain a family
of solutions parametrized by $\kappa$ and defined on a common time interval. (iv)
Our estimates show that on this common time interval the family of solutions
converges (up to a subsequence) to the incompressible solution.

\begin{rmk}
We stress that the uniformity of $T$ on the sound speed $\kappa$ comes
from the fact that we can close our estimate for $\mathcal{N}$ (defined in 
equation (1.13)) uniformly on $\kappa$ (for large $\kappa$). This can only
be done because of the use of $\mathfrak{R}_\kappa$-weights in our energy
which is, furthermore, tailored to the incompressible 
limit\footnote{In particular,
as we discuss in this section and in section \ref{section 1.5}, our energy is related to that used 
in \cite{coutand2013LWP0.5}, but it also differs from it in important aspects.}.
\emph{After} we have obtained such a uniform-in-$\kappa$ estimate, 
we can derive further estimates
which can in principle depend on $\kappa$. In fact, estimates of this type
are used below. They are harmless because they are used in arguments that only
require \emph{finiteness} of some quantities. However we again insist
that the entire argument given below relies on the fact that we are able to derive
estimates independent of $\kappa$ (or, more precisely, independent of $\kappa$ for all
$\kappa$ sufficiently large).
\end{rmk}

We will now elaborate on the argument
summarized in above. We will present its logic step-by-step, but for the
sake of brevity will not write down explicitly many of the estimates involved.
After that, we will show that this uniform control of $T$
that we obtained does not follow from the result in \cite{coutand2013LWP0.5}.

In what follows, denote by $E^{CHS}$ the energy used
in \cite{coutand2013LWP0.5}, i.e., equation (1.9) of \cite{coutand2013LWP0.5}. 

\paragraph*{Claim I. Continuation criteria.} 
We begin with the following statement. Let $(v,R)$ be a solution defined on a time
interval $[0,T_*)$ and with regularity given by the norms\footnote{To avoid confusion, 
we stress that by ``regularity
given by the norms" we mean that the maps belong to the function spaces
of the corresponding norms, but we do not mean finiteness of the corresponding
energy over the time interval. For example, if we say that $v$ has regularity given
by the norms $|| v ||_s + || \partial_t v ||_{s-1}$, we mean that $v \in H^s$ and
$\partial_t v \in H^{s-1}$, but we do not assert that the supremum in $t$ of 
$|| v ||_s + || \partial_t v ||_{s-1}$ is finite.} in $E^{CHS}$.
Set $\mathcal{M} := \sup_{0\leq t < T_*} E^{CHS}(t) $. We claim that if
$\mathcal{M}< \infty$, then the solution $(v,R)$ can 
be continued pass $T_*$. 

Suppose that $\mathcal{M} < \infty$.
Because $|| \partial_t v ||_3$ is controlled by $E^{CHS}$, the fundamental theorem
of calculus combined with $\mathcal{M} < \infty$ shows that $v \in C^0([0,T_*), H^3)$.
Let $\{t_\ell\}_{\ell=1}^\infty$ be a sequence of times such that $t_\ell \rightarrow T_*$.
Using again the fundamental theorem of calculus and 
the triangle inequality, we have
\begin{align}
|| v(t_{\ell + j}) - v(t_\ell) ||_3 \leq \mathcal{M} | t_{\ell + j} - t_\ell |,
\nonumber
\end{align}
showing that $v(t_\ell)$ is a Cauchy sequence in $H^3$ so it converges. Since this is
true for any sequence $t_\ell \rightarrow T_*$, we have that there exists
a $v_* \in C^0([0,T_*],H^3)$ that extends $v \in C^0([0,T_*),H^3)$.
Moreover, since $v(t_\ell)$ converges to $v_*(T_*)$ in $H^3$ and is 
bounded in $H^4$ (because the $H^4$-norm is controlled by $E^{CHS}$), 
we obtain that in fact $v_*(T_*) \in H^4$.
 Using the equations of motion, which  give $\partial_t R \sim \partial v$,
 and the fact that $|| v ||_4$ is controlled 
by $E^{CHS}$, we similarly obtain an extension of $R$ to the closed interval. 
The same argument also gives that the flow $\eta_*$ of $v_*$, whose $H^5$ norm is controlled
by $E^{CHS}$ on $[0,T_*)$, satisfies $\eta_*(T_*) \in H^5$.
Repeating exactly the same argument for the boundary norms in $E^{CHS}$ (i.e., 
the last sum of (1.9) in \cite{coutand2013LWP0.5} and the next-to-the-last term of (1.9) in \cite{coutand2013LWP0.5}), we finally
conclude that $v$ and $R$ extend to functions on the closed interval $[0,T_*]$
and that $E^{CHS}(T_*) < \infty$. We can now apply Theorem 1.6 of \cite{coutand2013LWP0.5}, 
which says that if $E^{CHS}(T_*) < \infty$, then a solution exists on $[T_*, T_* + \varepsilon)$
for some $\varepsilon>0$.

\paragraph*{Claim II. Control of $E^{CHS}$ from $\mathcal{N}$ for fixed $\kappa$.} 
Let $(v,R)$ be a solution defined on a time interval $[0,\mathcal{T})$ with the regularity
given by the norms in $E^{CHS}$. We will show that if 
$\sup_{0\leq t < \mathcal{T}} \mathcal{N}(t) < \infty$, then  
$\sup_{0\leq t < \mathcal{T}} E^{CHS}(t) < \infty$, where $\mathcal{N}$ is the quantity
introduced in \eqref{Nkk}.

We begin noticing that the solution $(v,R)$ has enough regularity so that
$\mathcal{N}$ is well-defined. Assume that 
$\mathcal{N}_0 := \sup_{0\leq t < \mathcal{T}} \mathcal{N}(t) < \infty$. This immediately
gives that $\sup_{0\leq t < \mathcal{T}} || \partial_t^k v ||_{4-k}$, $k=0,\dots, 4$,
is bounded in terms of $\mathcal{N}_0$. We remark that this bound, and the ones that follow in this part of the argument, may depend on $\kappa$. However, here
this is not a problem since we only want to show the finiteness of 
$\sup_{0\leq t < \mathcal{T}} E^{CHS}(t)$ for fixed $\kappa$.
Thus, we obtain that all terms in the first
sum of (1.9) in \cite{coutand2013LWP0.5} are controlled by $\mathcal{N}_0$,
 except for the case $a=0$, i.e., $|| \eta ||_5$.

We next bound $|| \eta ||_5$ by controlling its curl, divergence, and normal component.
 For the curl, we use the compressible Cauchy invariance, equation
\eqref{cauchy invarn}. We note that this requires having 
$\left. \operatorname{curl} v \right|_{t=0} \in H^4$, which is true 
in the assumptions of the existence result of \cite{coutand2013LWP0.5} 
(see equation (4.7) in \cite{coutand2013LWP0.5} and the discussion surrounding it). 
We therefore
obtain a bound of $|| \operatorname{curl} \eta ||_4$ in terms of $\mathcal{N}_0$
and $E^{CHS}(0)$. For the divergence, we apply estimate
(34) of \cite{disconzi2017prioriC}, except that instead of the $H^{2.5+\delta}$
norm on the LHS, we use the $H^{4}$ norm (it is not difficult to see that the same
argument as in \cite{disconzi2017prioriC} goes through with the $H^4$ norm on the LHS).
This gives control of $||\operatorname{div} \eta ||_4$ in terms of $\mathcal{N}$ and 
$\left. \operatorname{div} \eta \right|_{t=0}$, where the latter is smooth in view
of the initial condition $\eta(0)= \operatorname{id}$. Thus, we conclude that 
$||\operatorname{div} \eta ||_4$ can be controlled in terms of $\mathcal{N}_0$.
Finally, we need to control $\eta \cdot N$ on $\Gamma$. For this, we  
use the boundary condition in (1.3) and apply elliptic estimates. We remark that
the coefficients do not have enough regularity for an application of the ``standard"
elliptic estimates, but we can apply estimates with coefficients in Sobolev spaces
(see Theorem 4 and Remark 2 in \cite{milani1983regularity}). 
Invoking div-curl estimates, we conclude that we
can control $|| \eta ||_5$ in terms of $\mathcal{N}_0$ and $E^{CHS}(0)$.

\begin{rmk} The presence of $E^{CHS}(0)$ above comes from the fact that the energy
$E^{CHS}$ requires an extra derivative for $\eta$, which we do not include in
$\mathcal{N}$. 
In fact, we cannot include this term in $\mathcal{N}$, otherwise we would not be able
to close our estimates uniformly in $\kappa$, as we discuss in more detail
in section \ref{section 1.5}.
However, above we showed that if a solution with such extra differentiability
is given, then we can update our estimates to control such extra derivative in terms 
of $\mathcal{N}$, with bounds possibly depending on $\kappa$. 
A similar remark applies to the boundary terms
in the last sum of $E^{CHS}$
 controlled below, which are more regular
in the energy $E^{CHS}$ than in our $\mathcal{N}$.
\end{rmk}

The term $ v_{ttt} \cdot n \in H^1(\Gamma)$ appearing in $E^{CHS}$, where $n$ is
the unit outer normal to the moving boundary, is directly
controlled by the boundary term in our energy $E$ which enters  
in the definition of $\mathcal{N}$ (see Definition \ref{def E}; see
also Lemma \ref{prelim lemma b} for identities relating $n$ with the projection $\Pi$
appearing in $E$) with
$\mathfrak{D}^4=(\mathfrak{R}_\kappa)^2 \partial_t^4$ (again, with a bound
depending on $\kappa$). 

Finally, we need to bound the terms in the last sum of $E^{CHS}$ (equation (1.9) of \cite{coutand2013LWP0.5}).
For this, we time differentiate the boundary condition in (1.3) up to three times\footnote{Using
the boundary condition gives control over $v\cdot N$, whereas the corresponding
terms in $E^{CHS}$ are $v \cdot n$. But it is not difficult to see that control of the latter follows from control of the former; see the proof of Theorem \ref{thm estimates for closing}}, and apply
elliptic estimates for operators with coefficients in Sobolev spaces (see above).
We also need 
control of up to three time derivatives of $q$ restricted to $\Gamma$ in terms of 
$\mathcal{N}$.  Such bounds are immediately available (with constants depending
on $\kappa$) from the bounds for $R$ and its time derivatives provided by $\mathcal{N}$.
We thus conclude that the last sum in  $E^{CHS}$ can be controlled in terms
of $\mathcal{N}_0$.

From the foregoing, we conclude that
$\sup_{0\leq t < \mathcal{T}} E^{CHS}(t) < \infty$
if $\sup_{0\leq t < \mathcal{T}} \mathcal{N}(t) < \infty$, as desired.  

\paragraph*{Claim III. Uniformity of the interval where the a priori estimates hold.} 
This is just a restatement of Theorem \ref{main theorem 2}, but we include it here for clarity
of the presentation.
The time interval $[0,T]$ in Theorem \ref{main theorem 2} is 
uniform on $\kappa$ in the following sense: $T$ has to be
chosen sufficiently small, but the smallness of $T$ depends only on a fixed large
$\kappa_0$. In other words, Theorem \ref{main theorem 2} says that there exists a $\kappa_0$ 
and a $T = T(\kappa_0)$, such that if $\kappa \geq \kappa_0$ and
$(v,R)$ is a solution defined on $[0,T)$ and with the regularity
given by the norms in $\mathcal{N}$, then $\mathcal{N}\leq \mathfrak{M}$,
where the constant $\mathfrak{M}$ in Theorem \ref{main theorem 2} (which is a constant
depending on norms of the initial data).

\paragraph*{Claim IV. Existence of solutions in the spaces where one takes the incompressible
limit.}
The results of \cite{coutand2013LWP0.5}, in particular Theorem 1.6, say
that given data
such that $E^{CHS}(0) < \infty$, there exists a solution $(v_\kappa,R_\kappa)$ 
defined on a short
time interval with regularity given by the norms in $E^{CHS}$. 
Let $[0,T_\kappa)$ be the maximal interval where the solution $(v_\kappa,R_\kappa)$
exists and has the regularity given by the norms in $E^{CHS}$. We use
the subscript ${}_\kappa$ to indicate that the solution as well as the time interval in 
principle depend on the sound speed $\kappa$.
Let $T$ be given by Claim III (i.e., by Theorem 1.1) and $\kappa_0$ be as in Claim III.
We  will show that $T_\kappa \geq T$ for all $\kappa \geq \kappa_0$.

Suppose that $T_\kappa < T$.
We remark that the solutions $(v_\kappa,R_\kappa)$ have enough regularity
so that the estimates of our Theorem 1.1 can be applied, i.e., all quantities
entering in the definition of $\mathcal{N}$ (equation \eqref{Nkk}) are well defined
for the solutions $(v_\kappa,R_\kappa)$. Since $T_\kappa < T$, Claim
III implies that $\sup_{0\leq t < \mathcal{T_\kappa}} \mathcal{N}_\kappa(t) < \infty$,
where we write $\mathcal{N}_\kappa$ to emphasize that this corresponds
to the quantity $\mathcal{N}$ for the solution $(v_\kappa,R_\kappa)$. By Claim
II we then obtain $\sup_{0\leq t < \mathcal{T}} E_\kappa^{CHS}(t) < \infty$,
where we write $E_\kappa^{CHS}$ to emphasize that this corresponds to the energy
$E^{CHS}$ for the solution $(v_\kappa,R_\kappa)$ (as in Claim II, the resulting
bound on $\sup_{0\leq t < \mathcal{T}} E_\kappa^{CHS}(t)$ depends on $\kappa$,
but only the finiteness of this quantity matters here). By Claim I, the solution
$(v_\kappa,R_\kappa)$ can be extended pass $T_\kappa$, which contradicts
the maximality of $T_\kappa$.

\vskip 0.2cm
Thus, we obtained a family of solutions parametrized by $\kappa$, $\kappa \geq \kappa_0$
and defined on $[0,T)$; shrinking $T$ a bit if necessary we can consider the close
interval $[0,T]$. Moreover, the estimate of Theorem 1.1, $\mathcal{N}(t) \leq 
\mathfrak{M}$, holds on $[0,T]$ for each solution in this family, with $\mathfrak{M}$
independent of $\kappa$ in view of Theorem \ref{main theorem 2}.

\paragraph*{Existence of initial data compatible with the regularity of $E^{CHS}$.}
The constant $\mathfrak{M}$ in Theorem \ref{main theorem 2}
 depends on norms of the initial data.
It will be uniform on $\kappa$, for all $\kappa$ sufficiently large, if the corresponding
norms of the initial data are uniform on $\kappa$, as assumed in Theorem \ref{main theorem 2}. To show that this assumption is not empty, in Theorem \ref{data}, we constructed 
initial data satisfying such uniformity on $\kappa$. However, for the existence of solutions
given above, we actually need data such that $E^{CHS}(0) < \infty$. Since 
$E^{CHS}(0)$ requires more regularity than $(v_0,q_0) \in (H^4(\Omega) \cap 
H^4(\Gamma) ) \times ( H^4(\Omega) \cap H^4(\Gamma))$, which 
is what we stated in Theorem \ref{data}, we need to explain how data
satisfying $E^{CHS}(0) < \infty$ and that is uniform on the sound speed can be obtained.
This, however, follows from the proof of Theorem \ref{data}. Indeed, the data constructed
in Theorem \ref{data} is regular enough so that $E^{CHS}(0)$ is well-defined for it. 
However,
only the $( H^4(\Omega) \cap H^4(\Gamma))$ 
norms of this data are needed to be controlled uniformly on $\kappa$ for
Theorem \ref{main theorem 2}, so the statement of Theorem \ref{data}
is restricted to this situation.

\paragraph*{$E^{CHS}$ cannot be controlled uniformly on $\kappa$.} The foundation
of our result on the incompressible limit is the fact that we can derive estimates
that are uniform in $\kappa$ (for large $\kappa$). On the other hand, in order to obtain solutions to the equations
of motion with the desired regularity, we relied on \cite{coutand2013LWP0.5}. This raises the natural question
of whether the incompressible limit could not be obtained directly from the estimates
derived in \cite{coutand2013LWP0.5}. Here we show that this is not the case, i.e., that the energy
$E^{CHS}$ cannot be closed uniformly in $\kappa$ solely within the framework of \cite{coutand2013LWP0.5}.

The relevant fact is that the energy estimates in \cite{coutand2013LWP0.5} are non-uniform in $\kappa$ and diverge
when $\kappa \rightarrow \infty$. In particular, the interval of existence obtained 
in \cite{coutand2013LWP0.5} could in principle shrink to zero when $\kappa \rightarrow \infty$ (that this is \emph{not}
the case is what we showed above using our uniform-in-$\kappa$ estimates).

Let us now provide details. We will show that Proposition 4.1 in \cite{coutand2013LWP0.5} is not uniform in $\kk$. To see this, we take $q=\kk (R^2-\beta)$ and so $q'(R)=2\kk R$. Plugging the identity 
$R=\rho_0/J$ to the Euler's equations we obtain
\begin{align}
\rho_o v_t +\kk A^{\mu\alpha} \p_\mu (\rho_0^2 J^{-2})=0, \label{1}\\
R_t+R a^{\mu\alpha} \p_\mu v_\alpha=0.
\end{align}
Testing four time-derivatives of \eqref{1} against $\p_t^4 v$ in the $L^2$ and then integrating by parts yields the energy
\begin{align}
\sup_{t\in [0,T]} ||v_{tttt} (t)||_0^2 +\sup_{t\in [0,T]}||\sqrt{\kk}\p_t^4 J(t)||_{0}^2+\sup_{t\in [0,T]}||v_{ttt}\cdot n(t)||_{1,\Gamma}^2.
\label{CHS tangential energy}
\end{align}
However, to control this energy, we need to control 
\begin{align*}
\int_{\Omega} \kk \p_t^4 (\rho^2 J^{-2}) ([\p_t^4, a^{\mu\alpha}]\p_\mu v_\alpha).
\end{align*}
In \cite{coutand2013LWP0.5} this term is part of the error term $\mathcal{R}$ and it can be controlled directly by the energy by H\"older's inequality. But there is a mismatch 
of $\kk^{1/2}$ between this term and \eqref{CHS tangential energy}. In fact, this term is associated with our term $\mathcal{I}_3$
(see section \ref{section 3.3.3}), which requires our $\mathfrak{R}_\kk$-weighted energy
to be controlled.

\subsection{Strategy, organization of the paper, and discussion of the difficulties} \label{section 1.5}

In this section we overview the main arguments of the paper, summarize the main difficulties, 
and explain how they are confronted.

\subsubsection{Special cancellations\label{sec_special_cancellations}}

As mentioned, having $\sigma>0$ leads to several new  difficulties 
not present when $\sigma=0$. This can be immediately seen from the boundary terms appearing
in the energy estimates (see sections \ref{section B} and \ref{sec_B_weighted}), since all these
terms are proportional to $\sigma$ and, therefore, automatically vanish when $\sigma=0$.
(Incidentally, we do not set $\sigma$ to $1$ as it is customary but keep it
explicit in order to highlight all the terms that would be absent had $\sigma$ been zero.)
Not only are these terms present but, as we discuss below, they are some of the most
difficult terms to handle.  As a consequence, the methods used in the second author's previous papers
to study the problem with $\sigma = 0$ \cite{lindblad2018priori,luo2017motion} cannot be applied 
  when $\sigma>0$.  

At first sight one might think that the surface tension should help with closing a priori estimates
since it has a regularizing effect on the boundary. This regularization, however, it is not enough
to produce control of the velocity on the boundary. After differentiating the equations
with respect to $D^k$, where $D^k$ is a $k^{\text{th}}$ order derivative, possibly mixing space
and time derivatives, contracting with $D^k v$ and integrating by parts, one is left with a boundary
term that reads, schematically,
\begin{align}
\int_\Gamma D^k v D^k q \, dS.
\nonumber
\end{align}
It is not difficult to see that we can only hope to control this term by employing the boundary condition so that (again, schematically)
\begin{align}
\int_\Gamma D^k v D^k q \, dS \sim \int_\Gamma D^k v D^k (\Delta_g \eta) \, dS.
\label{schematic_bry}
\end{align}
The presence of the boundary Laplacian and the fact that $v = \partial_t \eta$ suggest that we
should integrate by parts in space and factor a $\partial_t$. Although this is the strategy, 
we end up with a commutator term that is not of lower order. This is because the coefficients
of $\Delta_g$ involve one derivative of $g$ which, in turn, involves one derivative of $\eta$
(so that the coefficients depend on as many derivatives of $\eta$ as the order of the equation).
Thus, commuting $D^k$ and $\Delta_g $ still leaves a top order term that \emph{cannot} be written
as a perfect derivative (in time or space) to be integrated away. Moreover, this top order term
does not seem to have any good structure. In fact, one should not expect such term to have
a good structure, since differentiating the  coefficients of $\Delta_g$ corresponds to differentiate
$g^{ij}$, and, thus, to take derivatives of some non-linear combinations of the components $g_{ij}$
and its determinant.

The above difficulties are overcome by observing some remarkable cancellations among the bad
top order terms in \eqref{schematic_bry}. Such cancellations are not visible in any way in the
expressions that appear by simply manipulating \eqref{schematic_bry}. Rather, they are
identified after some judicious and lengthy analysis that relies heavily on some
geometric properties, expressed in the form of several geometric identities, of the boundary.
The first cancellation appears in \eqref{cancellation_1}. The reader can check that the terms
that cancel out are top order and that there does not seem to be possible to bound them individually.
The second cancellation happens between a term in \eqref{cancellation_2_1} and
\eqref{cancellation_2_2}. This second cancellation is even more remarkable because the terms
involved come from completely different parts of $D^k \Delta_g \eta$: one from when all derivatives
fall on the coefficient $\sqrt{g} g^{ij}$ of $\Delta_g$, the other from when we integrate one derivative in $\Delta_g$ by parts.

We also need a special cancellation for interior terms. This comes from when we take 
$D^k$ of the first equation in \eqref{E} and all derivatives fall on $a$. Since the matrix
$a$ already involves one derivative of $\eta$, we find terms in $D^{k+1} \eta$, which have
one too many derivatives of the Lagrangian map. Exploiting the explicit structure of $a$, however,
we are able to show that, when appropriately grouped, these bad terms cancel each other after some careful integration by parts (see \eqref{cancellation} and what follows). 

As this point one may ask if all such cancellations are indeed necessary since 
a priori estimates for \eqref{E} have been derived in the literature. The relevant work in this
regard is \cite{coutand2013LWP0.5}. There, the authors construct initial data where
$\eta$ is everywhere one degree more differentiable than $v$, and then prove that this extra
regularity is propagated by the evolution. They rely on such extra regularity to 
close the estimates. However, this does not seem possible here because
 such an extra differentiability is not compatible with the $\rr$-weights we need to introduce in order
to obtain estimates uniform in the sound speed (see section \ref{sec_weights}). 

A crucial aspect of all the cancellations mentioned above is that they require the derivatives
$D^k$ to contain \emph{at least one time derivative.} As a consequence, only the Sobolev norms of
time-derivatives of $v$ on the boundary are controlled from the energy estimates (we remark
that the energy \emph{does} involve time derivatives of the variables; it does not seem possible
to close the estimates without time-differentiating the equations). To obtain control
of non-time differentiated $v$ on the boundary, we rely directly on the boundary condition which, after
a time derivative, produces an equation of the form $\Delta_g v = \dots$ which is amenable to 
elliptic estimates. (One might wonder why we do not take further time derivatives of the boundary
condition to obtain estimates for $\partial_t^k v$ on the boundary. The reason is that, as mentioned
above, $\Delta_g$ does not commute well with derivatives due to the dependence of the coefficients
on two derivatives of $\eta$, so that we obtain an equation of worsening structure with each derivative\footnote{Taking several time derivatives of the boundary condition, in particular, would lead to a source term that can only be bounded with
$\kappa$-dependent bounds, preventing us from closing the argument uniformly
on $\kappa$. See section \ref{sec_weights} below for more on the need for uniform bounds. Compare also
with the use of the boundary condition to derive estimates for boundary terms
in Claim II of section \ref{S:On_existence}, where
the resulting bounds depend on $\kappa$.}. 
However, for only one time derivative, the resulting equation still has some good structure
that can be used to derive estimates.)

\subsubsection{$\rr$-weighted estimates\label{sec_weights}}

Another difficulty to establish the incompressible limit is that one has to derive 
estimates that are uniform in the sound speed, since the goal is to take the sound speed to 
infinity. This is substantially different than estimates for \eqref{E} (with $\sigma > 0$ )
currently available \cite{coutand2013LWP0.5, disconzi2017prioriC}. Establishing the required
 uniform-in-$\kk$ a priori estimate does not seem to be possible solely by the methods used
to derive the currently available estimates. In particular, a crucial element to derive such uniform estimates
is the use of a non-linear wave equation satisfied by the density, whereas non-uniform-in-$\kk$
estimates have been proven without this wave equation.
In fact, the known a priori energy bounds rely heavily on the fact that when $\rr$ is bounded from below (as $\kk$ is bounded from above), $\p q \approx \p R$ and $||q||_r \approx ||R||_r$, which is a direct consequence of the equation of state.  In particular, the energy used in \cite{disconzi2017prioriC} controls $||\p_t^k q||_{3-k}$ for free as a lower order term. However, this fact no longer holds when $\rr \to 0$. Indeed, since $\p R = R' \p q$, $||R||_r$ is merely equivalent to $||\rr q||_r$; in other words, we have to take extra effort to control the full Sobolev norms of $\p_t^k q$.
In \cite{lindblad2018priori} and \cite{luo2017motion}, where $\sigma =0$, these norms are controlled by elliptic estimate. This relies on the fact that one is able to control $||q_t||_r$ by the $r$-th order energy $E_r$ since 
$$\p^r q_t \sim \cp^r q_t + \p^{r-2} q_t+\text{lower order terms},$$
where $\cp$ denotes derivatives tangent to the boundary. The first term, $\cp^r q_t$,
vanishes due to $q|_\Gamma=0$. However, this method does not work when $\sigma>0$, which is simply due to the fact that $q\sim \lap_g \eta$ on $\Gamma$, and so $\cp^r q_t \sim \cp^{r+2} v$ on the boundary which has two derivatives too many. 

To resolve the above difficulties, our energy is defined using the  $\rr$-weighted derivatives  $\dd^r$ ($1\leq r\leq 4$), where
\begin{align}
\dd =\cp, \p_t; \q
\dd^2 = \cp^2, \cp\p_t, \sqrt{\rr}\p_t^2;\q
\dd^3 =  \cp^2\p_t, \sqrt{\rr}(\cp\p_t^2), \rr\p_t^3;\\
\dd^4 = \rr(\cp^3\p_t), \rr(\cp^2\p_t^2), (\rr)^{\frac{3}{2}}(\cp\p_t^3), (\rr)^2\p_t^4.
\end{align}
The energy $E=E(t)$ is defined by employing these $\rr$-weighted derivatives, which is of the form:
\begin{align}
E = {\sum}_{1\leq \ell\leq 4}||\dd^\ell v||_{L^2(\Omega)}^2+{\sum}_{1\leq \ell\leq 4} \sqrt{\rr}||\dd^{\ell}q||_{L^2(\Omega)}^2+\sigma{\sum}_{1\leq \ell\leq 4} ||\Pi\cp\dd^{\ell}\eta||_{L^2(\Gamma)}^2+W,
\end{align}  
where $\Pi$ is the projection onto the normal to the moving boundary 
(see Lemma \ref{prelim lemma b}) and $W$ stands for the energy of the wave equation satisfied by $q$, which is defined in section \ref{section 2.3}-\ref{section 2.4}. 

The energy estimate for $E$ cannot be closed by itself; in fact, the energy estimate requires  control of 
\begin{align}||v||_4, ||\rr v_t||_3, ||\rr v_{tt}||_2, ||(\rr)^{\frac{3}{2}}v_{ttt}||_1,\label{weight v}\end{align}
and
\begin{align}
||R||_4, ||R_t||_3, ||\sqrt{\rr}R_{tt}||_2, ||\rr R_{ttt}||_1\label{weight R}.
\end{align}
These quantities are not part of the energy since $\dd^\ell$ for $\ell=1,2,3,4$ do not involve non-tangential derivatives, nor the full tangential spatial derivative $\cp^4$. 
Such missing derivatives, however, cannot be included in the energy because they would lead
to the presence of non-tangential derivatives on the boundary.
As a consequence, we need to estimate $E$ together with the quantities above in order to close the a priori estimate. This is done with the help of elliptic estimates.

We now schematically show how to get the correct $\rr$-weights for our energy, since they 
are crucial for the desired uniform-in-$\kk$ estimates. We differentiate the equations
\begin{equation}
R\p_t v_\alpha +q'(R)a^{\mu\alpha}\p_\mu R=0,\label{1st}
\end{equation}
and
\begin{equation}
\p_t R+Ra^{\mu\alpha}\p_\mu v_\alpha=0, \label{2nd}
\end{equation}
with respect to time. Since $R'=R'(q) = \frac{1}{q'(R)}$, equation \eqref{1st} implies 
\begin{align}
\p \p_t^k R \sim R'\p_t^{k+1} v; \label{trading R}
\end{align}
in other words, we can trade one (full) spatial derivative on $R$ by one time derivative of $v$ multiplied by $R'$. On the other hand, in view of the standard div-curl estimate (i.e., \eqref{div curl tang} in Appendix), $\p_t^k v$ is estimated via  $\di \p_t^k v$, $\curl \p_t^k v$ and $\p_t^k v\cdot N$. While in the reference domain $\Omega=\mathbb{T}^2\times (0,1)$,  $\p_t^k v\cdot N=\pm \p_t^k v^3$, which is almost $\Pi \p_t^k v$, where $\Pi$ denotes the projection to the normal direction, and hence this can be controlled by $E$. In addition, $\curl \p_t^k v$ is estimated via Cauchy invariance which can be treated by adapting the method introduced in \cite{disconzi2017prioriC}. Finally, the equation \eqref{2nd} yields
\begin{align}
a^{\mu\alpha}\p_\mu \p_t^k v_\alpha\sim \p_t^{k+1}R;
\end{align}
in other words, we can estimate $\di \p_t^k v$ using $\p_t^{k+1} R$. Hence, 
\begin{align}
\p ^4 v \xrightarrow{\di}\p^3 R_t\xrightarrow{\eqref{trading R}} R'\p^2\p_t^2 v \xrightarrow{\di} R'\p \p_t^3 R \xrightarrow{\eqref{trading R}} (R')^2 \p_t^4 v,\label{algo p^4 v}
\end{align} 
where $(R')^2 \p_t^4 v$ is part of $E$. In addition, we have
\begin{align}
R' \p^3 \p_t v \xrightarrow{\di} R'\p^2\p_t^2 R \xrightarrow{\eqref{trading R}}(R')^2\p \p_t^3 v \xrightarrow{\di} (R')^2\p_t^4 R.
\end{align}
This algorithm also provides
\begin{align}
R'\p^2\p_t^2 v \xrightarrow{\di} R' \p \p_t^3 R \xrightarrow{\eqref{trading R}} (R')^2 \p_t^4 v,\\
(R')^{\frac{3}{2}}\p\p_t^3 v \xrightarrow{\di} (R')^{\frac{3}{2}}\p_t^4 R.
\end{align}
Here, $(R')^{\frac{3}{2}} \p_t^4 R$ can be controlled directly by $E$ since it is equal to $(R')^{\frac{5}{2}} \p_t^4 q$ up to lower order terms. On the other hand, applying this algorithm starting from $\p^4 R$, we get
\begin{align}
\p^4 R \xrightarrow{\eqref{trading R}} R' \p^3 \p_t v\xrightarrow{\di}R' \p^2 \p_t^2 R \xrightarrow{\eqref{trading R}} (R')^2 \p \p_t^3 v \xrightarrow{\di} (R')^2\p_t^4 R,\\
\p^3 \p_t R \xrightarrow{\eqref{trading R}} R'\p^2\p_t^2 v \xrightarrow{\di} R'\p \p_t^3 R \xrightarrow{\eqref{trading R}} (R')^2 \p_t^4 v,\\
\sqrt{R'}\p^2\p_t^2 R \xrightarrow{\eqref{trading R}} (R')^{\frac{3}{2}} \p \p_t^3 v \xrightarrow{\di} (R')^{\frac{3}{2}} \p_t^4 R,\\
R'\p\p_t^3 R \xrightarrow{\eqref{trading R}} (R')^2 \p_t^4 v.\label{algo p p_t R}
\end{align}
The detailed analysis can be found in section \ref{section 4}. But the above algorithm provides good guideline for the choice of $\rr$-weights in \eqref{weight v} and \eqref{weight R} using \eqref{equivlent R_kk}, as well as in $\dd^r$.

\rmk 
The condition \eqref{equivlent R_kk} allows us to define the weighted Sobolev norms (e.g., \eqref{Nkk}) with constant $\rr$-weights. It is convenient to have constant weights for the boundary estimates in section \ref{section B} to avoid derivatives falling on $R'$.  In addition, the condition \eqref{equivlent R_kk} allows us to distribute $\rr$-weights in order to obtain an uniform control in $\kk$.

The definition of the $\rr$-weighted derivative $\dd^r$ allows us to control the highest order (i.e., $4$th order) mixed norms of $q$ directly by the energy. However, in order to pass to the incompressible limit, we have to control $||v||_4$ directly without $\rr$-weights,  and this requires the control of $||q_t||_2$. In section \ref{section 3.2}, we control $||q_t||_2$ by the elliptic estimate, which requires the control of $||q_t||_1$ first. This is indeed of lower order but we need to take extra effort to prove that they can be controlled uniformly as $\rr \to 0$.  In addition, we remark here that in \cite{disconzi2017prioriC}, the authors were able to close the a priori energy estimate in $H^3$. However, in our case,  the bound for $||q_t||_1$ require the control of $||v||_4$ and $||\eta||_4$. This is because control of $||\p q_t||_0^2$ requires integration by parts, which yields $||\cp^2 \eta||_{1.5,\Gamma}$ and $||\cp^2 v||_{1.5,\Gamma}$ at the top order, and these quantities require 
$H^4$ control of $v, \eta$.

\subsubsection{The initial data}
As with the estimates themselves, the initial data has to be constructed uniform in the sound
speed in order to allow the passage to the limit $\kk \rightarrow 0$.
This was done for $\sigma=0$ in  \cite{lindblad2018priori}, but that method relied heavily
on the fact that $q$ vanishes on the boundary when surface tension is absent. Instead, we employ the method used in \cite{coutand2013LWP0.5}: For each $1\leq k\leq 3$,  the data that satisfies the $k$-th order compatibility condition is obtained via solving an elliptic equation of order $2k$,  which is acquired by time differentiating the boundary condition $q=\sigma \mathcal{H}$ for $k$ times and then restrict at $t=0$,  where the previous $0, \cdots, k-1$-th compatibility conditions are served as the boundary conditions.     This construction process allows one to show that the initial data is uniformly bounded for all sound speed $\kk$, so that one can take the limit $\kk\rightarrow\infty$. 

\subsection*{Acknowledgments} We would like to thank Jared Speck for useful discussions.
We also would like to thank the anonymous referees for raising questions 
whose answers improved the quality of the manuscript.

\subsection{List of notations}
\begin{itemize}
\item $\nab$: Eulerian spatial derivative.
\item $\p$: Lagrangian spatial derivative.
\item $\cp$: Tangential spatial derivative. In particular, $\cp=(\p_1, \p_2)$ in $\Omega$ and we will 
emphasize that these derivatives are tangential by denoting $(\p_1, \p_2)=(\cp_1,\cp_2)$.  
\item $D$: Either $\cp$ or $\p_t$. 
\item $\Omega$ and $\Gamma$: The reference domain $(0,1)\times\mathbb{T}^2$ in Lagrangian coordinate, whose boundary $\p\Omega=\Gamma$.
\item The matrices $a$ and $A$: $a=(\p \eta)^{-1}$, and $A=Ja$, where $J=\det (\p\eta)$. 
\item $\kk$: The sound speed.
\item $\rr$: $\rr\approx R'_{\kk}\rightarrow 0$  as $\kk\to \infty$. 
\item $||\cdot||_s = ||\cdot||_{H^s(\Omega)}$ and $||\cdot||_{s,\Gamma} = ||\cdot||_{H^s(\Gamma)}$.
\item $P(\cdot)$: A smooth function expression in its arguments.   
\item $\lleq$: Equality modulo lower order terms that can be controlled appropriately.
\end{itemize}

\section{Preliminary results} \label{section 2}
In this section, we give some auxiliary results providing the bounds on the flow map $\eta$ and the matrix $a$. In addition, we record several facts, expressions and inequalities that will come in handy in the later sections.  These results will be employed in the proof of Theorem \ref{main theorem 2}.
\lem \label{prelim lemma a} Assume that $||v||_{L^{\infty}([0,T], H^4(\Omega))}+||R||_{L^{\infty}([0,T], H^4(\Omega))}\leq M$. Let $p\in [1,\infty)$, then there exists a sufficiently large constant $C>0$, such that if $T\in [0, \frac{1}{CM^2}]$ and $(v,q)$ is defined on $[0,T]$, the following statements hold:
\begin{enumerate}
\item $||\eta||_4 \leq C$.
\item $||a||_3 \leq C$.
\item $||a_t||_{L^p(\Omega)} \leq C||\p v||_{L^p(\Omega)}$, and $||a_t||_s \leq C||\p v||_s,\q 0\leq s\leq 3$.
\item $||\p_\alpha a_t||_{L^p(\Omega)} \leq C||\p v||_{L^{p_1}}||\p_\alpha a||_{L^{p_2}}+C||\p_\alpha\p v||_{L^p}$, where $\frac{1}{p}=\frac{1}{p_1}+\frac{1}{p_2}$.
\item $||a_{tt}||_s \leq C||\p v||_s||\p v||_{L^\infty}+C||\p v_t||_s, \q 0\leq s\leq 2$.
\item $||a_{ttt}||_s \leq C||\p v_t||_s||\p v||_{L^\infty}+C||\p v_{tt}||_s, \q 0\leq s\leq 1$.
\item $||\p_t^4 a||_{L^p(\Omega)} \leq C||\p v||_{L^p}||\p v||_{L^\infty}^2+C||\p v_t||_{L^p}||\p v_t||_{L^\infty}+C||\p v_{tt}||_{L^p}||\p v||_{L^\infty}+C||\p v_{ttt}||_{L^p}$.
\item $J\geq \frac{1}{2}$.
\item If $\epsilon$ is sufficiently small and for $t\in [0, \frac{\epsilon}{CM^2}]$, we have $||a^{\alpha\beta}-\delta^{\alpha\beta}||_3\leq \epsilon$, and $||a^{\alpha\mu}a^\beta_{\,\,\mu}-\delta^{\alpha\beta}||_3 \leq \epsilon$. In particular, the form $a^{\alpha\mu}a^\beta_{\,\,\mu}$ is elliptic, i.e.,  $a^{\alpha\mu}a^\beta_{\,\,\mu}\xi_\alpha\xi_\beta \geq C^{-1}|\xi|^2$.
\item $C^{-1}\leq R\leq C$.
\end{enumerate}
\begin{proof}
We refer \cite{disconzi2017prioriC} and \cite{ignatova2016local} for the detailed proof. 
We point out that the proof follows directly from the equations, interpolation, and the fundamental
theorem of calculus.
\end{proof}
We record here the explicit form of the matrix $a$ which will be needed.
\begin{align}
a= J^{-1}\begin{pmatrix}
\cp_2 \eta^2\p_3\eta^3-\p_3\eta^2\cp_2\eta^3\q \p_3 \eta^1\cp_2\eta^3-\cp_2\eta^1\p_3\eta^3\q \cp_2 \eta^1\p_3\eta^2-\p_3\eta^1\cp_2\eta^2 \\
\p_3 \eta^2\cp_1\eta^3-\cp_1\eta^2\p_3\eta^3\q \cp_1 \eta^1\p_3\eta^3-\p_3\eta^1\cp_1\eta^3\q \cp_1 \eta^1\cp_1\eta^2-\cp_1\eta^1\p_3\eta^2\\
\cp_1 \eta^2\cp_2\eta^3-\cp_2\eta^2\cp_1\eta^3\q \cp_2 \eta^1\cp_1\eta^3-\cp_1\eta^1\cp_2\eta^3\q \cp_1 \eta^1\cp_2\eta^2-\cp_2\eta^1\cp_1\eta^2
\end{pmatrix} \label{a}
\end{align}
Moreover, since $A=Ja$, and in view of \eqref{a}, we can write
\begin{align}
A^{1\alpha}=\epsilon^{\alpha\lambda\tau}\cp_2\eta_\lambda\p_3\eta_\tau,\q A^{2\alpha}=-\epsilon^{\alpha\lambda\tau}\cp_1\eta_\lambda\p_3\eta_\tau,\q
A^{3\alpha}=\epsilon^{\alpha\lambda\tau}\cp_1\eta_\lambda\cp_2\eta_\tau.\label{A}
\end{align}
Here, $\epsilon^{\alpha\lambda\tau}$ is the fully antisymmetric symbol with $\epsilon^{123}=1$. This representation will be used to create a special cancellation scheme that leads to control of the energy  when all derivatives fall on the cofactor matrix (recall the discussion in section
\ref{sec_special_cancellations}). 

We also need some geometric identities to treat the boundary terms in the energy estimate. We record these identities in the next lemma. 
\lem \label{prelim lemma b} Let $n$ be the outward unit normal to $\eta(\Gamma)$. Let $\tau$ be the tangent bundle of $\overline{\eta(\Omega)}$ and $\nu$ be the normal bundle of $\eta(\Gamma)$, the canonical projection is given by 
$$
\Pi^\alpha_\beta = \delta^\alpha_\beta-g^{kl}\cp_k \eta^\alpha \cp_l \eta_\beta,
$$
and on $\Gamma$ it holds that:
\begin{enumerate}
\item $-\Delta_g \eta^\alpha = \mathcal{H} \circ \eta \,  n^\alpha \circ \eta$.
\item $n \circ \eta = \frac{a^T N}{|a^TN|}$.
\item $J |a^T N| = \sqrt{g}$.
\end{enumerate}
Above, $a^T$ is the transpose of $a$.
Furthermore, setting $\hat{n}=n\circ \eta$, the following identities hold on $\Gamma$:
\begin{enumerate}[resume]
\item $\Pi_\beta^\alpha = \hat{n}_\beta\hat{n}^\alpha$.
\label{Pi_nn}
\item $\Pi_\lambda^\alpha\Pi_\beta^\lambda = \Pi_\beta^\alpha$.
\label{Pi_Pi_contraction}
\item  $\hat{n}_\alpha = \hat{n}_\tau \Pi^\tau_\alpha$.
\label{n_Pi_contraction}
\item $\sqrt{g}\lap_g \eta^\alpha = \sqrt{g}g^{ij}\Pi_\mu^\alpha\cp^2_{ij}\eta^\mu$.
\item $\p_t \hat{n}_\mu = -g^{kl}\cp_kv^\tau \hat{n}_\tau \cp_l \eta_\mu$.
\label{partial_t_n}
\item $\p_i \hat{n}_\mu = -g^{kl}\cp^2_{ik}\eta^\tau \hat{n}_\tau \cp_l \eta_\mu$.
\label{partial_i_n}
\item $\cp_i(\sqrt{g} g^{ik} ) = -\sqrt{g} g^{ij} g^{kl} \cp_i \cp_j \eta^\mu \cp_l \eta_\mu $.
\label{partial_i_ggij}
\item $\partial_t (\sqrt{g} g^{ij}) = \sqrt{g}(g^{ij}g^{kl} - 2 g^{lj} g^{ik}) \cp_k v^\lambda \cp_l 
\eta_\lambda$.
\label{partial_t_ggij}
\end{enumerate}
\begin{proof}
These identities are well-known. The interested reader can consult, e.g., 
\cite{disconzi2017prioriC} for their proof.  
\end{proof}
The equation of state $q=q(R)$ allows us to control $R' q$ and $R$ interchangeably:
\lem \label{q and R} Suppose $R':=R'(q)$ satisfies \eqref{R_kk assumption}, and let $\p$ be either $\p_t$ or $\p_\alpha$,  then for each $1\leq r\leq 4$, we have:
\begin{align}
|R' \p^r q| \lesssim |\p^r R|+{\sum}_{\substack{j_1+\cdots+j_k=r\\2\leq k\leq r}}|\p^{j_1}R|\cdots|\p^{j_k} R|.\label{R' p^r q}
\end{align} 
\begin{proof}
A direct computation yields:
\begin{align}
R' \p^r q = \p^r R + {\sum}_{\substack{j_1+\cdots+j_k=r\\2\leq k\leq r}}C_{j_1,\cdots, j_m, k}R^{(k)}\p^{j_1}q \cdots \p^{j_k}q,
\end{align}
and invoking \eqref{R_kk assumption} and the fact $R'\p q=\p R$, \eqref{R' p^r q} then follows.
\end{proof}

\subsection{The boundary condition\label{sec_the_bry_con}}
The identities of Lemma \ref{prelim lemma b} imply that 
the boundary condition
\begin{equation}
A^{\mu\alpha}N_\mu q +\sigma\sqrt{g}\lap_g\eta^\alpha=0,\q\text{on}\,\,\Gamma,
\label{boundary cond}
\end{equation}
can be expressed in the following equivalent ways:
\begin{enumerate}
\item $\sqrt{g}g^{ij}\cp_{ij}^2\eta^\alpha-\sqrt{g}g^{ij}g^{ij}\cp_k\eta^\alpha \cp_l \eta^\mu \cp_{ij}^2\eta_\mu = -\frac{1}{\sigma}A^{\mu\alpha}N_\mu q $, where $g^{kl}\cp_l \eta^\mu \cp_{ij}^2 \eta_\mu=\Gamma_{ij}^k$.
\item $\sqrt{g}g^{ij}\Pi_\mu^\alpha \cp_{ij}^2 \eta^\mu= -\frac{1}{\sigma}A^{\mu\alpha}N_\mu q $.
\item $q=-\sigma (A^{3\alpha}\n_\alpha)^{-1}\sqrt{g}g^{ij}\n_\mu \cp^2_{ij}\eta^\mu=-\sigma g^{ij}\n_\mu \cp^2_{ij}\eta^\mu$, 
since $(A^{3\alpha}\n_\alpha)^{-1}\sqrt{g}$
simplifies to $1$.
\label{q_cont_bry_con}
\end{enumerate}
These identities follow directly from the definition. Interested readers can consult  \cite{disconzi2017prioriC} for their proof.
The above expressions will be frequently used to deal with the boundary estimates. 

\subsection{The interpolation inequality}
Besides standard interpolation, we will also use the following interpolation inequality throughout this paper. 
\thm \label{interpolation} Let $u:\Omega\to \RR$ be a $H^1$ function. Then:
\begin{equation}
||u||_{L^4(\Omega)} \lesssim ||u||_0^{\frac{1}{2}}||u||_1^{\frac{1}{2}}.
\end{equation}
\begin{proof}
See Theorem 5.8 in \cite{Adams}.
\end{proof}

\subsection{The wave equations of order $3$ or less}
\label{section 2.3}
The second equation in \eqref{E} can be re-expressed as 
\begin{align}
a^{\mu\alpha}\p_\mu v_\alpha = -\frac{R' \p_t q}{R}, \label{divergence expression}
\end{align}
where $R'=R_\kk'(q)\sim \mathfrak{R}_\kk$ via assumption \eqref{R_kk assumption}. Identity \eqref{divergence expression} together with \eqref{E mod} yield, after commuting $\p_t^{r-1}$ for $1\leq r\leq 3$ and then $a^{\nu\alpha}\p_\nu$, that:
\begin{align}
JR'\p_t^{r+1} q - a^{\nu\alpha}A^\mu_{\,\,\alpha}\p_\nu\p_\mu \p_t^{r-1} q = \mathcal{F}_r, \label{wave equation}
\end{align}
where
\begin{align}
\mathcal{F}_r=
-{\sum}_{\substack{j_1+j_2=r\\j_1\geq 1}}\big(\p_t^{j_1}(JR')\big)(\p_t^{j_2+1}q)
+a^{\nu\alpha}(\p_\nu \rho_0)\p_t^r v_\alpha\no\\
+{\sum}_{\substack{j_1+j_2=r-1\\j_1\geq 1}}a^{\nu\alpha}\p_\nu(\p_t^{j_1}A^{\mu}_{\,\,\alpha}\cdot\p_\mu \p_t^{j_2}q)\nonumber\\
-\rho_0{\sum}_{j_1+j_2=r-1}(\p_t^{j_1+1} a^{\nu\alpha})(\p_t^{j_2}\p_\nu v_\alpha)
+ a^{\nu\alpha}(\p_\nu A^\mu_{\,\,\alpha})\p_\mu \p_t^{r-1}q. \label{F_r}
\end{align}
The wave equation \eqref{wave equation} yields an energy identity which is essential when estimating $||q||_2$ and $||q_t||_2$ in section \ref{section 3.2}:
\thm  For $1\leq r\leq 3$, let
\label{est W thm}
\begin{align}
W_r^2= \frac{1}{2} \int_\Omega \rho_0^{-1}(JR'\p_t^r q)^2\,dy+\frac{1}{2}\int_\Omega \rho_0^{-1}R'(A^{\nu\alpha}\p_\nu \p_t^{r-1}q)(A^\mu_{\,\,\alpha}\p_\mu \p_t^{r-1}q)\,dy\no\\
+\frac{\sigma}{2}\int_{\Gamma} \rr \sqrt{g}g^{ij}\Pi_\mu^\alpha(\cp_i\p_t^r \eta^\mu)(\cp_j \p_t^r \eta_\alpha)\,dS. \label{W}
\end{align}
Then, 
\begin{equation}
{\sum}_{1\leq r\leq 3} W_r^2 \leq  \epsilon P(\mathcal{N}) +\epsilon (||q||_2^2+||q_t||_2^2) +\PP_0 +\PP\int_0^t\PP, \q t\in[0,T], \label{est W}
\end{equation}
where $T>0$ is sufficiently small. 
\begin{proof}
See Appendix \ref{section B'}.
\end{proof}

\subsection{The $\rr$-weighted wave equations}
\label{section 2.4}
We consider the following $\rr$-weighted derivatives:
$$
\rr \p_t^3, \q \sqrt{\rr} \p_t^2\cp, \q \p_t \cp^2.
$$
Writing these derivatives as $\rr^{\ell}D^3$ ($\ell=1,\frac{1}{2},0$), and the identity \eqref{divergence expression} together with \eqref{E mod} yield, after commuting $\rr^{\ell}D^3$ and then $a^{\nu\alpha}\p_\nu$, that:

\begin{align}
\rr^\ell R' JD^3\p_t^{2} q - \rr^\ell a^{\nu\alpha}A^\mu_{\,\,\alpha}\p_\nu\p_\mu D^3 q = \widetilde{\mathcal{F}}, \label{wwave equation}
\end{align}
where
\begin{align}
\wf = -\rr^{\ell}[D^3\p_t, JR']\p_t q+\rr^\ell [D^3, \rho_0]\p_t (R^{-1}R'\p_t q)\no\\
+\rr^{\ell}a^{\nu\alpha}(\p_\nu \rho_0)D^3\p_t v_\alpha+\rr^{\ell}a^{\nu\alpha}\p_\nu\big([D^3, A^{\mu}_{\,\,\alpha}]\p_\mu q\big)
+\rr^{\ell}a^{\nu\alpha}\p_\nu\big([D^3,\rho_0]\p_t v_\alpha\big)\no\\
-\rr^{\ell}\rho_0[D^3\p_t, a^{\nu\alpha}]\p_\nu v_\alpha+\rr^{\ell}a^{\nu\alpha}(\p_\nu A^\mu_{\,\,\alpha})\p_\mu D^3 q.
\label{wF_r}
\end{align}
We need these $\rr$-weighted wave equations since their energies yield a better control of certain $\rr$-weighted energy terms. 

\thm  Let
\label{est W4 thm}
\begin{align}
W_4^2= \frac{1}{2} \int_\Omega \rho_0^{-1}\rr^{2\ell}(JR'D^3\p_t q)^2\,dy+\frac{1}{2}\int_\Omega \rho_0^{-1}\rr^{2\ell} R'(A^{\nu\alpha}\p_\nu D^3 q)(A^\mu_{\,\,\alpha}\p_\mu D^3 q)\,dy\no\\
+\frac{\sigma}{2}\int_{\Gamma} \rr^{2\ell+1} \sqrt{g}g^{ij}\Pi_\mu^\alpha(\cp_iD^3\p_t \eta^\mu)(\cp_j D^3 \p_t \eta_\alpha)\,dS. \label{W4}
\end{align}
Then, 
\begin{equation}
W_4^2 \leq  \epsilon P(\mathcal{N}) +\PP_0 +\PP\int_0^t\PP, \q t\in[0,T], \label{est W4}
\end{equation}
where $T>0$ is sufficiently small. 

\begin{proof}
See Appendix \ref{C}.
\end{proof}

\rmk The energy \eqref{W4} yields a better control of $q$ with $1/2$ less $\rr$-weights, e.g., when $D=\p_t$, $W_4$ controls $||\rr^2 \p_t^4 q||_0$ and $||\rr^{\frac{3}{2}} \p \p_t^3 q||_0$. The corresponding terms in $E$ control merely $||\rr^{\frac{5}{2}} \p_t^4 q||_0$ and $||\rr^{2} \p \p_t^3 q||_0$.  This observation is crucial to control $\mathcal{I}_3$ in section \ref{section 3.3} when $\dd^4 = \rr^{\frac{3}{2}}\cp \p_t^3$ or $\rr \cp^2 \p_t^2$. 

\subsection{The Cauchy invariance}
We conclude this section with a compressible version of the Cauchy invariance, which was introduced in \cite{disconzi2017prioriC}. 
\thm 
Let $(v, R)$ be a smooth solution to \eqref{E}, then 
\begin{align}
\epsilon^{\alpha\beta\gamma}\p_\beta v^\mu \p_\gamma\eta_\mu =\omega_0^\alpha+\int_0^t\epsilon^{\alpha\beta\gamma}a^{\lambda\mu}\p_\lambda q\p_\gamma\eta_\mu \frac{\p_\beta R}{R^2},\label{cauchy invarn}
\end{align}
for $t\in [0,T)$. Here, $e^{\alpha\beta\gamma}$ is the totally antisymmetric symbol with $\epsilon^{123}=1$ and $\omega_0$ us the vorticity at $t=0$.

\section{Energy estimates} \label{section 3}
In this section we provide estimates for $(v,q)$ and their time derivatives.  We shall make frequent use of the assumptions \eqref{equivlent R_kk}-\eqref{R_kk assumption} and of the two preliminary lemmas (i.e., Lemma \ref{prelim lemma a} and Lemma \ref{prelim lemma b}) in section \ref{section 2}
 throughout this section without mentioning them every time.

\nota \label{nota PP N}
Let $E$ be defined as in Definition \ref{def E}, and let 
\begin{align*}
\PP=P(||v||_4, ||\rr v_t||_3,  ||v_t||_2, ||\rr v_{tt}||_2, ||\sqrt{\rr} v_{tt}||_1, ||(\rr)^{\frac{3}{2}}v_{ttt}||_1, ||\rr v_{ttt}||_0, \\
 ||R||_4, ||R_t||_3, ||\sqrt{\rr}R_{tt}||_2, ||R_{tt}||_1, ||\rr R_{ttt}||_1, ||\sqrt{\rr} R_{ttt}||_0,\\ 
||\rr\Pi\cp^3 v_t||_{0,\Gamma}, ||\rr \Pi\cp^2  v_{tt}||_{0,\Gamma}, ||(\rr)^{\frac{3}{2}}\Pi \cp v_{ttt}||_{0,\Gamma}, ||\Pi \cp^2 v_t||_{0,\Gamma}, ||\sqrt{\rr}\Pi \cp v_{tt}||_{0,\Gamma}) 
\end{align*}
and $\PP_0=P(||\eta_0||_{7.5}, ||v_0||_4, ||v_0||_{4,\Gamma}, ||q_0||_4, ||q_0||_{4,\Gamma}, ||\di v_0|_{\Gamma}||_{3,\Gamma}, ||\lap v_0|_\Gamma||_{2,\Gamma})$, where we abbreviate 
$$
||\Pi  w||_{0,\Gamma}^2 = \int_\Gamma \Pi_\mu^{\beta} w ^\mu 
 \Pi_\beta^{\alpha} w_\alpha.
$$ 
Here (and throughout this paper), we use $P(\cdot)$ to denote a smooth function in its arguments. 
In addition, we define $\mathcal{N}$ to be
\begin{align}
\mathcal{N}(t)\! =\!  ||v||_4^2+||\rr v_t||_3^2+||\rr v_{tt}||_2^2+||(\rr)^{\frac{3}{2}}v_{ttt}||_1^2
+||R||_4^2+||R_{t}||_3^2\no\\
+||\sqrt{\rr}R_{tt}||_2^2+||\rr R_{ttt}||_1^2
+||v_t||_2^2+||\sqrt{\rr}v_{tt}||_1^2+||R_{tt}||_1+E. 
\end{align}
The rest of this section is devoted to prove:
\thm \label{energy estimate E}(Energy estimate for $E$) For sufficiently large $\kk>0$,  we have
\begin{align}
E(t) \leq \epsilon P(\mathcal{N}(t))+\PP_0+\PP\int_0^t \PP,\label{estimate E}
\end{align}
where $t\in [0,T]$ for some $T>0$ chosen sufficiently small, provide that the a priori assumption
\begin{align}
||\p \eta||_{L^{\infty}}+||\p^2 \eta||_{L^{\infty}}+||g^{ij}||_{L^{\infty}} \leq M,\label{apriori 1}
\end{align} 
hold. 
\nota Here and thereafter, we use $\epsilon$ to denote a small positive constant which may very from expression to expression. Typically $\epsilon$ comes from choosing the time sufficiently small (e.g., Lemma \ref{prelim lemma a} (9)) and the Young's inequality with $\epsilon$. When all estimates are obtained, we can fix $\epsilon$ sufficiently small in order to close the estimates.
\subsection{The energy identity for the Euler equations}
\nota(Weighted tangential mixed derivatives)
We let $\dd^r, r=1,2,3,4$ to be the mixed tangential differential operator defined as
\begin{align}
\begin{cases}
\dd =\cp, \p_t,\\
\dd^2 = \cp^2, \cp\p_t, \sqrt{\rr}\p_t^2,\\
\dd^3 =  \cp^2\p_t, \sqrt{\rr}(\cp\p_t^2),\rr\p_t^3,\\
\dd^4 = \rr(\cp^3\p_t), \rr(\cp^2\p_t^2), (\rr)^{\frac{3}{2}}(\cp\p_t^3), (\rr)^2\p_t^4.
\end{cases}
\label{D}
\end{align}

\nota \label{notation R} Here and in sequel, we use $\mathcal{R}$ to denote lower order terms whose time integral $\int_0^t \mathcal{R}$ can be controlled by the right hand side of \eqref{estimate E}.

\mydef\label{def E} For each fixed $1\leq r\leq 4$, let $E=\sum_{r=1,2,3,4}(E_r+W_r^2)$, where
\begin{align}
E_r = \frac{1}{2}\int_\Omega \rho_0(\dd^r v_\alpha)(\dd^r v^\alpha)\,dy+\frac{1}{2}\int_\Omega JR'R^{-1}(\dd^r q)^2\,dy\nonumber\\
+\frac{\sigma}{2}\int_{\Gamma}\sqrt{g}g^{ij}\Pi^\alpha_\mu (\cp_i \dd^r \eta^\mu)(\cp_j \dd^r\eta_\alpha)\,dS.  \label{E_r}
\end{align} 
Here,  $W_r^2$ ($1\leq r\leq 4$) is defined as \eqref{W} and \eqref{W4},  and $\Pi$ is the normal projection operator defined in Lemma \ref{prelim lemma b}.
\rmk
 We use throughout that
$\norm{\rr^\ell\Pi \cp^m \partial_t^l \eta }^2_{0,\bou}$ is comparable with the coercive term
coming from the boundary part of the energy. We use that 
$g^{ij}$ is almost the Euclidean metric to make this comparison.  For example, in the boundary estimates (section \ref{section B}) we control $\norm{\rr^2\Pi \cp \partial_t^3 v}_{0,\bou}^2$ by $E$. 

The energy defined above is derived by differentiating $\frac{1}{2}\int_\Omega R(\dd^r v_\alpha)(\dd^r v^\alpha)\,dy$ in time, invoking \eqref{E}, \eqref{p_t J}, \eqref{RJ=rho_0}, \eqref{E mod}, \eqref{R_kk assumption}, \eqref{divergence expression} and the Piola's identity
\begin{align}
\p_\mu A^{\mu\alpha}=\p_\mu (J a^{\mu\alpha})=0, \label{Piola}
\end{align}
which follows from a direct computation using \eqref{a}, we have:
\begin{align}
\frac{d}{dt}\frac{1}{2}\int_\Omega \rho_0(\dd^r v_\alpha)(\dd^r v^\alpha)\,dy = -\int_{\Omega}JR\big(\dd^r v_\alpha\big)\big(\dd^r(a^{\mu\alpha}\frac{\p_\mu q}{R})\big) \,dy\nonumber\\
=-\int_{\Omega}(\dd^r v_\alpha)\Big(\dd^r(A^{\mu\alpha}\p_\mu q)\Big) \,dy+\underbrace{\int_{\Omega}(\dd^r v_\alpha)\Big([\dd^r, RJ](a^{\mu\alpha}\frac{\p_\mu q}{R})\Big) \,dy}_{\mathcal{I}_1}\nonumber\\
=\int_{\Omega}( \dd^r \p_\mu v_\alpha)\Big( \dd^r(A^{\mu\alpha} q)\Big) \,dy\underbrace{-\int_{\Gamma}(\dd^r  v_\alpha)\Big( N_\mu \dd^r(A^{\mu\alpha} q)\Big) \,dy}_{BD}+\mathcal{I}_1\nonumber\\
=\int_{\Omega}(\dd^r \p_\mu v_\alpha)(A^{\mu\alpha}\dd^r q) \,dy+\underbrace{\int_{\Omega}(\dd^r \p_\mu v_\alpha)\Big([\dd^r, A^{\mu\alpha}] q\Big) \,dy}_{\mathcal{I}_2}+BD+\mathcal{I}_1. \label{3.6}
\end{align}
The term $\int_{\Omega}(\dd^r \p_\mu v_\alpha)(A^{\mu\alpha}\dd^r q) \,dy$ is equal to
\begin{equation}
\int_{\Omega}\dd^r (A^{\mu\alpha} \p_\mu v_\alpha)\dd^r q\,dy+\underbrace{\int_{\Omega}\big([\dd^r ,A^{\mu\alpha}] \p_\mu v_\alpha\big)\dd^r q\,dy}_{\mathcal{I}_3},
\end{equation} 
where, after invoking \eqref{divergence expression}, we obtain
\begin{align}
\int_{\Omega}\dd^r (A^{\mu\alpha} \p_\mu v_\alpha)\dd^r q\,dy=-\int_{\Omega} \dd^r( \frac{JR'\p_t q}{R})\dd^r q\,dy\nonumber\\
= -\int_{\Omega} JR'R^{-1}(\p_t\dd^r q) \dd^r q\,dy+\underbrace{\int_{\Omega}\big([\dd^r, JR'R^{-1}] \p_t q\big)\dd^r q\,dy}_{\mathcal{I}_4}.\label{3.8}
\end{align}
The first term in the second line of \eqref{3.8} is equal to $$-\frac{d}{dt}\frac{1}{2}\int_\Omega JR'R^{-1}(\dd^r q)^2\,dy+\mathcal{R}, $$ where the main term is moved to the left hand side of \eqref{3.6}.

 On the other hand, invoking the boundary condition $A^{\mu\alpha}N_\mu q =-\sigma\sqrt{g}\lap_g\eta^\alpha$, as well as the seventh identity in Lemma \ref{prelim lemma b}, $BD$ is equal to:
\begin{align}
BD=-\int_{\Gamma} \dd^r  v_\alpha\dd^r(A^{\mu\alpha}N_\mu q) \,dy\no\\
=\sigma\int_{\Gamma} \dd^r  v_\alpha \dd^r(\sqrt{g}\lap_g\eta^\alpha) \,dy
 =\sigma\int_{\Gamma}\dd^r v_\alpha \dd^r (\sqrt{g}g^{ij}\Pi^\alpha_\mu\cp^2_{ij}\eta^\mu)\,dS\no\\
 =\sigma\int_{\Gamma}\sqrt{g}g^{ij}\Pi^\alpha_\mu (\dd^r v_\alpha)( \dd^r\cp^2_{ij}\eta^\mu)\,dS
 +\underbrace{\sigma\int_{\Gamma}\dd^r v_\alpha \big([\dd^r, \sqrt{g}g^{ij}\Pi^\alpha_\mu]\cp^2_{ij}\eta^\mu\big)\,dS}_{\mathcal{B}_1}.\label{3.9}
\end{align}

Integrating by parts the first term in the very last line of \eqref{3.9}, we have
\begin{align}
\sigma\int_{\Gamma}\sqrt{g}g^{ij}\Pi^\alpha_\mu (\dd^r v_\alpha)( \dd^r\cp^2_{ij}\eta^\mu)\,dS=-\sigma\int_{\Gamma}\sqrt{g}g^{ij}\Pi^\alpha_\mu (\cp_i \p_t \dd^r \eta_\alpha) (\cp_j \dd^r\eta^\mu)\,dS\nonumber\\
-\underbrace{\sigma\int_\Gamma \cp_i (\sqrt{g}g^{ij}\Pi_\mu^\alpha)(\p_t \dd^r \eta_\alpha)(\cp_j \dd^r \eta^\mu)\,dS}_{\mathcal{B}_2}. \label{3.11}
\end{align}
The first term on the right hand side of \eqref{3.11} is equal to 
$$
-\frac{d}{dt}\frac{\sigma}{2}\int_{\Gamma}\sqrt{g}g^{ij}\Pi^\alpha_\mu (\cp_i \dd^r \eta^\mu)(\cp_j \dd^r\eta_\alpha)\,dS
+{\frac{1}{2}}\underbrace{\sigma\int_\Gamma \p_t(\sqrt{g}g^{ij}\Pi^\alpha_\mu)(\cp_i \dd^r\eta_\alpha)(\cp_j \dd^r \eta^\mu)\,dS}_{\mathcal{B}_3},
$$
where the main term is moved to the left hand side of \eqref{3.6}. Summing things up, we have shown:
\begin{align}
\frac{dE_r}{dt} = {\sum}_{1\leq j \leq 4} \mathcal{I}_j+{\sum}_{j=1,2,3}\mathcal{B}_j+\mathcal{R}.
\end{align}
Thus, Theorem \ref{energy estimate E} follows if the terms $\mathcal{I}_{1,2,3,4}$ and $\mathcal{B}_{1,2,3}$ can be controlled by the right hand side of \eqref{estimate E}, which shall be treated in sections \ref{section 3.3}-\ref{section B} below. However, before doing this, we need to control $||q||_2$ and $||q_t||_2$.  

\subsection{Bounds for $||q||_2$ and $||q_t||_2$}\label{section 3.2}
Since $\dd^r$ symbolizes both $\rr$-weighted and non-$\rr$-weighted derivatives,
we need to bound $||q||_2$ and $||q_t||_2$ in order to control $\mathcal{I}_3$. Also, the bound for $||q_t||_2$ is required to control $||v||_4$ in section \ref{section 4}.  Taking $X=\p q$ and $X=\p q_t$, $s=1$, the standard div-curl estimate \eqref{div curl tang} yields that we need to control the lower order terms $||\p q||_0$ and $||\p q_t||_0$.  We remark here that in the case when $\sigma=0$ (e.g., \cite{lindblad2018priori}), these terms are controlled via $||\lap q||_0$ and $||\lap q_t||_0$, respectively, after integrating by parts and applying the Poincar\'e's inequality. However, we need to work a bit harder in order to control these quantities when $\sigma>0$. 

\nota We write $X\lesssim Y$ to mean $X\leq CY$, where $C>0$ is a large constant. 
\nota We are going to identify $\mathcal{P}^n=\mathcal{P}$ ($n\geq 1$) by a slight abuse of notations. Also, when $0\leq t<1$,
$
(\int_0^t \PP)^n \leq t^{n-1}\int_0^t \PP^n \leq t\int_0^t \PP,
$
via Jensen's inequality.

\lem \label{estimate for F_1} Let $\mathcal{F}_r$ be defined as \eqref{F_r}. Assuming the a priori assumption \eqref{apriori 1} holds, then for sufficiently large $\kk>0$ (i.e., $\rr \ll 1$), we have:
\begin{align}
||\mathcal{F}_1||_0 \lesssim \epsilon\mathcal{N}+\PP_0+\PP\int_0^t\PP.
\end{align}
\begin{proof}
 First, invoking \eqref{p_t J} and the assumption \eqref{R_kk assumption}, we have:
 \begin{align*}
 ||\p_t(JR')(\p_t q)||_0
\lesssim \PP_0+\PP\int_0^t\PP.
 \end{align*}
 Second, invoking Lemma \ref{prelim lemma a}(1-4), since $\p_\mu q = R (a^{-1})_{\mu\beta}\p_t v^\beta $ and $\p_\nu \rho_0\lesssim \rr |\p_\nu q_0|\leq \epsilon |\p_\nu q_0|$ for sufficiently small $\rr$, we get:
\begin{align*}
||(\p_t  a^{\nu\alpha})(\p_\nu v_\alpha)||_0 +
||a^{\nu\alpha}(\p_\nu \rho_0)\p_t v_\alpha||_0+||a^{\nu\alpha}(\p_\nu A^\mu_{\,\,\alpha})\p_\mu q||_0 \no\\
\lesssim  \epsilon\mathcal{N}+\PP_0+\PP\int_0^t\PP.
\end{align*}

\end{proof}
\lem \label{estimate for F_r} Let $\mathcal{F}_r$ be defined as \eqref{F_r}, Assuming the a priori assumption \eqref{apriori 1} holds,  then for sufficiently large $\kk>0$ (i.e., $\rr<<1$), we have:
\begin{align}
||\mathcal{F}_2||_0 \lesssim \epsilon||q_t||_2+\epsilon(\sqrt{\mathcal{N}}+\mathcal{N})+\PP_0+\PP\int_0^t \PP.
\end{align}

\begin{proof}
 First, there is no problem to control ${\sum}_{\substack{j_1+j_2=2\\j_1\geq 1}}||\p_t^{j_1}(JR')(\p_t^{j_2+1}q)||_0$ appropriately when $j_1=1$ using \eqref{p_t J} and the assumption \eqref{R_kk assumption}. Moreover, when $j_1=2$, one writes $J=\rho_0 R^{-1}$ and then
$
||\p_t^2 (\rho_0 R^{-1} R') q_t||_0 = ||\rho_0 R' \p_t^2 (R^{-1}) q_t||_0
$
modulo controllable terms, where
\begin{align*}
 ||\rho_0 R' \p_t^2 (R^{-1}) q_t||_0 \lesssim ||\p_t^2 (R^{-1})||_1^{\frac{1}{2}}||\p_t^2 (R^{-1})||_0^{\frac{1}{2}}||R_t||_1^{\frac{1}{2}}||R_t||_0^{\frac{1}{2}}\leq \epsilon(\sqrt{\mathcal{N}}+\mathcal{N})+\PP_0+\PP\int_0^t \PP.
\end{align*}
Here, we have applied the interpolation inequality (i.e., Theorem \ref{interpolation}) and the fact $R'\p_t q= \p_t R$.  Second, invoking Lemma \ref{prelim lemma a}(1-6) we get:
\begin{align*}
{\sum}_{j_1+j_2=1}||(\p_t^{j_1+1} a^{\nu\alpha})(\p_t^{j_2}\p_\nu v_\alpha)||_0 \lesssim \epsilon\mathcal{N}+\PP_0+\PP\int_0^t\PP,
\end{align*}
and since $\p_\nu A^{\mu}_{\,\,\alpha} = O(\epsilon)$ for small time and $\p_\nu \rho_0\lesssim \rr |\p_\nu q_0|\leq \epsilon |\p_\nu q_0|$ for sufficiently small $\rr$, we have:
\begin{align*}
||a^{\nu\alpha}(\p_\nu \rho_0)\p_t^2 v_\alpha||_0+||a^{\nu\alpha}(\p_\nu A^\mu_{\,\,\alpha})\p_\mu \p_tq||_0 \lesssim \epsilon ||q_t||_2+\epsilon\sqrt{\mathcal{N}}+\PP_0+\PP\int_0^t \PP.
\end{align*}
Third, since $\p_\mu q = R (a^{-1})_{\mu\beta}\p_t v^\beta $, $||a^{\nu\alpha}\p_\nu(\p_t a^{\mu}_{\,\,\alpha}\cdot\p_\mu q)||_0$ can be controlled appropriately by interpolation. 

\end{proof}
\lem \label{esi lower order lemma} We have
\begin{align}
||\p q(t,\cdot)||_0^2+||\p q_t(t,\cdot)||_0^2 \leq \epsilon ||q_t(t,\cdot)||_2^2+\epsilon P(\mathcal{N})+W_3^2+\PP_0+\mathcal{P}\int_0^t\mathcal{P},\label{estimate lower order q and q_t}
\end{align}
for $t\in[0,T]$ where $T>0$ is chosen sufficiently small. 
\begin{proof}
It suffices to consider $||\p q_t||_0$ only. Integrating by parts yields:
\begin{equation}
||\p q_t||_0^2 = \int_{\Omega} (\p_\mu q_t)(\p^\mu q_t) = -\int_\Omega q_t \lap q_t+\int_\Gamma (N^\mu\p_\mu q_t) q_t,
\end{equation}
and so we need to bound $\int_\Omega q_t \lap q_t$ and $\int_\Gamma (N^\mu\p_\mu q_t) q_t$, respectively. 
\paragraph*{Bound for $\int_\Omega q_t \lap q_t$:} Since $t\in[0,T]$ and $T>0$ is small, as well as
$$
\lap q_t = (\delta^{\mu\nu}-a^{\mu\alpha}a^\nu_{\,\,\alpha})\p_\mu\p_\nu q_t +a^{\mu\alpha}a^\nu_{\,\,\alpha}\p_\mu\p_\nu q_t,
$$
Lemma \ref{prelim lemma a} implies that 
\begin{align}
\int_\Omega q_t(\lap q_t) \leq \epsilon||q_t||_2^2+ \int_\Omega q_t(a^{\mu\alpha}a^\nu_{\,\,\alpha}\p_\mu\p_\nu q_t).
\end{align}
Now, invoking the wave equation \eqref{wave equation} and Lemma \ref{estimate for F_r},  we have:
\begin{align}
\int_\Omega q_t(a^{\mu\alpha}a^\nu_{\,\,\alpha}\p_\mu\p_\nu q_t) = \int_\Omega R'q_t q_{ttt}-\int_\Omega (q_t\mathcal{F}_2)J^{-1}
\lesssim ||q_t||_0(||R' q_{ttt}||_0 + ||\mathcal{F}_2||_0)\no\\
\leq ||q_t||_0 \Big(W_3+\epsilon ||q_t||_2+\epsilon(\sqrt{\mathcal{N}}+\mathcal{N})+\PP_0+\mathcal{P}\int_0^t \mathcal{P}\big).
\end{align}
On the other hand, since $$||q_t||_0\leq ||\p q_t||_0+\int_\Omega q_t$$ by the Poincar\'e's inequality, if we let $Y=(0,0,y^3)$, then

\begin{align}
||q_t||_0 \leq ||\p q_t||_0+\int_\Omega \p_\mu Y^\mu q = ||\p q_t||_0-\int_\Omega Y^\mu \p_\mu q_t+\int_{\Gamma} N_\mu Y^\mu q_t \no\\
\leq C_{\vol \Omega}||\p q_t||_0+ \int_\Gamma y^3 q_t.\label{3.18}
\end{align}
To control the last integral $\int_\Gamma y^3 q_t$, time differentiating the boundary condition $$q=-\sigma g^{ij}\n^\mu \cp^2_{ij}\eta_\mu,\q \text{on}\,\,\Gamma$$  gives
\begin{align}
q_t = -\sigma g^{ij}\n^\mu \cp^2_{ij}v_\mu+\mathcal{R}_{q_t}, \q\text{on}\,\,\Gamma\label{q_t on Gamma}
\end{align}
where $\mathcal{R}_{q_t}$ consists of terms of the form either 
$$\sigma g^{ij}g^{kl}(\cp_k v^\tau\n_\tau\cp_l\eta_\mu)\cp_{ij}^2\eta^\mu \q\text{or}\q \sigma(\cp^i v_\nu)(\cp^j \eta^\nu)\n_\mu \cp^2_{ij}\eta^\mu.$$
Now, invoking Lemma \ref{prelim lemma a}, Lemma \ref{prelim lemma b}, and the a priori assumption \eqref{apriori 1},  we have:
\begin{align}
\int_\Gamma y^3 q_t\lesssim \epsilon \mathcal{N}+\PP_0+\mathcal{P}\int_0^t \mathcal{P}. \label{3.20}
\end{align}
 Wrapping these up, we get:
\begin{align}
\int_\Omega q_t\lap q_t \lesssim \epsilon ||q_t||_2^2+\epsilon ||\p q_t||_0^2+\epsilon P(\mathcal{N})+W_3^2+\PP_0+\mathcal{P}\int_0^t \mathcal{P}.
\end{align}

\paragraph*{Bound for $\int_\Gamma (N^\mu\p_\mu q_t) q_t$:} We have
\begin{align}
\int_\Gamma (N^\mu\p_\mu q_t) q_t \leq ||q_t||_{0,\Gamma}||\p_3 q_t||_{0,\Gamma} \leq C(\epsilon^{-1})||q_t||_{0,\Gamma}^2+\epsilon ||\p q_t||_{0,\Gamma}^2.
\end{align}
Here,  we bound $\epsilon||\p q_t||_{0,\Gamma}^2$ by $\epsilon||q_t||_2^2$ using the trace lemma, which is part of the right hand side of \eqref{estimate lower order q and q_t}. On the other hand, invoking \eqref{q_t on Gamma}, we have:
\begin{align}
||q_t||_{0,\Gamma}^2 \lesssim \epsilon(\mathcal{N}^2+\mathcal{N})+ \PP_0+ \PP\int_0^t\PP.
\end{align}
To see this, note that in $||q_t||_{0,\Gamma}^2$, the top order term is $\sqrt{g}g^{ij} \hat{n}^\mu 
\cp_{ij}^2 v_\mu$. Using the trace inequality, it suffices to bound
 $||\sqrt{g}g^{ij}\n^\mu \p\cp^2_{ij}v_\mu||_0^2$. 
 We control this top order term by the Young's inequality, which leads to the appearance of $\epsilon \mathcal{N}^2$.  In addition, the lower order terms are controlled by $\epsilon\mathcal{N}+ \PP_0+ \PP\int_0^t\PP$ using the interpolation. 
 
Hence,
\begin{align}
\int_\Gamma (N^\mu\p_\mu q_t) q_t \lesssim \epsilon ||q_t||_{2}^2+\epsilon(\mathcal{N}^2+\mathcal{N})+\PP_0+\PP\int_0^t\PP.
\end{align}
Therefore, 
\begin{align}
||\p q_t||_0^2 = -\int_\Omega q_t\lap q_t+\int_\Gamma q_t(N^\mu \p_\mu q_t)\no\\
\lesssim \epsilon ||q_t||_2^2+\epsilon(\mathcal{N}+\mathcal{N}^2)+W_3^2+\PP_0+\mathcal{P}\int_0^t \mathcal{P}.
\end{align}
In addition, we are able to control $||\p q||_0^2$ appropriately by integrating $||\p q_t||_0^2$ in time, 
which, together with the estimate for $||\p q_t||_0^2$, conclude the proof of \eqref{estimate lower order q and q_t}.
\end{proof}

In fact, the above proof implies the control for the lowest order norms $||q||_0$ and $||q_t||_0$.
\cor \label{esi lowest order q and q_t lemma} We have
\begin{align}
||q||_0^2+||q_t||_0^2 \lesssim ||\p q||_0^2+||\p q_t||_0^2+\epsilon\mathcal{N}+\PP_0 + \PP\int_0^t \PP.\label{estimate lowest order q and q_t}
\end{align}

\begin{proof}
Let $Y=(0,0,y^3)$, the Poincar\'e's inequality implies
\begin{align}
||q||_0 +||q_t||_0\lesssim ||\p q||_0 +||\p q_t||_0+\int_\Omega \p_\mu Y^\mu q+\int_\Omega \p_\mu Y^\mu q_t.
\end{align}
Now, we proceed as in \eqref{3.18}-\eqref{3.20} and get
\begin{align}
\int_\Omega \p_\mu Y^\mu q +\int_\Omega \p_\mu Y^\mu q_t \lesssim ||\p q||_0+||\p q_t||_0+\epsilon\sqrt{\mathcal{N}}+\PP_0+ \PP\int_0^t\PP,
\end{align}
and hence \eqref{estimate lowest order q and q_t} follows after squaring the above estimate.
\end{proof}

\thm \label{esi ||q||_2 and ||q_t||_2 theo} We have 
\begin{align}
||q(t,\cdot)||_2^2+||q_t(t,\cdot)||_2^2 \lesssim \epsilon P(\mathcal{N})+ \PP_0+ \PP\int_0^t\PP, \label{estimate ||q||_2 and ||q_t||_2}
\end{align}
for $t\in[0,T]$ where $T>0$ is chosen sufficiently small. 
\begin{proof}
It suffices to control $||q_t||_2^2$ by the right hand side of \eqref{estimate ||q||_2 and ||q_t||_2} since the control of $||q||_2^2$ follows from  time integrating $||q_t||_2^2$. To control $||q_t||_2^2$, it suffices to consider $||\p q_t||_1^2$ only thanks to Lemma \ref{esi lower order lemma} and Corollary \ref{esi lowest order q and q_t lemma}. Now, invoking the div-curl estimate \eqref{div curl tang} with $X=\p q_t$ and $s=1$, we have
\begin{align}
||\p q_t||_1^2 \lesssim ||\lap q_t||_0^2 + ||\cp q_t||_{0.5,\Gamma}^2+||\p q_t||_0^2.
\end{align}
\paragraph*{Bound for $||\lap q_t||_0^2$:} Invoking Lemma \ref{prelim lemma a}, since $t\in[0,T]$ and $T$ is sufficiently small, we have
\begin{align}
||\lap q_t||_0^2 \leq ||a^{\mu\alpha}a^\nu_{\,\,\alpha}\p_\mu\p_\nu q_t||_0^2+||(\delta^{\mu\nu}-a^{\mu\alpha}a^\nu_{\,\,\alpha})\p_\mu\p_\nu q_t||_0^2 \leq ||a^{\mu\alpha}a^\nu_{\,\,\alpha}\p_\mu\p_\nu q_t||_0^2+\epsilon ||\p^2 q_t||_0.
\end{align}
 Furthermore, the wave equation \eqref{wave equation}, as well as Lemma \ref{estimate for F_r} yield:
\begin{align}
||a^{\mu\alpha}a^\nu_{\,\,\alpha}\p_\mu\p_\nu q_t||_0^2 \leq ||R'q_{ttt}||_0^2+ ||\mathcal{F}_2 J^{-1}||_0^2 \leq  W_3^2+\epsilon(\mathcal{N}+\mathcal{N}^2)+\PP_0+ \PP\int_0^t\PP + \epsilon||q_t||_2^2. 
\end{align}
\paragraph*{Bound for $||\cp q_t||_{0.5,\Gamma}^2$:} Invoking \eqref{q_t on Gamma} and taking one more tangential derivative, we have
\begin{align}
\cp q_t=\sigma g^{ij}\n^\mu \cp\cp^2_{ij}v_\mu-\sigma g^{ij}g^{kl}(\cp_k v^\tau\n_\tau\cp_l\eta_\mu)\cp\cp_{ij}^2\eta^\mu+\mathcal{R}_{\cp q_t},
\end{align}
where $\mathcal{R}_{\cp q_t}$ consists products of $\cp^k \eta$ and $\cp^k v$, $k=1,2$. To be more specific, $\mathcal{R}_{\cp q_t}$ consists terms of the forms
\begin{align*}
\sigma g^{ij}g^{kl}(\cp_k \cp v^\tau\n_\tau\cp_l\eta_\mu)\cp_{ij}^2\eta^\mu,\q
 \sigma g^{ij}g^{kl}(\cp_k \cp \eta^\tau\n_\tau\cp_l\eta_\mu)\cp_{ij}^2 v^\mu,\\
\sigma (\cp^i v^\mu)( \cp^j \eta_\mu) g^{kl}(\cp_k \cp \eta^\tau\n_\tau\cp_l\eta_\mu)\cp_{ij}^2 \eta^\mu,\q
\sigma (\cp \cp^i\eta^\mu)( \cp^j \eta_\mu) g^{kl}(\cp_k \cp \eta^\tau\n_\tau\cp_l\eta_\mu)\cp_{ij}^2 v^\mu.
\end{align*}
Given these, we have:
\begin{align}
||\cp q_t||_{0.5,\Gamma}^2 \lesssim \epsilon(\mathcal{N}^2+\mathcal{N})+ \PP_0+ \PP\int_0^t \PP,
\end{align}
by interpolation and the Young's inequality. Here, $\epsilon \mathcal{N}^2$ appears since
 \begin{equation}
 ||\sqrt{g}g^{ij}\n^\mu \p\cp^3v_\mu||_0^2\lesssim \epsilon ||v||_4^4+||\sqrt{g}g^{ij}\n||_2^4\lesssim \epsilon \mathcal{N}^2+\PP_0+\PP\int_0^t\PP,
 \label{quadratic N}
 \end{equation}
and we remark here that the interpolation cannot be applied since $\p\cp^3 v$ is of the top order.  
Wrapping these up and invoking Lemma \ref{esi lower order lemma} and Corollary \ref{esi lowest order q and q_t lemma}, we get
\begin{align}
||q_t||_2^2 \lesssim W_3^2+\PP_0+\PP\int_0^t\PP + \epsilon||q_t||_2^2+\epsilon(\mathcal{N}^2+\mathcal{N}),
\end{align}
which proves the estimate for $||q_t||_2^2$ by invoking \eqref{est W} and then absorbing $\epsilon||q_t||_2^2$ to the left hand side. 
\end{proof}

\rmk We are unable to control $||\p q_{tt}||_1$ when surface tension is present. This is due to that the div-curl estimate yields the boundary term $||\cp q_{tt}||_{0.5,\Gamma}$, where $\cp q_{tt}\sim \cp^3 v_t$ on $\Gamma$, and hence $||\cp q_{tt}||_{0.5,\Gamma}$ yields a loss of derivative. Therefore, one has to define the energy using the $\rr$-weighted derivatives and so the corresponding term can then be controlled by the energy. 

\subsection{Bounds for $\int_0^t\mathcal{I}_{1,2,3,4}$}
\label{section 3.3}
This section is devoted to control $\int_0^t\mathcal{I}_{1,2,3,4}$. We recall
\begin{align*}
\mathcal{I}_1=\int_{\Omega}(\dd^r v_\alpha)\Big([\dd^r, RJ](a^{\mu\alpha}\frac{\p_\mu q}{R})\Big) ,\q
\mathcal{I}_2=\int_{\Omega}(\dd^r \p_\mu v_\alpha)\Big([\dd^r, A^{\mu\alpha}] q\Big),\\
\mathcal{I}_3= \int_{\Omega}\big([\dd^r ,A^{\mu\alpha}] \p_\mu v_\alpha\big)\dd^r q,\q
\mathcal{I}_4=\int_{\Omega}\Big([\dd^r, JR'R^{-1}] \p_t q\Big)\dd^r q.
\end{align*}
\nota In what follows, we use $D$ to denote either $\cp$ or $\p_t$. This allows us to represent $\dd^r$ as $(\rr)^\ell D^r$, where $r=\frac{1}{2},1,\frac{3}{2}, 2$.

\subsubsection{Control of $\int_0^t \mathcal{I}_1$}
\paragraph*{For non-$\rr$-weighted $\dd^r$:} We recall that there are four mixed derivatives which are not $\rr$-weighted, which are $\cp$, $\cp^2$,$\cp \p_t$ and $\cp^2\p_t$. Hence, it suffices to consider only the case when $\dd^r=\cp^2 \p_t$. Invoking \eqref{RJ=rho_0} and Theorem \ref{esi ||q||_2 and ||q_t||_2 theo}, We have:
\begin{align}
\mathcal{I}_1={\sum}_{\substack{j_1+j_2=2\\j_1\geq 1}}\int_{\Omega}(\cp^2 \p_t v_\alpha)(\cp^{j_1}\rho_0)(\cp^{j_2}\p_t(a^{\mu\alpha}\frac{\p_\mu q}{R})).
\end{align}
Since, to the highest order, the last term on the right hand side is $R^{-1}a^{\mu\alpha}\cp\p q_t$, which can be controlled by invoking Theorem \ref{esi ||q||_2 and ||q_t||_2 theo}. Therefore, 
\begin{equation}
\int_0^t \mathcal{I}_1 \leq  \PP_0+\PP\int_0^t\PP.
\end{equation}
The $\epsilon P(N)$ term introduced in Theorem \ref{esi ||q||_2 and ||q_t||_2 theo} does not 
figure here since $\mathcal{I}_1$ is estimated under the time integral.

\paragraph*{For $\rr$-weighted $\dd^r$:} It suffices to consider derivatives of the form $(\rr)^{\ell}D^{r-2}\cp\p_t$, where $\ell =\frac{1}{2}, 1, \frac{3}{2}$ and $r\leq 4$, since otherwise $\mathcal{I}_1$ would be $0$ due to \eqref{RJ=rho_0}. 
\begin{align}
\mathcal{I}_1 = {\sum}_{j_1+j_2=r-2}\int_\Omega (\rr)^{2\ell} (D^{r-2}\cp\p_t v_\alpha)(\cp D^{j_1}\rho_0)( D^{j_2}\p_t(a^{\mu\alpha}\frac{\p_\mu q}{R})).
\end{align}
We henceforth adopt:
\nota \label{notation lleq} We use $\lleq$ to denote equality modulo lower order terms that can be controlled, i.e., $A\lleq B$ mean $A= B +\text{error terms}$, where the ``error terms" can be controlled by the bound of $B$ plus $\PP_0+\PP\int_0^t\PP$. 

Invoking \eqref{equivlent R_kk} and \eqref{R_kk assumption} at $t=0$ lead to
\begin{align}
\int_0^t\mathcal{I}_1 \lleq {\sum}_{j_1+j_2=r-2}\int_0^t\int_\Omega (\mathfrak{R}_\kk)^{2\ell+1} (D^{r-2}\cp\p_t v_\alpha)(\cp D^{j_1} p_0)( D^{j_2}\p_t(a^{\mu\alpha}\frac{\p_\mu q}{R}))\no\\
=\int_0^t\int_\Omega  (\cp D^{j_1} p_0)\Big((\rr)^{\ell+\frac{1}{2}}D^{r-2}\cp\p_t v_\alpha\Big)\Big((\rr)^{\ell+\frac{1}{2}} D^{j_2}\p_t(a^{\mu\alpha}\frac{\p_\mu q}{R})\Big)
\leq \PP_0+\PP\int_0^t\PP.
\end{align}
\rmk The above expression yields a slightly better bound for $D^{r-2}\cp\p_t v$, since $\PP$ requires only $||(R')^{\ell}D^{r-2}\cp\p_t v||_0$.
\subsubsection{Control of $\int_0^t\mathcal{I}_2$}
For each $r$, $\int_0^t\mathcal{I}_2$ contains a term which is of the order $r+1$, i.e., 
\begin{align*}
\int_0^t\mathfrak{T}=\int_\Omega (\dd^r\p_\mu v_\alpha)(\dd^r A^{\mu\alpha}) q.
\end{align*}
There is no problem to control $\mathfrak{T}$ when $r\leq 2$, and when $r=3$, we need to put extra effort to control $\mathfrak{T}$ when $\dd^3 = \cp^2 \p_t$ since there are terms which cannot be controlled directly without $\rr$-weights, and one needs to integrate by parts in (tangential) spatial derivative and time derivative, respectively.  
On the other hand, when $r=4$, this term is of above the top order, but it can be controlled using 
one of the special cancellations referred to in section
\ref{sec_special_cancellations}, as we now show.
\paragraph*{For non-$\rr$-weighted $\dd^r$:}As mentioned above, we consider only the case when $r=3$ and $\dd^3=\cp^2\p_t$. In this case, 
\begin{align}
\mathfrak{T} = \int_\Omega (\cp^2\p_\mu \p_t v_\alpha)(\cp^2\p_t A^{\mu\alpha}) q.
\end{align}
Although this term is of the correct order, $\cp^2 \p_\mu \p_t v$ cannot be controlled without $\rr$-weights. Hence, we integrate by parts with respect to the tangential derivative and get: 
\begin{align}
\mathfrak{T} = -\Big(\int_\Omega (\cp \p_\mu \p_t v_\alpha)(\cp^3\p_t A^{\mu\alpha}) q+\int_\Omega (\cp \p_\mu \p_t v_\alpha)(\cp^2\p_t A^{\mu\alpha}) \cp q\Big)\no\\
\leq ||\cp \p v_t||_0||\cp^3 \p_t A||_0||q||_{L^{\infty}}+||\cp \p v_t||_0||\cp^2\p_t A||_{L^4}||\cp q||_{L^4}.
\label{T lower order}
\end{align}
Here, one adapts Theorem \ref{esi ||q||_2 and ||q_t||_2 theo} to control $||q||_2$.
Integrating with respect to time, we obtain:
\begin{align}
\int_0^t \mathfrak{T} \lesssim \PP_0+\PP\int_0^t\PP.
\end{align}

We next consider $\mathcal{I}_2-\mathfrak{T}$.  All terms involved in $\mathcal{I}_2-\mathfrak{T}$ can be controlled straightforwardly after integrating by part with respect to $\cp$ thanks to Theorem \ref{esi ||q||_2 and ||q_t||_2 theo}, except for 
$$
\int_\Omega (\cp^2 \p_\mu \p_t v_\alpha)(\p_t A^{\mu\alpha})(\cp^2 q).
$$
This is due to that integrating by part in $\cp$ yields $\cp^3 q$ which cannot be controlled without $\rr$-weights. To deal with this issue, we consider 
\begin{equation}
\int_0^t \int_\Omega (\cp^2 \p_\mu \p_t v_\alpha)(\p_t A^{\mu\alpha})(\cp^2 q).\label{low 1}
\end{equation}
Integrating by part in time, we get:
\begin{align}
\int_\Omega (\cp^2 \p_\mu v_\alpha)(\p_t A^{\mu\alpha})(\cp^2 q)|^{t}_{0}-\int_0^t \int_\Omega (\cp^2 \p_\mu  v_\alpha)(\p_t^2 A^{\mu\alpha})(\cp^2 q)-\int_0^t \int_\Omega (\cp^2 \p_\mu v_\alpha)(\p_t A^{\mu\alpha})(\cp^2 \p_t q).
\end{align}
The last two terms are bounded by $\PP_0+\PP\int_0^t \PP$ thanks to Theorem \ref{esi ||q||_2 and ||q_t||_2 theo}, while the pointwise term at $t=0$ by $\PP_0$. The pointwise term at $t$ is bounded by
\begin{align}
||\eta||_3||\cp^2 q||_0||\cp^2 v||_1^{\frac{1}{2}}||\cp^2 v||_2^{\frac{1}{2}}||v||_1^{\frac{1}{2}}||v||_2^{\frac{1}{2}}
\lesssim \epsilon P(\mathcal{N})+\PP_0+\PP\int_0^t\PP,\label{low 2}
\end{align}
which is controlled by the right hand side of \eqref{estimate E}.

\paragraph*{For $\rr$-weighted $\dd^r$:}$\mathcal{I}_2$ contains a term above the top order, i.e., 
\begin{align*}
\mathfrak{T}=\int_\Omega (\dd^r\p_\mu v_\alpha)(\dd^r A^{\mu\alpha}) q.
\end{align*}
This term is controlled using the aforementioned special cancellation
(see section \ref{sec_special_cancellations}).
For $\rr$-weighted derivatives, it suffices to consider only the case when $r=4$, i.e., the derivatives are of the form $(\rr)^{\ell}D^3 \p_t$, for $\ell=1,\frac{3}{2}, 2$. Then the ``tricky" term 
to be bounded is:
\begin{align}
\mathcal{L}=\int_0^t\mathfrak{T}=\int_0^t \int_\Omega (\rr)^{2\ell}(D^3 \p_t\p_\mu v_\alpha)(D^3\p_t A^{\mu\alpha}) q. \label{tricky}
\end{align}
In view of \eqref{A}, expanding the index $\mu$ in \eqref{tricky}, we have
\begin{align}
\mathcal{L} = \underbrace{\int_0^t \int_\Omega (\rr)^{2\ell}q\epsilon^{\alpha\lambda\tau}\cp_2D^3v_\lambda\p_3 \eta_\tau\cp_1D^3\p_tv_\alpha}_{L_1}+\underbrace{\int_0^t\int_\Omega(\rr)^{2\ell}q\epsilon^{\alpha\lambda\tau}\cp_2\eta_\lambda\p_3D^3v_\tau\cp_1D^3\p_tv_\alpha}_{L_2}\no\\
\underbrace{-\int_0^t \int_\Omega (\rr)^{2\ell}q\epsilon^{\alpha\lambda\tau}\cp_1D^3v_\lambda\p_3 \eta_\tau\cp_2D^3\p_tv_\alpha}_{L_3}\underbrace{-\int_0^t\int_\Omega(\rr)^{2\ell}q\epsilon^{\alpha\lambda\tau}\cp_1\eta_\lambda\p_3D^3v_\tau\cp_2D^3\p_tv_\alpha}_{L_4}\no\\
\underbrace{\int_0^t \int_\Omega (\rr)^{2\ell}q\epsilon^{\alpha\lambda\tau}\cp_1D^3v_\lambda\cp_2 \eta_\tau\p_3D^3\p_tv_\alpha}_{L_5}+\underbrace{\int_0^t\int_\Omega(\rr)^{2\ell}q\epsilon^{\alpha\lambda\tau}\cp_1\eta_\lambda\cp_2D^3v_\tau\p_3D^3\p_tv_\alpha}_{L_6}+L_{low}, \label{cancellation}
\end{align}
where $L_{low}$ are lower order terms, which are all of the form
\begin{align}
{\sum}_{\substack{j_1+j_2=3\\j_1,j_2\leq 2}}\int_0^t \int_\Omega (\rr)^{2\ell}q (\p D^{j_1}v)(\p D^{j_2}\eta)(\p D^3 \p_t v)={\sum}_{\substack{j_1+j_2=3\\j_1,j_2\leq 2}}\Big(\int_\Omega (\rr)^{2\ell} q (\p D^{j_1}v)(\p D^{j_2}\eta)(\p D^3 v)\no\\
-\int_0^t\int_\Omega (\rr)^{2\ell} q_t (\p D^{j_1}v)(\p D^{j_2}\eta)(\p D^3 v)-\int_0^t\int_\Omega (\rr)^{2\ell} q (\p D^{j_1}v_t)(\p D^{j_2}\eta)(\p D^3 v)\no\\
\int_0^t\int_\Omega (\rr)^{2\ell} q (\p D^{j_1}v)(\p D^{j_2}v)(\p D^3 v)  \Big)\no.
\end{align}
Invoking Theorem \ref{esi ||q||_2 and ||q_t||_2 theo}, it is easy to see that the last three terms are controlled by $\PP_0+\PP\int_0^t\PP$, while the pointwise term at $t$ is treated similar to \eqref{ptwise I}-\eqref{ptwise II}, after distributing correct amount of $\rr$-weights to each term. We omit the detail here. But it is worth noting that there are more than enough $\rr$-weights for the pointwise term since there is one time derivative less. 

Next, integrating by part in time in $L_3$, we find
\begin{align}
-\int_0^t \int_\Omega (\rr)^{2\ell}q\epsilon^{\alpha\lambda\tau}\cp_1D^3v_\lambda\p_3 \eta_\tau\cp_2D^3\p_tv_\alpha \lleq -\int_\Omega(\rr)^{2\ell}q\epsilon^{\alpha\lambda\tau}\cp_1D^3v_\lambda\p_3 \eta_\tau\cp_2D^3v_\alpha\no\\
\int_0^t \int_\Omega (\rr)^{2\ell}q\epsilon^{\alpha\lambda\tau}\cp_1D^3\p_tv_\lambda\p_3 \eta_\tau\cp_2D^3v_\alpha.
\end{align}
Adding $L_1$, we get:
\begin{align}
L_1+L_3\lleq \int_0^t \int_\Omega (\rr)^{2\ell}q\epsilon^{\alpha\lambda\tau}\cp_2D^3v_\lambda\p_3 \eta_\tau\cp_1D^3\p_tv_\alpha\no\\+\int_0^t \int_\Omega (\rr)^{2\ell}q\epsilon^{\alpha\lambda\tau}\cp_1D^3\p_tv_\lambda\p_3 \eta_\tau\cp_2D^3v_\alpha\no\\-\int_\Omega(\rr)^{2\ell}q\epsilon^{\alpha\lambda\tau}\cp_1D^3v_\lambda\p_3 \eta_\tau\cp_2D^3v_\alpha|^{t}_{0}\no\\
=-\int_\Omega(\rr)^{2\ell}q\epsilon^{\alpha\lambda\tau}\cp_1D^3v_\lambda\p_3 \eta_\tau\cp_2D^3v_\alpha|^{t}_{0}=L_{13},
\end{align}
since first and the second term cancels with each other by the antisymmetry of $\epsilon^{\alpha\lambda\tau}$. Similarly, we have
\begin{align}
L_4+L_6\lleq L_{46}=\int_\Omega(\rr)^{2\ell}q\epsilon^{\alpha\lambda\tau}\cp_1\eta_\lambda\cp_2D^3v_\tau\p_3D^3v_\alpha|^{t}_{0},\\
L_2+L_5\lleq L_{25}= \int_\Omega (\rr)^{2\ell}q\epsilon^{\alpha\lambda\tau}\cp_1D^3v_\lambda\cp_2 \eta_\tau\p_3D^3v_\alpha|^{t}_{0}.
\end{align}
\paragraph*{Bounds for $L_{13}$, $L_{46}$ and $L_{25}$}
Since $L_{13}$ is pointwise in $t$, it suffices to consider 
$$
\int_\Omega(\rr)^{2\ell}q\epsilon^{\alpha\lambda\tau}\cp_1D^3v_\lambda\p_3 \eta_\tau\cp_2D^3v_\alpha|_t
$$
only, since the other part is controlled directly by $\PP_0$. In addition, since $D^3$ corresponds to $\p_t^3$, $\cp \p_t^2$, $\cp^2\p_t$ and $\cp^3$, associated with weights $(\rr)^2$, $(\rr)^{\frac{3}{2}}$, $\rr$ and $\rr$, respectively, we have:
\begin{align}
\int_\Omega(\rr)^{4}q\epsilon^{\alpha\lambda\tau}\cp_1\p_t^3v_\lambda\p_3 \eta_\tau\cp_2\p_t^3v_\alpha|_t \leq \int_\Omega(\rr) q\epsilon^{\alpha\lambda\tau}((\rr)^{\frac{3}{2}}\cp_1\p_t^3v_\lambda)\p_3 \eta_\tau((\rr)^{\frac{3}{2}}\cp_2\p_t^3v_\alpha)|_t\no\\
\leq ||R||_2||\eta||_3 ||(\rr)^{\frac{3}{2}}v_{ttt}||_1^2\leq \epsilon P(\mathcal{N})+\PP_0+\PP\int_0^t\PP,
\end{align} 
where we have used $||\rr q||_2\lesssim ||R' q||_2=||R||_2$.
Similarly, we have
\begin{align}
\int_\Omega(\rr)^{3}q\epsilon^{\alpha\lambda\tau}\cp_1\cp\p_t^2v_\lambda\p_3 \eta_\tau\cp_2\cp\p_t^2v_\alpha|_t+\int_\Omega(\rr)^{2}q\epsilon^{\alpha\lambda\tau}\cp_1\cp^2\p_t v_\lambda\p_3 \eta_\tau\cp_2\cp^2\p_t v_\alpha|_t\no\\
+\int_\Omega(\rr)^{2}q\epsilon^{\alpha\lambda\tau}\cp_1\cp^3v_\lambda\p_3 \eta_\tau\cp_2\cp^3v_\alpha|_t \leq \epsilon P(\mathcal{N})+\PP_0+\PP\int_0^t\PP. 
\end{align}
Moreover, this method can be adapted to control $L_{46}$ and $L_{25}$, and we omit the details. 
Therefore, 
\begin{align}
(L_1+L_3)+(L_4+L_6)+(L_2+L_5) \leq \epsilon P(\mathcal{N})+\PP_0+\PP\int_0^t\PP.
\end{align}

Now, we complete the treatment of $\mathcal{I}_2$ by estimating the rest of the terms, i.e., $\mathcal{I}_2-\mathfrak{T}$, for $\rr$-weighted forth order derivatives.  
Expressing: 
\begin{align}
\mathcal{I}_2-\mathfrak{T}= \int_\Omega (\rr)^{2\ell}(D^3\p_t\p_\mu v_\alpha)\Big(D^3\p_t (A^{\mu\alpha}q)-A^{\mu\alpha}D^3\p_t q-(D^3\p_t A^{\mu\alpha})q\Big) , 
\end{align}
and similar to the non-$\rr$-weighted case, we consider $\int_0^t \mathcal{I}_2-\mathfrak{T}$ and integrate by part in time to get
\begin{align}
\int_0^t \mathcal{I}_2-\mathfrak{T} \lleq \int_\Omega (\rr)^{2\ell}(D^3\p_\mu v_\alpha)\Big(D^3\p_t (A^{\mu\alpha}q)-A^{\mu\alpha}D^3\p_t q-(D^3\p_t A^{\mu\alpha})q\Big)\Big|^{t}_{0}\no\\
-\int_0^t\int_\Omega (\rr)^{2\ell}(D^3\p_\mu v_\alpha)\p_t\Big(D^3\p_t (A^{\mu\alpha}q)-A^{\mu\alpha}D^3\p_t q-(D^3\p_t A^{\mu\alpha})q\Big).
\end{align} 
First, it is easy to check that 
$$
\int_0^t\int_\Omega (\rr)^{2\ell}(D^3\p_\mu v_\alpha)\p_t\Big(D^3\p_t (A^{\mu\alpha}q)-A^{\mu\alpha}D^3\p_t q-(D^3\p_t A^{\mu\alpha})q\Big)\leq \int_0^t \PP,
$$
Second, for the pointwise terms at $t$, it suffices to consider the case when $D^3=\p_t^3$ and $\ell=2$, since the bounds for the other (easier) cases follow from the same method. There are three terms, i.e., 
\begin{align}
\int_\Omega (\rr)^{4}(\p_t^3\p_\mu v_\alpha)(\p_t^3 A)\p_t q,\q \int_\Omega (\rr)^{4}(\p_t^3\p_\mu v_\alpha)(\p_t^2 A)\p_t^2 q,\q \int_\Omega (\rr)^{4}(\p_t^3\p_\mu v_\alpha)(\p_t A)\p_t^3 q. \label{ptwise I}
\end{align}
These terms are treated as
\begin{align}
\int_\Omega (\rr)^{4}(\p_t^3\p_\mu v_\alpha)(\p_t^3 A)\p_t q&\approx\int_\Omega (\rr)^{4}(\p_t^3\p_\mu v_\alpha)\Big((\p v_{tt})(\p \eta)+(\p v_t)(\p v)\Big)\p_t q\no\\
&\lesssim ||(\rr)^{\frac{3}{2}}v_{ttt}||_1||(\rr)^{\frac{3}{2}}v_{tt}||_1^{\frac{1}{2}}||(\rr)^{\frac{3}{2}}v_{tt}||_2^{\frac{1}{2}}||\eta||_3||\rr q_t||_0^{\frac{1}{2}}||\rr q_t||_1^{\frac{1}{2}}\no\\
&+||(\rr)^{\frac{3}{2}}v_{ttt}||_1||\rr v_{t}||_1^{\frac{1}{2}}||\rr v_{t}||_2^{\frac{1}{2}}||(\rr)^{\frac{1}{2}}v||_3||\rr q_t||_0^{\frac{1}{2}}||\rr q_t||_1^{\frac{1}{2}} \no\\
&\lesssim \epsilon P(\mathcal{N})+\PP_0+\PP\int_0^t\PP ,
\end{align}
 and
\begin{align}
\int_\Omega (\rr)^{4}(\p_t^3\p_\mu v_\alpha)(\p_t^2 A)\p_t^2 q= \int_\Omega (\rr)^{4}(\p_t^3\p_\mu v_\alpha)\Big((\p v)^2+(\p v_t)(\p \eta)\Big)\p_t^2 q\no\\
\lesssim ||(\rr)^{\frac{3}{2}}v_{ttt}||_1\Big(||v||_2||\rr v||_3+||\rr v_t||_1^{\frac{1}{2}}||\rr v_t||_2^{\frac{1}{2}}||\eta||_3\Big)||(\rr)^{\frac{3}{2}}q_{tt}||_0^{\frac{1}{2}}||(\rr)^{\frac{3}{2}}q_{tt}||_1^{\frac{1}{2}}\no\\
\lesssim \epsilon P(\mathcal{N})+\PP_0+\PP\int_0^t\PP .
\end{align}
Finally, we have
\begin{align}
\int_\Omega (\rr)^{4}(\p_t^3\p_\mu v_\alpha)(\p_t A)\p_t^3 q=\int_\Omega (\rr)^{4}(\p_t^3\p_\mu v_\alpha)(\p v)(\p \eta)\p_t^3 q\no\\
\lesssim ||(\rr)^{\frac{3}{2}}v_{ttt}||_1||v||_1^{\frac{1}{2}}||v||_2^{\frac{1}{2}}||\eta||_3||(\rr)^{\frac{5}{2}}q_{ttt}||_0^{\frac{1}{2}}||(\rr)^{\frac{5}{2}}q_{ttt}||_1^{\frac{1}{2}}\no\\
\lesssim \epsilon P(\mathcal{N})+\PP_0+\PP\int_0^t\PP. \label{ptwise II}
\end{align}

\subsubsection{Control of $\int_0^t\mathcal{I}_3$}
\label{section 3.3.3}
\paragraph*{For non-$\rr$-weighted $\dd^r$:}
Expressing these derivatives as $D^r$ where $r\leq 3$, we have:
\begin{align}
\mathcal{I}_3 = {\sum}_{\substack{j_1+j_2=r\\j_1\geq 1}}\int_{\Omega}(D^{j_1}A^{\mu\alpha}) (\p_\mu D^{j_2}v_\alpha)(D^r q)
\leq {\sum}_{\substack{j_1+j_2=r\\j_1\geq 1}}||(D^{j_1}A^{\mu\alpha}) (\p_\mu D^{j_2}v_\alpha)||_0||D^r q||_0,
\end{align}
and so $\int_0^t \mathcal{I}_3\leq \PP_0+\PP\int_0^t\PP$ in light of Theorem \ref{esi ||q||_2 and ||q_t||_2 theo}. 

\paragraph*{For $\rr$-weighted $\dd^r$:} It suffices to consider only the case when $r=4$, i.e., the derivatives are of the form $(\rr)^{\ell}D^3 \p_t$, for $\ell=1,\frac{3}{2},2$. Now, 
\begin{align}
\mathcal{I}_3 = \int_{\Omega}(\rr)^{2\ell}(\p_t A^{\mu\alpha})(D^3 \p_\mu v_\alpha)(D^3 \p_t q)+\int_{\Omega}(\rr)^{2\ell}(D A^{\mu\alpha})(D^2 \p_t \p_\mu v_\alpha)(D^3 \p_t q)
\no\\+\int_{\Omega}(\rr)^{2\ell}(\p_t D^3 A^{\mu\alpha})(\p_\mu v_\alpha)(D^3 \p_t q)+ \text{error terms},
\end{align}
where the main term is equal to
\begin{align}
\int_{\Omega}(\p_t A^{\mu\alpha})\Big((\rr)^{\ell-\frac{1}{2}}D^3 \p_\mu v_\alpha\Big)\Big((\rr)^{\ell+\frac{1}{2}}D^3 \p_t q\Big)\no\\
+\int_{\Omega}(\cp A^{\mu\alpha})\Big((\rr)^{\ell}D^2\p_t \p_\mu v_\alpha\Big)\Big((\rr)^{\ell}D^3 \p_t q\Big)\no\\
+ \int_{\Omega}(\p_\mu v_\alpha)\Big((\rr)^{\ell-\frac{1}{2}}\p_t D^3 A^{\mu\alpha}\Big)\Big((\rr)^{\ell+\frac{1}{2}}D^3 \p_t q\Big)=\mathcal{I}_{3,1}+\mathcal{I}_{3,2}+\mathcal{I}_{3,3},
\end{align}
where $\mathcal{I}_{3,2}$ does not appear when $\dd^4= \rr^2 \p_t^4$. 

$\int_0^t\mathcal{I}_{3,1}+\mathcal{I}_{3,3}$ can be controlled directly by $\PP_0+\PP\int_0^t\PP$. For $\mathcal{I}_{3,2}$, one requires the wave energy \eqref{W4} to control $||(\rr)^{\ell}D^3 \p_t q||_0$ when $D^3$ contains at least one\footnote{This is explained in the remark after Theorem \ref{est W4 thm}.} $\p_t$, and \eqref{q and R} to control this term when $D^3=\cp^3$ (i.e., $\rr \cp^3 \p_t q\sim \cp^3 \p_t R$), and so $\int_0^t\mathcal{I}_{3,2}$ can be controlled appropriately by  $\PP_0+\PP\int_0^t\PP$.
Furthermore, the (time integrated) error terms are of the form
\begin{align}
{\sum}_{\substack{j_1+j_2+j_3=3\\j_1+j_2\geq 1}}\int_0^t\int_\Omega (\rr)^{2\ell}(\p D^{j_1} \eta)(\p D^{j_2}v)(\p D^{j_3} v)(D^4 q)\no\\
 ={\sum}_{\substack{j_1+j_2+j_3=3\\j_1+j_2\geq 1}}\int_0^t\int_\Omega \Big((\rr)^{\ell-\frac{1}{2}}(\p D^{j_1} \eta)(\p D^{j_2} v)(\p D^{j_3} v)\Big)\Big((\rr)^{\ell+\frac{1}{2}}D^4 q\Big)\leq \PP_0+\PP\int_0^t\PP.
\end{align}

\subsubsection{Control of $\int_0^t\mathcal{I}_4$}
$\mathcal{I}_4$ is the easiest one to control among the other $\mathcal{I}$ terms. This is due to the assumption \eqref{R_kk assumption}, which implies that there are ``sufficient" $\rr$-weights that can be distributed for all terms. In addition to this, we can also use the fact $DR = R'Dq$ to get an extra $\rr$-weights if necessary. 
\paragraph*{For non-$\rr$-weighted $\dd^r$:} 
By \eqref{R_kk assumption} and since $r\leq 3$,  invoking Theorem \ref{esi ||q||_2 and ||q_t||_2 theo}, we have:
\begin{align}
\int_0^t\mathcal{I}_4 \lleq {\sum}_{\substack{j_1+j_2=r\\j_1\geq 1}}\int_0^t\int_{\Omega}\rr\Big(D^{j_1}(\rho_0 R^{-2})\Big)(D^{j_2}\p_t q)(D^r q)\leq \PP_0+\PP\int_0^t\PP.
\end{align}
\paragraph*{For $\rr$-weighted $\dd^r$:} For $\ell=\frac{1}{2},1,\frac{3}{2},2$, we have:
\begin{align}
\int_0^t\mathcal{I}_4 \lleq \sum_{\substack{j_1+j_2=r\\j_1\geq 1}}\int_0^t\int_{\Omega}(\rr)^{2\ell+1}\Big(D^{j_1}(\rho_0R^{-2})\Big)(D^{j_2}\p_t q)(D^r q)=\no\\
{\sum}_{\substack{j_1+j_2=r\\j_1\geq 1}}\int_0^t\int_{\Omega}\Big((\rr)^{\ell+\frac{1}{2}}\Big(D^{j_1}(\rho_0R^{-2})\Big)(D^{j_2}\p_t q)\Big)\Big((\rr)^{\ell+\frac{1}{2}}D^r q\Big)\leq \PP_0+\PP\int_0^t\PP,
\end{align}
where the fact $DR = R' Dq$ is used if $j_1=1$.

\subsection{Control of $\int_0^t\mathcal{B}$ for non-$
\rr$-weighted $\dd^r$\label{section B}}
This section is devoted to control the boundary terms 
\begin{align}
\mathcal{B}_1= \sigma\int_{\Gamma}(\dd^r v_\alpha) \big([\dd^r, \sqrt{g}g^{ij}\Pi^\alpha_\mu]\cp^2_{ij}\eta^\mu\big)\,dS,\q \mathcal{B}_2 =
{-} \sigma\int_\Gamma \cp_i (\sqrt{g}g^{ij}\Pi_\mu^\alpha)(\p_t \dd^r \eta_\alpha)(\cp_j \dd^r \eta^\mu)\,dS,\no\\
\mathcal{B}_3= {\frac{1}{2}}\sigma\int_\Gamma \p_t(\sqrt{g}g^{ij}\Pi^\alpha_\mu)(\cp_i \DD^r\eta_\alpha)(\cp_j \DD^r \eta^\mu)\,dS,
\nonumber
\end{align} 
which appears in the energy estimate when $\dd^r$ is non-$\rr$-weighted.
The $\rr$-weighted cases are treated in section \ref{sec_B_weighted}.

We recall that if $\dd^r$ is non-$\rr$-weighted, then $r\leq 3$, i.e., the corresponding term is of lower order. Because of this, it would be suffice to consider the case when $\dd^r=\cp^2\p_t$. Now, since $\Pi_\mu^\alpha = \n_\mu\n^\alpha$, we have:
\begin{align*}
\mathcal{B}_1=\sigma{\sum}_{\substack{j_1+j_2=3\\j_1\geq 1}}\int_{\Gamma}(\cp^2\p_t v_\alpha)D^{j_1}(\sqrt{g}g^{ij}\n_\mu\n^\alpha)(D^{j_2}\cp^2_{ij}\eta^\mu)\,dS,\\
\mathcal{B}_2 = \sigma\int_\Gamma \cp_i (\sqrt{g}g^{ij}\Pi_\mu^\alpha)(\p_t \cp^2 v_\alpha)(\cp_j \cp^2 v^\mu)\,dS,\q \mathcal{B}_3=\frac{1}{2}\sigma\int_\Gamma \p_t(\sqrt{g}g^{ij}\Pi^\alpha_\mu)(\cp_i \cp^2 v_\alpha)(\cp_j \cp^2 v^\mu)\,dS
\end{align*}
Invoking  Lemma \ref{prelim lemma b}, we get
\begin{align*}
\p_t(\sqrt{g}g^{ij}\n_\mu\n^\alpha)=Q(\cp \eta)\cp v,\q
\cp (\sqrt{g}g^{ij}\n_\mu\n^\alpha) = Q(\cp \eta)\cp^2 \eta,
\end{align*}
where $Q$ is a rational function, and hence
\begin{align*}
\cp \p_t (\sqrt{g}g^{ij}\n_\mu\n^\alpha) = Q(\cp \eta, \cp v)\cp^2\eta+Q(\cp \eta, \cp v)\cp^2 v,\q
\cp^2 (\sqrt{g}g^{ij}\n_\mu\n^\alpha) = Q(\cp \eta, \cp^2\eta)\cp^3\eta,\\
\cp^2\p_t (\sqrt{g}g^{ij}\n_\mu\n^\alpha) =Q(\cp \eta, \cp v, \cp^2\eta, \cp^2 v) (\cp^3\eta+ \cp^3 v).
\end{align*}
In light of these, we have
\begin{align*}
\int_0^t\mathcal{B}_2 = \sigma \int_0^t\int_\Gamma Q(\cp \eta)\cp^2 \eta (\p^2 \p_t v)(\p^3 v)\,dS\leq \PP_0+\PP\int_0^t\PP,
\end{align*}
via $(H^{-\frac{1}{2}}, H^{\frac{1}{2}})$ duality. Moreover, $\int_0^t\mathcal{B}_3$ is treated similarly.  On the other hand, 
\begin{align*}
\mathcal{B}_1\lleq \sigma \int_\Gamma (\cp^2\p_t v)Q(\cp\eta)(\cp v)(\cp^4 \eta)+ \sigma \int_\Gamma (\cp^2\p_t v)Q(\cp \eta, \cp v, \cp^2\eta, \cp^2 v)(\cp^3\eta+\cp^3 v)(\cp^2 \eta)\no\\
+\sigma \int_\Gamma (\cp^2\p_t v)Q(\cp \eta, \cp^2\eta)(\cp^3\eta)(\cp^2 v)+\sigma \int_\Gamma (\cp^2\p_t v)Q(\cp \eta, \cp v)(\cp^2 v+\cp^2 \eta)(\cp^3 \eta).
\end{align*}
The last three terms can be controlled in a routine fashion.
However,  $\sigma \int_\Gamma (\cp^2\p_t v)Q(\cp\eta)(\cp v)(\cp^4 \eta)$ cannot be controlled directly since $(H^{-\frac{1}{2}},H^{\frac{1}{2}})$ duality requires the control $||v_t||_3$ which is not part of $\PP$, and so we consider
$$
\sigma\int_0^t \int_\Gamma (\cp^2\p_t v)Q(\cp\eta)(\cp v)(\cp^4 \eta)
$$
and then integrate by parts in $t$. This yields
\begin{align}
\sigma\int_0^t \int_\Gamma (\cp^2\p_t v)Q(\cp\eta)(\cp v)(\cp^4 \eta)=\sigma\int_\Gamma (\cp^2 v)Q(\cp\eta)(\cp v)(\cp^4 \eta)|_0^t\no\\
-\sigma\int_0^t \int_\Gamma (\cp^2 v)Q(\cp\eta,\cp v)(\cp^4 \eta)-\sigma\int_0^t \int_\Gamma (\cp^2 v)Q(\cp\eta)(\cp v_t)(\cp^4 \eta)\no\\
-\sigma\int_0^t \int_\Gamma (\cp^2 v)Q(\cp\eta)(\cp v)(\cp^4 v).
\end{align} 
The last three term on the right hand side can be controlled directly by $\PP_0+\PP\int_0^t\PP$ via $(H^{-\frac{1}{2}},H^{\frac{1}{2}})$ duality. Moreover, the pointwise term is bounded by 
\begin{align}
 &\PP_0+\sigma Q(||\eta||_4) ||(\cp^2 v)(\cp v)||_1\no\\
 &\lesssim\PP_0+\sigma Q(||\eta||_4)(||v||_3^{\frac{1}{2}}||v||_4^{\frac{1}{2}}||v||_1^{\frac{1}{2}}||v||_2^{\frac{1}{2}}+||v||_2||v||_3)\no\\
&\lesssim \epsilon\mathcal{N}+ \PP_0+\PP\int_0^t\PP.
\end{align}

\subsection{Control of $\mathcal{B}$ for $\rr$-weighted $\dd^r$\label{sec_B_weighted}}
Here we show how to control $\mathcal{B}$ when
$\dd^2 =  \sqrt{\rr}\p_t^2$,
$\dd^3 =   \sqrt{\rr}(\cp\p_t^2), $ $\dd^3 = \rr\p_t^3$,
$\dd^4 = \rr(\cp^3\p_t)$, $\dd^4 = \rr(\cp^2\p_t^2) $, $\dd^4 = (\rr)^{\frac{3}{2}}(\cp\p_t^3)$
and $\dd^4 = (\rr)^2\p_t^4$.

\subsubsection{Case $\dd^4 = (\rr)^2\p_t^4$\label{sec_weighted_main_case}}

The ensuing calculations produce a series of terms. In what follows we focus on the 
most delicate ones, in particular those leading to special cancellations.
The remaining  terms will either be of lower order 
or can be controlled by arguments similar to the ones presented for the aforementioned
main terms. Therefore, all such remainders are collected and estimated at the very end
in section \ref{sec_remainders}. We note that
certain cancellations are only visible after a series of manipulations have been made, requiring
us to keep track of the explicit form of most terms in our calculations.

The following remark will be used throughout below.
In view of identity Lemma \ref{prelim lemma b}--\ref{n_Pi_contraction},
we have
$\hat{n}^\alpha \cp^m \partial_t^k v_\alpha = \hat{n}_\tau \Pi^{\tau \alpha} 
\cp^m \partial_t^k v_\alpha$, so that an estimate for
$\hat{n} \cdot \cp^m \partial_t^k v$ can controlled by
 $ \Pi \cp^m \partial_t^k v$.

We shall also need the following identity
\begin{align}
\partial_t v^\alpha \cp_l \eta_\alpha = -\frac{J}{\rho_0} \cp_l q,\,
\text{ on } \Gamma,
\label{identity_partial_t_v_q} 
\end{align}
which is obtained upon contracting the first equation in \eqref{E}
with $\cp_l \eta_\alpha$, using the definition of $a$, and
\eqref{RJ=rho_0}.

\paragraph{Estimate for $\int_0^t \mathcal{B}_3$ with $\dd^4 = (\rr)^2\p_t^4$\label{sec_B3}}

Using  $\dd^4 = \bw^2 \partial_t^4$ in $\mathcal{B}_3$ gives

\begin{align}
\mathcal{B}_3 &= \frac{1}{2}\sigma\int_\Gamma \p_t(\sqrt{g}g^{ij}\Pi^\alpha_\mu)
\cp_i (\bw^2 \partial_t^3 v_\alpha) \cp_j (\bw^2 \partial_t^3 v^\mu)\,dS
\nonumber
\\
& =
\frac{1}{2}\sigma\int_\Gamma \p_t(\sqrt{g}g^{ij}) \bw^4\Pi^\alpha_\mu
\cp_i  \partial_t^3 v_\alpha\cp_j  \partial_t^3 v^\mu\,dS
+ 
\sigma\int_\Gamma \sqrt{g}g^{ij} \bw^4
\partial_t \Pi^\alpha_\lambda \Pi^\lambda_\mu
\cp_i  \partial_t^3 v_\alpha\cp_j  \partial_t^3 v^\mu\,dS
\nonumber
\\
& = \mathcal{B}_{31} + \mathcal{B}_{32},
\nonumber
\end{align} 
where we used Lemma \ref{prelim lemma b}--\ref{Pi_Pi_contraction}.
It is immediate
to estimate
\begin{align}
\norm{\mathcal{B}_{31}} \leq \poly \norm{ \Pi \bw^2 \cp  \partial_t^3 v}_{0,\bou}.
\nonumber
\end{align}
For $\mathcal{B}_{32}$, use 
 Lemma \ref{prelim lemma b}--\ref{Pi_nn}
 to find
\begin{align}
\mathcal{B}_{32} & = 
\sigma\int_\Gamma \sqrt{g}g^{ij} \bw^4
\partial_t \hat{n}^\alpha \hat{n}_\lambda \Pi^\lambda_\mu
\cp_i  \partial_t^3 v_\alpha\cp_j  \partial_t^3 v^\mu\,dS
+ 
\sigma\int_\Gamma \sqrt{g}g^{ij} \bw^4
 \hat{n}^\alpha \partial_t\hat{n}_\lambda \Pi^\lambda_\mu
\cp_i  \partial_t^3 v_\alpha\cp_j  \partial_t^3 v^\mu\,dS
\nonumber
\\
&= \mathcal{B}_{321} + \mathcal{B}_{322}.
\nonumber
\end{align}
We have
\begin{align}
\norm{\mathcal{B}_{322}} \leq \poly \norm{ \Pi \bw^2 \cp  \partial_t^3 v}_{0,\bou}.
\nonumber
\end{align}
Using Lemma \ref{prelim lemma b}--\ref{partial_t_n} we can write
\begin{align}
\mathcal{B}_{321} & = 
-\sigma\int_\Gamma \bw^4 \sqrt{g}g^{ij} 
g^{kl}   \hat{n}_\lambda  \hat{n}_\tau\cp_k v^\tau \cp_l \eta^\alpha
\Pi^\lambda_\mu
\cp_i  \partial_t^3 v_\alpha\cp_j  \partial_t^3 v^\mu\,dS.
\nonumber
\end{align}
From \eqref{identity_partial_t_v_q} we have
\begin{align}
\cp_l \eta^\alpha \cp_i \partial_t^3 v_\alpha 
& = 
-\frac{J}{\rho_0} \cp_i \cp_l \partial_t^2 q
+ [ \cp_i\partial_t^2, -\frac{J}{\rho_0} \cp_l] q
- [ \cp_i \partial_t^2, \cp_l \eta_\alpha \partial_t ] v^\alpha.
\nonumber
\end{align}
Thus,
\begin{align}
\mathcal{B}_{321} & = 
 \sigma \int_\Gamma \frac{J}{\rho_0} \bw^4 \sqrt{g} g^{ij} g^{kl} \hat{n}_\lambda
\hat{n}_\tau \cp_k v^\tau \cp_i \cp_l \partial_t^2 q \Pi^\lambda_\mu \cp_j \partial_t^3 v^\mu 
\, dS
\nonumber
\\
&
- 
 \sigma \int_\Gamma \frac{1}{\rho_0} \bw^4 \sqrt{g} g^{ij} g^{kl} \hat{n}_\lambda
\hat{n}_\tau \cp_k v^\tau ( [ \cp_i\partial_t^2, -\frac{J}{\rho_0} \cp_l]q )
\Pi^\lambda_\mu \cp_j \partial_t^3 v^\mu 
\, dS
\nonumber
\\
& 
+
 \sigma \int_\Gamma \frac{1}{\rho_0} \bw^4 \sqrt{g} g^{ij} g^{kl} \hat{n}_\lambda
\hat{n}_\tau \cp_k v^\tau ([ \cp_i \partial_t^2, \cp_l \eta_\alpha \partial_t ] v^\alpha)
\Pi^\lambda_\mu \cp_j \partial_t^3 v^\mu 
\, dS
\nonumber
\\
& = \mathcal{B}_{3211} + \mathcal{B}_{3212}+\mathcal{B}_{3213}.
\nonumber
\end{align}
Integrating $\cp_i$ by parts in $\mathcal{B}_{3211}$,
\begin{align}
\mathcal{B}_{3211} & = -\sigma \int_\Gamma \frac{J}{\rho_0} \bw^4 \sqrt{g} g^{ij} g^{kl} \hat{n}_\lambda
\hat{n}_\tau \cp_k v^\tau  \cp_l \partial_t^2 q \Pi^\lambda_\mu \cp_i \cp_j \partial_t^3 v^\mu 
\, dS
\nonumber
\\
&-\sigma \int_\Gamma \cp_i(\frac{J}{\rho_0} \bw^4 \sqrt{g} g^{ij} g^{kl} \hat{n}_\lambda
\hat{n}_\tau \cp_k v^\tau  \cp_l \partial_t^2 q \Pi^\lambda_\mu ) \cp_j \partial_t^3 v^\mu 
\, dS.
\nonumber
\end{align}
From section \ref{sec_the_bry_con}, item \ref{q_cont_bry_con}, we have
\begin{align}
\cp_l \partial_t^2 q= - \sigma g^{mn} \hat{n}_\beta \cp_m \cp_n \cp_l \partial_t v^\beta
-[\cp_l \partial_t^2,\sigma g^{mn} \hat{n}_\beta \cp_m \cp_n ] \eta^\beta,
\nonumber
\end{align}
so that
\begin{align}
\mathcal{B}_{3211} & = \sigma^2 \int_\Gamma \frac{J}{\rho_0} \bw^4 \sqrt{g} g^{ij} g^{kl} \hat{n}_\lambda
\hat{n}_\tau \cp_k v^\tau  
g^{mn} \hat{n}_\beta \cp_m \cp_n \cp_l \partial_t v^\beta
 \Pi^\lambda_\mu \cp_i \cp_j \partial_t^3 v^\mu 
 \nonumber
 \\
 &+ \sigma
 \int_\Gamma \frac{J}{\rho_0} \bw^4 \sqrt{g} g^{ij} g^{kl} \hat{n}_\lambda
\hat{n}_\tau \cp_k v^\tau  
[\cp_l \partial_t^2,\sigma g^{mn} \hat{n}_\beta \cp_m \cp_n ] \eta^\beta
 \Pi^\lambda_\mu \cp_i \cp_j \partial_t^3 v^\mu 
 \nonumber
\, dS
\nonumber
\\
&-\sigma \int_\Gamma \cp_i(\frac{J}{\rho_0} \bw^4 \sqrt{g} g^{ij} g^{kl} \hat{n}_\lambda
\hat{n}_\tau \cp_k v^\tau  \cp_l \partial_t^2 q \Pi^\lambda_\mu ) \cp_j \partial_t^3 v^\mu 
\, dS
\nonumber
\\
& = \mathcal{B}_{32111}+\mathcal{B}_{32112}+\mathcal{B}_{32113}.
\nonumber
\end{align}
In $\mathcal{B}_{32111}$, we 
use Lemma \ref{prelim lemma b}--\ref{n_Pi_contraction} 
 and
factor a $\partial_t$ from $\partial_t^3$ to obtain
\begin{align}
\mathcal{B}_{32111} & = 
\sigma^2 \partial_t \int_\Gamma \frac{J}{\rho_0} \bw^4 \sqrt{g} g^{ij} g^{kl} 
\hat{n}_\tau \cp_k v^\tau  
g^{mn} \hat{n}_\beta \cp_m \cp_n \cp_l \partial_t v^\beta
 \hat{n}_\mu \cp_i \cp_j \partial_t^2 v^\mu 
\nonumber
\\
& - 
\sigma^2  \int_\Gamma \frac{J}{\rho_0} \bw^4 \sqrt{g} g^{ij} g^{kl} 
\hat{n}_\tau \cp_k v^\tau  
g^{mn} \hat{n}_\beta \cp_m \cp_n \cp_l \partial^2_t v^\beta
 \hat{n}_\mu  \cp_i \cp_j \partial_t^2 v^\mu 
 \nonumber
 \\
& -
\sigma^2  \int_\Gamma \partial_t( \frac{J}{\rho_0} \bw^4 \sqrt{g} g^{ij} g^{kl} 
\hat{n}_\tau \cp_k v^\tau  
g^{mn} \hat{n}_\beta 
 \hat{n}_\mu  )\cp_m \cp_n \cp_l \partial_t v^\beta \cp_i \cp_j \partial_t^2 v^\mu 
 \nonumber
 \\
& = 
\mathcal{B}_{321111}+ \mathcal{B}_{321112} + \mathcal{B}_{321113}.
\nonumber
\end{align}
For the first term, i.e., $\mathcal{B}_{321111}$, we have
\begin{align}
\int_0^t \mathcal{B}_{321111} \leq \poly_0 
+  \bw^{\frac{1}{2}}\poly (\bw \norm{\Pi \bw \cp^3 \partial_t v}_{0,\bou})
(\norm{\Pi \bw^{\frac{3}{2}} \cp^2 \partial_t^2 v}_{0,\bou}).
\nonumber
\end{align}
Using Young's inequality and the fact that $\bw $ can be made very small for large 
$\kappa$, we can 
bound the right-hand side by $\poly_0 + \epsilon P(\mathcal{N}) + \epsilon \mathcal{N}$.

For $\mathcal{B}_{321112}$, write
\begin{align}
g^{mn} \hat{n}_\beta \cp_m \cp_n \cp_l \partial^2_t v^\beta
 \hat{n}_\mu g^{ij} \cp_i \cp_j \partial_t^2 v^\mu 
 & =
 \cp_l(g^{mn} \hat{n}_\beta \cp_m \cp_n  \partial^2_t v^\beta)
 \hat{n}_\mu  g^{ij} \cp_i \cp_j \partial_t^2 v^\mu 
 \nonumber
 \\
 & - [ \cp_l , g^{mn} \hat{n}_\beta \cp_m \cp_n  \partial^2_t] v^\beta 
 \hat{n}_\mu  g^{ij} \cp_i \cp_j \partial_t^2 v^\mu 
 \nonumber
\\
 & =
 \frac{1}{2} \cp_l(  \hat{n}_\mu  g^{ij} \cp_i \cp_j \partial_t^2 v^\mu )^2
 \nonumber
 \\
 & - [ \cp_l , g^{mn} \hat{n}_\beta \cp_m \cp_n  \partial^2_t] v^\beta 
 \hat{n}_\mu  g^{ij} \cp_i \cp_j \partial_t^2 v^\mu ,
 \nonumber
\end{align}
so that 
\begin{align}
\mathcal{B}_{321112} & =
 - \frac{1}{2}\sigma^2  \int_\Gamma \frac{J}{\rho_0} \bw^4 \sqrt{g}  g^{kl} 
\hat{n}_\tau \cp_k v^\tau  
 \cp_l(  \hat{n}_\mu  g^{ij} \cp_i \cp_j \partial_t^2 v^\mu )^2
 \nonumber
 \\
 &+ \sigma^2  \int_\Gamma \frac{J}{\rho_0} \bw^4 \sqrt{g}  g^{kl} 
\hat{n}_\tau \cp_k v^\tau 
[ \cp_l , g^{mn} \hat{n}_\beta \cp_m \cp_n  \partial^2_t] v^\beta 
 \hat{n}_\mu  g^{ij} \cp_i \cp_j \partial_t^2 v^\mu .
 \nonumber
\end{align}
Integarting $\cp_l$ by parts in the first integral,
\begin{align}
\mathcal{B}_{321112} & =
  \frac{1}{2}\sigma^2  \int_\Gamma 
  \cp_l(\frac{J}{\rho_0} \bw^4 \sqrt{g}  g^{kl} 
\hat{n}_\tau \cp_k v^\tau  )
(  \hat{n}_\mu  g^{ij} \cp_i \cp_j \partial_t^2 v^\mu )^2
 \nonumber
 \\
 &+ \sigma^2  \int_\Gamma \frac{J}{\rho_0} \bw^4 \sqrt{g}  g^{kl} 
\hat{n}_\tau \cp_k v^\tau 
[ \cp_l , g^{mn} \hat{n}_\beta \cp_m \cp_n  \partial^2_t] v^\beta 
 \hat{n}_\mu  g^{ij} \cp_i \cp_j \partial_t^2 v^\mu 
 \nonumber
 \\
& = \mathcal{B}_{3211121}+ \mathcal{B}_{3211122}.
\nonumber
\end{align}
Writing 
\begin{align}
\mathcal{B}_{3211121} & =
 \bw \frac{1}{2}\sigma^2  \int_\Gamma 
  \cp_l(\frac{J}{\rho_0}  \sqrt{g}  g^{kl} 
\hat{n}_\tau \cp_k v^\tau  )
(  \hat{n}_\mu  g^{ij} \bw^{\frac{3}{2}} \cp_i \cp_j \partial_t^2 v^\mu )^2,
\nonumber
\end{align}
 we have
\begin{align}
\mathcal{B}_{3211121} \leq \epsilon P(\mathcal{N}).
\nonumber
\end{align}

This concludes the estimate for the most delicate terms in $\int_0^t \mathcal{B}_3$.
The remaining terms in $\mathcal{B}_3$, i.e., 
$\mathcal{B}_{3211122}$, $\mathcal{B}_{321113}$, $\mathcal{B}_{32113}$,
$\mathcal{B}_{32112}$, $\mathcal{B}_{3212}$, and $\mathcal{B}_{3213}$,
are treated in section \ref{sec_remainders} below.

\paragraph{Estimate for $\int_0^t \mathcal{B}_2$ with $\dd^4 = (\rr)^2\p_t^4$\label{sec_B2}}
We now move to estimate $\mathcal{B}_2$:
\begin{align}
\mathcal{B}_2 & =
-\sigma\int_\Gamma \cp_i (\sqrt{g}g^{ij}\Pi_\mu^\alpha)(\bw^2   \partial_t^4 v_\alpha)
( \bw^2\cp_j  \partial_t^3 v^\mu)\,dS
\nonumber
\\
& =-\sigma\int_\Gamma \cp_i (\sqrt{g}g^{ij}) \Pi_\mu^\alpha 
\bw^2  \partial_t^4 v_\alpha \bw^2\cp_j  \partial_t^3 v^\mu\,dS
-\sigma\int_\Gamma \sqrt{g}g^{ij}  \cp_i\Pi_\mu^\alpha 
\bw^2   \partial_t^4 v_\alpha \bw^2\cp_j  \partial_t^3 v^\mu\,dS
\nonumber
\\
& = \mathcal{B}_{21} + \mathcal{B}_{22}.
\label{cancellation_2_1}
\end{align}
We show below that  $\mathcal{B}_{21}$ exactly cancels with a term coming from 
$\mathcal{B}_1$. Here we move to estimate $\mathcal{B}_{22}$. Using
Lemma \ref{prelim lemma b}--\ref{Pi_nn},
\begin{align}
\mathcal{B}_{22} & = 
-\sigma\int_\Gamma \bw^4 \sqrt{g}g^{ij}  \cp_i\hat{n}_\mu \hat{n}^\alpha 
 \partial_t^4 v_\alpha \cp_j  \partial_t^3 v^\mu\,dS
 -\sigma\int_\Gamma \bw^4 \sqrt{g}g^{ij}  \hat{n}_\mu \cp_i\hat{n}^\alpha 
  \partial_t^4 v_\alpha \cp_j  \partial_t^3 v^\mu\,dS
 \nonumber
\\
&= \mathcal{B}_{221} + \mathcal{B}_{222}.
\nonumber
\end{align}
We use Lemma \ref{prelim lemma b}--\ref{partial_i_n}
to write
\begin{align}
\mathcal{B}_{221} & = 
\sigma\int_\Gamma \bw^4 \sqrt{g}g^{ij}  g^{kl} \cp_i \cp_k \eta^\tau \hat{n}_\tau \cp_l \eta_\mu
\hat{n}^\alpha 
 \partial_t^4 v_\alpha \cp_j  \partial_t^3 v^\mu\,dS.
\nonumber
\end{align}
From \eqref{identity_partial_t_v_q} we have
\begin{align}
 \cp_l \eta_\mu \cp_j \partial^3_t v^\mu
& =
-\frac{J}{\rho_0} \cp_j \cp_l \partial_t^2 q + [\cp_j \partial_t^2, -\frac{J}{\rho_0}\cp_l]q
-[\cp_j \partial_t^2, \cp_l \eta_\mu \partial_t ] v^\mu ,
\nonumber
\end{align}
whence
\begin{align}
\mathcal{B}_{221} & = 
-\sigma\int_\Gamma \bw^4 \frac{J}{\rho_0} \sqrt{g}g^{ij}  g^{kl} \cp_i \cp_k \eta^\tau \hat{n}_\tau 
\hat{n}^\alpha 
 \partial_t^4 v_\alpha 
  \cp_j \cp_l \partial_t^2 q
 \,dS
\nonumber
\\
&
+\sigma\int_\Gamma \bw^4 \sqrt{g}g^{ij}  g^{kl} \cp_i \cp_k \eta^\tau \hat{n}_\tau 
\hat{n}^\alpha 
  \partial_t^4 v_\alpha 
 [\cp_j \partial_t^2, -\frac{J}{\rho_0}\cp_l]q
 \,dS
 \nonumber
 \\
 &
-\sigma\int_\Gamma \bw^4 \sqrt{g}g^{ij}  g^{kl} \cp_i \cp_k \eta^\tau \hat{n}_\tau 
\hat{n}^\alpha 
 \partial_t^4 v_\alpha 
 [\cp_j \partial_t^2, \cp_l \eta_\mu \partial_t ] v^\mu 
 \,dS 
 \nonumber
 \\
 & = \mathcal{B}_{2211}+\mathcal{B}_{2212}+\mathcal{B}_{2213}.
 \nonumber
\end{align}
In $\mathcal{B}_{2211}$, we factor a $\partial_t$ in $\partial_t^4 v_\alpha$ to obtain
\begin{align}
\mathcal{B}_{2211} & =
- \sigma \partial_t \int_\Gamma \bw^4 \frac{J}{\rho_0} \sqrt{g}g^{ij}  g^{kl} \cp_i \cp_k \eta^\tau \hat{n}_\tau 
\hat{n}^\alpha 
 \partial_t^3 v_\alpha 
  \cp_j \cp_l \partial_t^2 q
 \,dS
\nonumber
\\
&+ \sigma\int_\Gamma \bw^4 \partial_t( \frac{J}{\rho_0} \sqrt{g}g^{ij}  g^{kl} \cp_i \cp_k \eta^\tau \hat{n}_\tau 
\hat{n}^\alpha 
  \cp_j \cp_l \partial_t^2 q )  \partial_t^3 v_\alpha 
 \,dS
\nonumber
\\
& = \mathcal{B}_{22111} + \mathcal{B}_{22112}.
\nonumber
\end{align}
For $\mathcal{B}_{22111}$, we integrate $\cp_j$ by parts to produce
\begin{align}
\int_0^t \mathcal{B}_{22111} \leq \poly 
(\norm{ \bw^2 \Pi \cp \partial_t^3 v}_{0,\bou})
(\norm{\bw^2\cp \partial_t^2 q}_{0,\bou} )
+ \int_0^t \poly,
\nonumber 
\end{align}
where 
\begin{align}
\norm{ \bw^2 \Pi \cp \partial_t^3 v}_{0,\bou}\norm{\bw^2\cp \partial_t^2 q}_{0,\bou} \lesssim \tilde{\epsilon} \norm{ \bw^2 \Pi \cp \partial_t^3 v}_{0,\bou}+ \bw^{\frac{1}{2}}\tilde{\epsilon}^{-1}\norm{\bw^{\frac{3}{2}}\cp \partial_t^2 q}_{1}\\
\lesssim \epsilon (\norm{ \bw^2 \Pi \cp \partial_t^3 v}_{0,\bou}+\norm{\bw^{\frac{3}{2}}\cp \partial_t^2 q}_{1}) \lesssim \epsilon P(\mathcal{N}),
\end{align}
after choosing $\bw$ sufficiently small and replacing $q$ by $R$.

For $\mathcal{B}_{22112}$, write
\begin{align}
\mathcal{B}_{22112} & =  \sigma\int_\Gamma \bw^4 \partial_t( \frac{J}{\rho_0} \sqrt{g}g^{ij}  g^{kl} \cp_i \cp_k \eta^\tau \hat{n}_\tau 
\hat{n}^\alpha 
   )\cp_j \cp_l \partial_t^2 q  \partial_t^3 v_\alpha 
 \,dS
 \nonumber
 \\
&+ 
   \sigma\int_\Gamma \bw^4  \frac{J}{\rho_0} \sqrt{g}g^{ij}  g^{kl} \cp_i \cp_k \eta^\tau \hat{n}_\tau 
\hat{n}^\alpha 
   \cp_j \cp_l \partial_t^3 q \partial_t^3 v_\alpha 
 \,dS
 \nonumber
 \\
 & = \mathcal{B}_{221121} + \mathcal{B}_{221122}.
 \nonumber
\end{align}
The term $\mathcal{B}_{221121}$ can be handled with integration by parts with respect
to $\cp_j$ (it yields a term in $\norm{\Pi \bw^2 \cp \partial_t^3 v}_{0,\bou}$). For
$\mathcal{B}_{221122}$, we use
section \ref{sec_the_bry_con}, item \ref{q_cont_bry_con},
to write
\begin{align}
\mathcal{B}_{221122} & = 
-\sigma^2\int_\Gamma \bw^4  \frac{J}{\rho_0} \sqrt{g}g^{ij}  g^{kl} 
\cp_i \cp_k \eta^\tau \hat{n}_\tau 
\hat{n}^\alpha \cp_j \cp_l ( g^{mn} \hat{n}_\beta \partial_m \partial_n \partial_t^2 v^\beta) 
 \partial_t^3 v_\alpha 
\,dS
\nonumber
\\
& 
-\sigma^2\int_\Gamma \bw^4  \frac{J}{\rho_0} \sqrt{g}g^{ij}  g^{kl} 
\cp_i \cp_k \eta^\tau \hat{n}_\tau 
\hat{n}^\alpha (\cp_j \cp_l [\partial^3_t, g^{mn}\hat{n}_\beta \cp_m \cp_n ] \eta^\beta)  \partial_t^3 v_\alpha 
\nonumber
\\
&=\mathcal{B}_{2211221}+\mathcal{B}_{2211222}.
\nonumber
\end{align}
Integrating by parts $\cp_l$ in
$\mathcal{B}_{2211221}$,
\begin{align}
\mathcal{B}_{2211221} & = 
\sigma^2\int_\Gamma \bw^4  \frac{J}{\rho_0} \sqrt{g}g^{ij}  g^{kl} 
\cp_i \cp_k \eta^\tau \hat{n}_\tau 
\hat{n}^\alpha \cp_j  ( g^{mn} \hat{n}_\beta \partial_m \partial_n \partial_t^2 v^\beta) 
\cp_l \partial_t^3 v_\alpha 
\,dS
\nonumber
\\
& +
\sigma^2\int_\Gamma \cp_l(\bw^4  \frac{J}{\rho_0} \sqrt{g}g^{ij}  g^{kl} 
\cp_i \cp_k \eta^\tau \hat{n}_\tau 
\hat{n}^\alpha) \cp_j  ( g^{mn} \hat{n}_\beta \partial_m \partial_n \partial_t^2 v^\beta)   \partial_t^3 v_\alpha 
\,dS
\nonumber
\\
& =
\sigma^2\int_\Gamma \bw^4  \frac{J}{\rho_0} \sqrt{g}g^{ij}  g^{kl} 
\cp_i \cp_k \eta^\tau \hat{n}_\tau 
\hat{n}^\alpha \cp_j  ( g^{mn} \hat{n}_\beta \partial_m \partial_n \partial_t^2 v^\beta) 
\cp_l \partial_t^3 v_\alpha 
\,dS
\nonumber
\\
& +
\sigma^2\int_\Gamma \bw^4 \frac{J}{\rho_0} \sqrt{g}g^{ij}  g^{kl} 
\cp_l\cp_i \cp_k \eta^\tau \hat{n}_\tau 
\hat{n}^\alpha \cp_j  ( g^{mn} \hat{n}_\beta \partial_m \partial_n \partial_t^2 v^\beta)   \partial_t^3 v_\alpha 
\,dS
\nonumber
\\
& +
\sigma^2\int_\Gamma \cp_l(\bw^4  \frac{J}{\rho_0} \sqrt{g}g^{ij}  g^{kl} 
 \hat{n}_\tau 
\hat{n}^\alpha) \cp_i \cp_k \eta^\tau \cp_j  ( g^{mn} \hat{n}_\beta \partial_m \partial_n \partial_t^2 v^\beta)   \partial_t^3 v_\alpha 
\,dS,
\nonumber
\end{align}
and then integrating by parts $\cp_i$ on the second integral,
\begin{align}
\mathcal{B}_{2211221} &=
\sigma^2\int_\Gamma \bw^4 \frac{J}{\rho_0} \sqrt{g}g^{ij}  g^{kl} 
\cp_i \cp_k \eta^\tau \hat{n}_\tau 
\hat{n}^\alpha \cp_j  ( g^{mn} \hat{n}_\beta \partial_m \partial_n \partial_t^2 v^\beta) 
\cp_l \partial_t^3 v_\alpha 
\,dS
\nonumber
\\
& -
\sigma^2\int_\Gamma \bw^4  \frac{J}{\rho_0} \sqrt{g}g^{ij}  g^{kl} 
\cp_l \cp_k \eta^\tau \hat{n}_\tau 
\hat{n}^\alpha \cp_i\cp_j  ( g^{mn} \hat{n}_\beta \partial_m \partial_n \partial_t^2 v^\beta)   \partial_t^3 v_\alpha 
\,dS
\nonumber
\\
& -
\sigma^2\int_\Gamma \bw^4  \frac{J}{\rho_0} \sqrt{g}g^{ij}  g^{kl} 
\cp_l \cp_k \eta^\tau \hat{n}_\tau 
\hat{n}^\alpha \cp_j  ( g^{mn} \hat{n}_\beta \partial_m \partial_n \partial_t^2 v^\beta)  
\cp_i \partial_t^3 v_\alpha 
\,dS
\nonumber
\\
& -
\sigma^2\int_\Gamma \cp_i( \bw^4  \frac{J}{\rho_0} \sqrt{g}g^{ij}  g^{kl} 
\hat{n}_\tau 
\hat{n}^\alpha )
\cp_l \cp_k \eta^\tau \cp_j  ( g^{mn} \hat{n}_\beta \partial_m \partial_n \partial_t^2 v^\beta)  
 \partial_t^3 v_\alpha 
\,dS
\nonumber
\\
& +
\sigma^2\int_\Gamma \cp_l(\bw^4  \frac{J}{\rho_0} \sqrt{g}g^{ij}  g^{kl} 
 \hat{n}_\tau 
\hat{n}^\alpha) \cp_i \cp_k \eta^\tau \cp_j  ( g^{mn} \hat{n}_\beta \partial_m \partial_n \partial_t^2 v^\beta)   \partial_t^3 v_\alpha 
\,dS
\nonumber
\\
& = 
\mathcal{B}_{22112211}+
\mathcal{B}_{22112212}+
\mathcal{B}_{22112213}+
\mathcal{B}_{22112214}+
\mathcal{B}_{22112215}.
\label{cancellation_1}
\end{align}
Note that the first and third terms, i.e., 
$\mathcal{B}_{22112211}$ and $\mathcal{B}_{22112213}$, cancel each other in view of the
following identity, which can be verified by inspection,
\begin{align}
\sum_{i,k,l=1}^2 (g^{ij} g^{kl} - g^{ik} g^{lj} ) = 0.
\nonumber
\end{align}
For  the second term,
$\mathcal{B}_{22112212}$, integrate $\cp_j \cp_j$ by parts:
\begin{align}
\mathcal{B}_{22112212} & = 
-
\sigma^2\int_\Gamma \bw^4  \frac{J}{\rho_0} \sqrt{g}g^{ij}  g^{kl} 
\cp_l \cp_k \eta^\tau \hat{n}_\tau 
  g^{mn} \hat{n}_\beta \partial_m \partial_n \partial_t^2 v^\beta 
\hat{n}^\alpha \cp_i\cp_j \partial_t^3 v_\alpha 
\,dS
\nonumber
\\
&-
\sigma^2\int_\Gamma \cp_j(\bw^4  \frac{J}{\rho_0} \sqrt{g}g^{ij}  g^{kl} 
\cp_l \cp_k \eta^\tau \hat{n}_\tau 
) g^{mn} \hat{n}_\beta \partial_m \partial_n \partial_t^2 v^\beta  \hat{n}^\alpha\cp_i   \partial_t^3 v_\alpha 
\,dS
\nonumber
\\
&-
\sigma^2\int_\Gamma \cp_i (\bw^4  \frac{J}{\rho_0} \sqrt{g}g^{ij}  g^{kl} 
\cp_l \cp_k \eta^\tau \hat{n}_\tau 
\hat{n}^\alpha )\cp_j(  g^{mn} \hat{n}_\beta \partial_m \partial_n \partial_t^2 v^\beta  ) \partial_t^3 v_\alpha 
\,dS 
\nonumber
\\
& = 
\mathcal{B}_{221122121}+
\mathcal{B}_{221122122}+
\mathcal{B}_{221122123}.
\nonumber
\end{align}
Factoring a $\partial_t$ from $\cp_i \cp_j \partial_t^3 v_\alpha$ in $\mathcal{B}_{221122121}$, we find

\begin{align}
\mathcal{B}_{221122121} & =
-
\frac{1}{2}\sigma^2\int_\Gamma \bw^4  \frac{J}{\rho_0} \sqrt{g}g^{ij}  g^{kl} 
\cp_l \cp_k \eta^\tau \hat{n}_\tau 
 g^{mn} \hat{n}_\beta \partial_m \partial_n \partial_t^2 v^\beta  
\hat{n}^\alpha \cp_i\cp_j \partial_t^2 v_\alpha 
\,dS
\nonumber
\\
& +
\frac{1}{2}\sigma^2\int_\Gamma \partial_t (\bw^4  \frac{J}{\rho_0} \sqrt{g}g^{ij}  g^{kl} 
\cp_l \cp_k \eta^\tau \hat{n}_\tau 
 g^{mn} \hat{n}_\beta \hat{n}^\alpha) \partial_m \partial_n \partial_t^2 v^\beta  
 \cp_i\cp_j \partial_t^2 v_\alpha 
\,dS
\nonumber
\\
& = 
\mathcal{B}_{2211221211}+ 
\mathcal{B}_{2211221212}. 
\nonumber
\end{align}
The first term, $\mathcal{B}_{2211221211}$, can be estimated by 
$\epsilon P(\mathcal{N})$. Here, the small number $\epsilon$ comes form estimating 
$\cp_l \cp_k \eta^\tau$ in $L^\infty$ and using that $\eta(0)$ is the identity diffeomorphism.

Now we move to $\mathcal{B}_{222}$. Factoring a $\partial_t$ from 
$\partial_t^4 v_\alpha$, we find
\begin{align}
\mathcal{B}_{222} & =
-\sigma\int_\Gamma \bw^4 \sqrt{g}g^{ij}  \hat{n}_\mu \cp_i\hat{n}^\alpha 
  \partial_t^4 v_\alpha \cp_j  \partial_t^3 v^\mu\,dS
\nonumber
\\
& =
- \sigma \partial_t\int_\Gamma \bw^4 \sqrt{g}g^{ij}  \hat{n}_\mu \cp_i\hat{n}^\alpha 
  \partial_t^3 v_\alpha \cp_j  \partial_t^3 v^\mu\,dS
\nonumber
  \\
&
+\sigma\int_\Gamma \bw^4 \sqrt{g}g^{ij}  \hat{n}_\mu \cp_i\hat{n}^\alpha 
  \partial_t^3 v_\alpha \cp_j  \partial_t^4 v^\mu\,dS
\nonumber
\\
& 
+\sigma\int_\Gamma \partial_t( \bw^4 \sqrt{g}g^{ij}  \hat{n}_\mu \cp_i\hat{n}^\alpha) 
  \partial_t^3 v_\alpha \cp_j  \partial_t^3 v^\mu\,dS.
\nonumber
\end{align}
Integrating $\cp_j$ by parts in the second integral,
\begin{align}
\mathcal{B}_{222} 
& =
- \sigma \partial_t\int_\Gamma \bw^4 \sqrt{g}g^{ij}  \hat{n}_\mu \cp_i\hat{n}^\alpha 
  \partial_t^3 v_\alpha \cp_j  \partial_t^3 v^\mu\,dS
\nonumber
  \\
&
-\sigma\int_\Gamma \bw^4 \sqrt{g}g^{ij}  \hat{n}_\mu \cp_i\hat{n}^\alpha 
  \cp_j\partial_t^3 v_\alpha   \partial_t^4 v^\mu\,dS
\nonumber
\\
& 
+\sigma\int_\Gamma \partial_t( \bw^4 \sqrt{g}g^{ij}  \hat{n}_\mu \cp_i\hat{n}^\alpha) 
  \partial_t^3 v_\alpha \cp_j  \partial_t^3 v^\mu\,dS.
\nonumber
\\
&-\sigma\int_\Gamma   \cp_j(\bw^4 \sqrt{g}g^{ij}  \hat{n}_\mu \cp_i\hat{n}^\alpha )
\partial_t^3 v_\alpha   \partial_t^4 v^\mu\,dS
\nonumber
\\
&= 
\mathcal{B}_{2221}+
\mathcal{B}_{2222}+
\mathcal{B}_{2223}+
\mathcal{B}_{2224}.
\nonumber
\end{align}
Note that $\mathcal{B}_{2222} = \mathcal{B}_{221}$, so this term is estimated as above.
The term $\mathcal{B}_{2221}$ can, after time integration, be estimated
using Young's inequality and interpolation.

With exception of $\mathcal{B}_{21}$, which, as said, involves a special cancellation 
showed below, this concludes the estimate of the most delicate terms in 
$\int_0^t \mathcal{B}_2$. The remaining terms
$\mathcal{B}_{2212}$,
$\mathcal{B}_{2213}$,
$\mathcal{B}_{2211222}$,
$\mathcal{B}_{22112214}$,
$\mathcal{B}_{22112215}$,
$\mathcal{B}_{221122122}$,
$\mathcal{B}_{221122123}$,
$\mathcal{B}_{2211221212}$,
$\mathcal{B}_{2223}$,
and
$\mathcal{B}_{2224}$,
are treated in section \ref{sec_remainders} below.

\paragraph{Estimate for $\int_0^t \mathcal{B}_1$ with $\dd^4 = (\rr)^2\p_t^4$\label{sec_B1}}
We now move to estimate $\mathcal{B}_1$:
\begin{align}
\mathcal{B}_1 & = \sigma\int_{\Gamma}(\bw^2 \partial_t^4 v_\alpha) \big([\bw^2 \partial_t^4, \sqrt{g}g^{ij}\Pi^\alpha_\mu]\cp^2_{ij}\eta^\mu\big)\,dS
\nonumber
\\
& = 
4\sigma\int_{\Gamma}\bw^4 \partial_t(\sqrt{g} g^{ij} \Pi^\alpha_\mu ) \cp_i \cp_j \partial_t^3 \eta^\mu 
\partial_t^4 v_\alpha
\nonumber
\\
&+6\sigma\int_{\Gamma}
\bw^4 \partial^2_t(\sqrt{g} g^{ij} \Pi^\alpha_\mu ) \cp_i \cp_j \partial_t^2 \eta^\mu 
\partial_t^4 v_\alpha
\nonumber
\\
&+4\sigma\int_{\Gamma}
\bw^4 \partial^3_t(\sqrt{g} g^{ij} \Pi^\alpha_\mu ) \cp_i \cp_j \partial_t \eta^\mu 
\partial_t^4 v_\alpha
\nonumber
\\
&+\sigma\int_{\Gamma}
\bw^4 \partial^4_t(\sqrt{g} g^{ij} \Pi^\alpha_\mu ) \cp_i \cp_j \eta^\mu 
\partial_t^4 v_\alpha
\nonumber
\\
& = \mathcal{B}_{11}+\mathcal{B}_{12}+\mathcal{B}_{13}+\mathcal{B}_{14}.
\nonumber
\end{align}
We have
\begin{align}
\mathcal{B}_{14} & =
\sigma\int_{\Gamma}\bw^4 
\sqrt{g} g^{ij} \partial_t^4 \Pi^\alpha_\mu 
 \cp_i \cp_j \eta^\mu 
\partial_t^4 v_\alpha
\nonumber
\\
&+
4\sigma\int_{\Gamma}\bw^4 
\partial_t(\sqrt{g} g^{ij}) \partial_t^3 \Pi^\alpha_\mu 
 \cp_i \cp_j \eta^\mu 
\partial_t^4 v_\alpha
\nonumber
\\
&+
6\sigma\int_{\Gamma}\bw^4
\partial^2_t(\sqrt{g} g^{ij}) \partial_t^2 \Pi^\alpha_\mu 
  \cp_i \cp_j \eta^\mu 
\partial_t^4 v_\alpha
\nonumber
\\
&+
4\sigma\int_{\Gamma}\bw^4 
\partial_t^3(\sqrt{g} g^{ij}) \partial_t \Pi^\alpha_\mu 
 \cp_i \cp_j \eta^\mu 
\partial_t^4 v_\alpha
\nonumber
\\
&+
\sigma\int_{\Gamma}\bw^4  
\partial_t^4(\sqrt{g} g^{ij})  \Pi^\alpha_\mu 
\cp_i \cp_j \eta^\mu 
\partial_t^4 v_\alpha
\nonumber
\\
&=\mathcal{B}_{141}+\mathcal{B}_{142}+\mathcal{B}_{143}+\mathcal{B}_{144}+\mathcal{B}_{145}.
\nonumber
\end{align}
Using Lemma \ref{prelim lemma b}--\ref{Pi_nn},
we have
\begin{align}
\mathcal{B}_{141} & =
\sigma\int_{\Gamma} \bw^4 
\sqrt{g} g^{ij}  \hat{n}^\alpha \partial_t^4 \hat{n}_\mu
 \cp_i \cp_j \eta^\mu 
\partial_t^4 v_\alpha
\nonumber
\\
&
+
4\sigma\int_{\Gamma} \bw^4 
\sqrt{g} g^{ij}  \partial_t\hat{n}^\alpha \partial_t^3 \hat{n}_\mu
 \cp_i \cp_j \eta^\mu 
\partial_t^4 v_\alpha
\nonumber
\\
&
+
6\sigma\int_{\Gamma} \bw^4 
\sqrt{g} g^{ij}  \partial^2_t\hat{n}^\alpha \partial_t^2 \hat{n}_\mu
 \cp_i \cp_j \eta^\mu 
\partial_t^4 v_\alpha
\nonumber
\\
&
+
4\sigma\int_{\Gamma} \bw^4 
\sqrt{g} g^{ij}  \partial^3_t\hat{n}^\alpha \partial_t \hat{n}_\mu
 \cp_i \cp_j \eta^\mu 
\partial_t^4 v_\alpha
\nonumber
\\
&
+
\sigma\int_{\Gamma} \bw^4 
\sqrt{g} g^{ij}  \partial^4_t\hat{n}^\alpha  \hat{n}_\mu
 \cp_i \cp_j \eta^\mu 
\partial_t^4 v_\alpha
\nonumber
\\
&= \mathcal{B}_{1411}+\mathcal{B}_{1412}+\mathcal{B}_{1413}+\mathcal{B}_{1414}+\mathcal{B}_{1415}.
\nonumber
\end{align}
From Lemma \ref{prelim lemma b}--\ref{partial_t_n} we have
\begin{align}
\partial_t^4 \hat{n}_\mu = -g^{kl} \cp_k \partial_t^3 v^\tau \hat{n}_\tau 
\cp_l \eta_\mu -[\partial_t^3, g^{kl} \hat{n}_\tau \cp_l \eta_\mu \cp_k]v^\tau,
\nonumber
\end{align}
and thus
\begin{align}
\mathcal{B}_{1411} & = 
-\sigma\int_{\Gamma} \bw^4 
\sqrt{g} g^{ij}  \hat{n}^\alpha
g^{kl} \cp_k \partial_t^3 v^\tau \hat{n}_\tau 
\cp_l \eta_\mu
 \cp_i \cp_j \eta^\mu 
\partial_t^4 v_\alpha
\nonumber
\\
&
-\sigma\int_{\Gamma} \bw^4 
\sqrt{g} g^{ij}  \hat{n}^\alpha
[\partial_t^3, g^{kl} \hat{n}_\tau \cp_l \eta_\mu \cp_k]v^\tau
 \cp_i \cp_j \eta^\mu 
\partial_t^4 v_\alpha
\nonumber
\\
& = \mathcal{B}_{14111} +\mathcal{B}_{14112}.
\nonumber
\end{align}
We now invoke Lemma \ref{prelim lemma b}--\ref{partial_i_ggij}, 
to replace $\sqrt{g} g^{ij} g^{kl} \cp_i \cp_j \eta^\mu \cp_l \eta_\mu $ in 
$\mathcal{B}_{14111}$ by 
$-\cp_i(\sqrt{g} g^{ik} )$, obtaining
\begin{align}
\mathcal{B}_{14111} & = 
\sigma\int_{\Gamma} \bw^4 
\partial_i(\sqrt{g} g^{ik} )
 \cp_k \partial_t^3 v^\tau \hat{n}_\tau 
\hat{n}^\alpha\partial_t^4 v_\alpha.
\label{cancellation_2_2}
\end{align}
We see that this term exactly cancels $\mathcal{B}_{21}$, as mentioned earlier. The other
terms in the estimate of $\int_0^t \mathcal{B}_1$ are treated in section
\ref{sec_remainders}.

\paragraph{Remainders in $\int_0^t \mathcal{B}$ with $\dd^4 = \rr^2\p_t^4$\label{sec_remainders}}
Above we have showed how to control the most delicate terms in the estimate for $\int_0^t \mathcal{B}$
when $\dd^4 = \bw^2 \partial_t^4$. In particular, we have showed how some
top order terms, which seemingly cannot be individually bounded, cancel out when taken together.
Now we consider the remaining terms, which we list here for the reader's convenience. They are,
for
$\mathcal{B}_3$, 
\begin{gather}
\mathcal{B}_{3211122}, \, \mathcal{B}_{321113}, \, \mathcal{B}_{32113},\,
\mathcal{B}_{32112},\,  \mathcal{B}_{3212}, \text{ and } \mathcal{B}_{3213}
\nonumber
\end{gather}
from section \ref{sec_B1}; for $\mathcal{B}_2$
\begin{gather}
\mathcal{B}_{2212}, \, 
\mathcal{B}_{2213},\, 
\mathcal{B}_{2211222},\, 
\mathcal{B}_{22112214},\, 
\mathcal{B}_{22112215},
\nonumber
\\
\mathcal{B}_{221122122},\, 
\mathcal{B}_{221122123},\, 
\mathcal{B}_{2211221212},\, 
\mathcal{B}_{2223},\, 
\text{ and }
\mathcal{B}_{2224}
\nonumber
\end{gather}
from section \ref{sec_B2};
for $\mathcal{B}_1$
\begin{gather}
\mathcal{B}_{11}, \, 
\mathcal{B}_{12},\, 
\mathcal{B}_{13},\, 
\mathcal{B}_{142},\, 
\mathcal{B}_{143},\, 
\mathcal{B}_{144},\, 
\mathcal{B}_{145},\, 
\mathcal{B}_{1412},\, 
\mathcal{B}_{1413},\, 
\mathcal{B}_{1414},\, 
\mathcal{B}_{1415}, \, 
\text{ and }
\mathcal{B}_{14112}
\label{B_1_remaining}
\end{gather}
from section \ref{sec_B1}. 
Not all these terms are immediately of lower order, but they can be estimated using the same kind
of ideas that have already been employed. Therefore, it suffices to briefly indicate how this is done.

The terms 
$\mathcal{B}_{3212}$, 
$\mathcal{B}_{3213}$, 
$\mathcal{B}_{32112}$, 
$\mathcal{B}_{321113}$, 
and
$\mathcal{B}_{3211122}$
can be bounded directly.
The term $\mathcal{B}_{32113}$ is bounded upon replacing $q$ by $R$ and estimating in routine 
fashion. 

The terms 
$\mathcal{B}_{2212}$ and 
$\mathcal{B}_{2213}$ can be estimated with integration by parts in time.
The terms $\mathcal{B}_{2211222}$, 
$\mathcal{B}_{22112214}$, 
$\mathcal{B}_{22112215}$,  
$\mathcal{B}_{221122122}$, 
$\mathcal{B}_{2211221212}$
$\mathcal{B}_{2223}$, and 
$\mathcal{B}_{2224}$
can be estimated directly.
The term $\mathcal{B}_{221122123}$  requires integration by parts in space
and then using arguments similar to above, with one extra step: after integrating $\cp_j$ by parts,
we obtain a term with four derivatives of $\eta$. This term, however, 
has the form $\hat{n}_\tau  g^{ij} \cp^2 \cp_i \cp_j \eta^\tau$, which allows us to use
section \ref{sec_the_bry_con}, item \ref{q_cont_bry_con}, to eliminate two derivatives of $\eta$.
(Alternatively, we can use elliptic estimates for equations with Sobolev coefficients, as, e.g., Theorem 4 and Remark 2 in \cite{milani1983regularity}).

The terms listed in \eqref{B_1_remaining} are again handled by a repetition of ideas used above 
(without requiring special cancellations). In particular, 
Lemma \ref{prelim lemma b}--\ref{partial_t_n} is used heavily and 
Lemma \ref{prelim lemma b}--\ref{partial_t_ggij} is employed to estimate
$\mathcal{B}_{145}$.

Combining these observations with the estimates of section
\ref{sec_B3}, \ref{sec_B2}, and \ref{sec_B3}, we finally obtain
\begin{align}
\int_0^t (\mathcal{B}_1 + \mathcal{B}_2 + \mathcal{B}_3 ) 
\leq \poly_0  + \epsilon \mathcal{N} + \epsilon P(\mathcal{N}) +
\poly \int_0^t \poly, \text{ when } \dd^4 = \bw^2 \partial_t^4.
\nonumber
\end{align}

\subsubsection{Estimate of the remaining weighed boundary terms}

It remains to carry out control of 
$\int_0^t (\mathcal{B}_1 + \mathcal{B}_2 + \mathcal{B}_3 ) $ when
$\dd^2 =  \sqrt{\rr}\p_t^2$,
$\dd^3 =   \sqrt{\rr}(\cp\p_t^2), $ $\dd^3 = \rr\p_t^3$,
$\dd^4 = \rr(\cp^3\p_t)$, $\dd^4 = \rr(\cp^2\p_t^2) $, and $\dd^4 = (\rr)^{\frac{3}{2}}(\cp\p_t^3)$.
These cases are treated in an almost identical fashion as the case
$\dd^4 = (\rr)^2\p_t^4$ from section \ref{sec_weighted_main_case}. In this regard, we note 
that a crucial requirement to carry these estimates is that $\dd$ contains at least one time
derivative, which is the case for all the $\rr$-weighted derivatives we need to 
consider\footnote{Incidentally, this is why an estimate for the normal component of $v$ with
no time derivatives has to be obtained in a different way, see section \ref{sec_curl_boundary_bounds}}.
We therefore conclude
\begin{align}
\int_0^t (\mathcal{B}_1 + \mathcal{B}_2 + \mathcal{B}_3 ) 
\leq \poly_0  + \epsilon \mathcal{N} + \epsilon P(\mathcal{N}) +
\poly \int_0^t \poly, 
\nonumber
\end{align}
for $\dd^2 =  \sqrt{\rr}\p_t^2$,
$\dd^3 =   \sqrt{\rr}(\cp\p_t^2), $ $\dd^3 = \rr\p_t^3$,
$\dd^4 = \rr(\cp^3\p_t)$, $\dd^4 = \rr(\cp^2\p_t^2) $,  $\dd^4 = (\rr)^{\frac{3}{2}}(\cp\p_t^3)$,
and $\dd^4 = (\rr)^2\p_t^4$.

\section{Closing the estimate} \label{section 4}
In this section, we prove:
\thm \label{close esi thm} Let $\mathcal{N}(t)$ and $\PP(t)$ be defined as  Notation \ref{nota PP N}, then for sufficiently large $\kk$ (i.e., $\rr\ll 1$), we have:
\begin{align}
\mathcal{N}(t) \leq C(M)\bigg(\epsilon P(\mathcal{N}(t))+\PP_0+\PP\int_0^t\PP\bigg),\q t\in[0,T],\label{close est}
\end{align}
where $T>0$ is chosen sufficiently small, provided that:
\begin{align}
||\p \eta||_{L^{\infty}}+||\p^2 \eta||_{L^{\infty}}\leq M,\label{apriori v eta}\\
||g^{ij}||_{L^{\infty}}+||\Gamma_{ij}^k||_{L^{\infty}} \leq M,\label{a priori g gamma}
\end{align}
hold a priori for some large constant $M$.

Since the energy estimate for $E$ is established in the previous section(i.e., Theorem \ref{energy estimate E}), we only need to show 
\begin{align}
||v||_4^2+||\rr v_t||_3^2+||\rr v_{tt}||_2^2+||(\rr)^{\frac{3}{2}}v_{ttt}||_1^2
+||R||_4^2+||R_{t}||_3^2+||\sqrt{\rr}R_{tt}||_2^2+||\rr R_{ttt}||_1^2\no\\
+||v_t||_2^2+||\sqrt{\rr}v_{tt}||_1^2+||\rr v_{ttt}||_0^2+||R_{tt}||_1+||\sqrt{\rr} R_{ttt}||_0\no\\
 \leq C(M)\bigg(\epsilon P(\mathcal{N}(t))+\PP_0+\PP\int_0^t\PP\bigg).
 \label{closing}
\end{align}
This is proved via an iterated argument using div-curl estimate \eqref{div curl normal}.  It suffices to consider the first line in \eqref{closing}, since the second line consists lower order terms and can be treated by the same method. Taking $X=v$ and $s=4$, \eqref{div curl normal} yields
\begin{align}
||v||_4^2 \lesssim ||\di v||_3^2+||\curl v||_3^2+ ||v^3||_{3.5, \p}^2+||v||_0^2.\label{v_4}
\end{align}
On the other hand, taking $X=\rr\p_t v$ and $s=3$, we have:
\begin{align}
||\rr v_t||_3^2 \lesssim ||\rr \di  v_t||_2^2+||\rr \curl  v_t||_2^2+||\rr v_t^3||_{2.5,\Gamma}^2+||\rr v_t||_0^2. 
\end{align}
Similarly, by taking $X=R'v_{tt}$, $s=2$ and $X=(R')^{\frac{3}{2}}v_{ttt}$, $s=1$, we get
\begin{align}
||\rr v_{tt}||_2^2 \lesssim ||\rr\di v_{tt}||_1^2+||\rr\curl v_{tt}||_1^2+||\rr v_{tt}^3||_{1.5,\Gamma}^2+||\rr v_{tt}||_0^2,\label{v_tt_2}\\
||(\rr)^{\frac{3}{2}}v_{ttt}||_1^2 \lesssim ||(\rr)^{\frac{3}{2}}\di v_{ttt}||_0^2+||(\rr)^{\frac{3}{2}}\curl v_{ttt}||_0^2+||(\rr)^{\frac{3}{2}}v_{ttt}^3||_{0.5,\Gamma}^2+||(\rr)^{\frac{3}{2}}v_{ttt}||_0^2,\label{v_ttt_1}
\end{align}
respectively. In light of \eqref{v_4}-\eqref{v_ttt_1}, in order to estimate $v$ and its time derivative, we need to bound $\di \p_t^k v$, $\curl \p_t^k v$ and $\p_t^k v^3$, for $k=0,1,2,3$, respectively. 
\subsection{Bounds for the curl and the boundary term of $v$\label{sec_curl_boundary_bounds}}
In this section we prove:
\thm 
\begin{align}
||\curl  v||_{3}^2+||\rr\curl v_t||_2^2+||\rr\curl v_{tt}||_1^2+||(\rr)^{\frac{3}{2}}\curl v_{ttt}||_0^2 \lesssim \epsilon P(\mathcal{N})+\PP_0+\PP\int_0^t \PP.\label{curl}
\end{align}
\begin{proof}
The proof is almost identical to section 4 of \cite{disconzi2017prioriC}, and so we omit the details. The only modification is that the weights $\rr$ or $(\rr)^{\frac{3}{2}}$ are used to compensate $q'(R)\sim \rr^{-1}$, which allows us to get an uniform control. 
\end{proof}
On the other hand, we have:
\thm 
\label{thm estimates for closing}
\begin{align}
||v^3||_{3.5,\Gamma}^2 \lesssim \epsilon P(\mathcal{N})+\PP_0+\PP\int_0^t \PP,\label{v^3}
\end{align}
and
\begin{align}
||\rr v_t^3||_{2.5,\Gamma}^2 \lesssim   \epsilon\mathcal{N}+||\rr \Pi \cp^3 v_t||_{0,\Gamma}^2 +\PP_0+\PP\int_0^t \PP\label{v_t^3},\\
||\rr v_{tt}^3||_{1.5,\Gamma}^2 \lesssim   \epsilon\mathcal{N}+||\rr \Pi \cp^2 v_{tt}||_{0,\Gamma}^2 +\PP_0+\PP\int_0^t \PP,\label{v_tt^3}\\
||(\rr)^{\frac{3}{2}}v_{ttt}^3||_{0.5,\Gamma}^2 \lesssim   \epsilon\mathcal{N}+||(\rr)^{\frac{3}{2}}\Pi \cp^3 v_{ttt}||_{0,\Gamma}^2 +\PP_0+\PP\int_0^t \PP\label{v_ttt^3}.
\end{align}
\begin{proof}
For any vector field $X$, the following identity allows one to compare $(\Pi \cp X)^3$ and $\cp X^3$:
\begin{equation}
(\Pi \cp X)^3 = \Pi_\lambda^3\cp X^{\lambda}=\cp X^3-g^{kl}\cp_k\eta^3\cp_l\eta_\lambda\cp X^{\lambda}.
\label{id comparsion}
\end{equation}
Invoking \eqref{id comparsion}, let $X=\rr^{\frac{3}{2}}\p_t^3 v$ and then taking $H^{-0.5}(\Gamma)$ norm yields
\begin{equation}
||\rr^{\frac{3}{2}}\cp\p_t^3 v^3||_{-0.5, \Gamma}^2 \lesssim ||\rr^{\frac{3}{2}}\Pi\cp\p_t^3 v||_{0,\Gamma}^2+||g^{kl}\cp_k\eta^3\cp_l\eta_\lambda||_{1.5,\Gamma}^2||\rr^{\frac{3}{2}}\p_t^3 v^{\lambda}||_{0.5,\Gamma}^2.
\end{equation}
We add $||\rr^{\frac{3}{2}}\p_t^2 v^3||_{-0.5, \Gamma}^2$ to both sides, use the fact that $||\rr^{\frac{3}{2}}\p_t^2 v^3||_{-0.5, \Gamma}^2+||\rr^{\frac{3}{2}}\cp\p_t^3 v^3||_{-0.5, \Gamma}^2$  is equivalent to $||\rr^{\frac{3}{2}}\p_t^2 v^3||_{0.5, \Gamma}^2$, invoke $\cp_k \eta^3 =\int_0^t\cp_k v^3$ (which is true since $\eta^3(0)=1$), to conclude \eqref{v_ttt^3}, where the term $||\rr^{\frac{3}{2}}\p_t^2 v^3||_{-0.5, \Gamma}^2$ on the right hand side is estimated using interpolation, Young's inequality, and the fundamental theorem of calculus.

Similarly, using \eqref{id comparsion} with $X=\rr \cp\p_t^2 v$ and $X=\rr \cp^2 \p_t v$, estimating in $H^{-0.5}(\Gamma)$ yields \eqref{v_tt^3} and \eqref{v_t^3}, respectively. Now, we need to control $||v^3||_{3.5,\Gamma}$. This cannot be controlled using the above method since $||\Pi \cp^4 v||_{0,\Gamma}^2$ is not part of the energy $E$. Nevertheless,  we recall the boundary condition
\begin{equation}
\sqrt{g}\Delta_g \eta^\alpha = \sqrt{g}g^{ij}\cp_{ij}^2 \eta^\alpha-\sqrt{g}g^{ij}\Gamma_{ij}^k\cp_k \eta^\alpha= -\sigma^{-1}A^{\mu \alpha}N_\mu q,\q \text{on}\,\,\Gamma
\label{bdy ful}
\end{equation}
where $\Gamma_{ij}^k = g^{kl}\cp_l\eta^{\mu}\cp_{ij}^2\eta_\mu$. Time differentiating \eqref{bdy ful} with $\alpha=3$ gives:
\begin{align}
 \sqrt{g}g^{ij}\cp_{ij}^2 v^3-\sqrt{g}g^{ij}\Gamma_{ij}^k\cp_k v^3=-\p_t(\sqrt{g}g^{ij})\cp_{ij}^2 \eta^3-\p_t(\sqrt{g}g^{ij}\Gamma_{ij}^k)\cp_k \eta^3\no\\
 -\sigma^{-1}\p_tA^{\mu 3}N_\mu q- \sigma^{-1}A^{\mu 3}N_\mu \p_t q
 \label{lap g v^3 = qt}
\end{align}
holds on $\Gamma$. Because $g^{ij}\in H^{2.5}(\Gamma)$ and $\Gamma_{ij}^k\in H^{1.5}(\Gamma)$, invoking the elliptic estimate for rough coefficients (see, e.g, Theorem 4 and Remark 2 in Milani \cite{milani1983regularity}), we obtain:
\begin{align}
||v^3||_{3.5,\Gamma}^2 \lesssim_M ||\p_t(\sqrt{g}g^{ij})\cp_{ij}^2 \eta^3||_{1.5,\Gamma}^2+||\p_t(\sqrt{g}g^{ij}\Gamma_{ij}^k)\cp_{k} \eta^3||_{1.5,\Gamma}^2\no\\
+||\p_t A^{\mu 3}N_\mu q||_{1.5, \Gamma}^2+||A^{\mu 3}N_\mu \p_tq||_{1.5,\Gamma}^2,
\end{align}
which can be controlled appropriately by the right hand side of \eqref{v^3}, where the last two terms can be controlled by with the help of Theorem \ref{esi ||q||_2 and ||q_t||_2 theo}.
\end{proof}

\subsection{Bounds for $v$, $R$ and their time derivatives}
Let $k=1,2,3$, commuting $\p_t^k$ to the second equation of \eqref{E}, we get
\begin{align}
\p^\alpha \p_t^k v_\alpha = (\delta^{\mu\alpha}-a^{\mu\alpha})\p_\mu \p_t^k v_\alpha-{\sum}_{\substack{j_1+j_2=k\\j_1\geq 1}}R^{-1}\p_t^{j_1}(Ra^{\mu\alpha})(\p_\mu \p_t^{j_2} v_\alpha)-R^{-1}\p_t^{k+1}R.\label{time derivatives density eq}
\end{align}
In addition, the first equation of \eqref{E} can be re-written as
\begin{align}
R'R\p_t v^\alpha+a^{\mu\alpha}\p_\mu R=0.
\end{align}
Commuting $\p_t^k$ to this equation and invoking \eqref{R_kk assumption}, we get
\begin{align}
\p^\alpha \p_t^k R = (\delta^{\mu\alpha}-a^{\mu\alpha})\p_\mu \p_t^k R - R'R\p_t^{k+1} v^\alpha\no\\
-{\sum}_{\substack{j_1+j_2=k\\j_1\geq 1}}[(\p_t^{j_1}a^{\nu\alpha})(\p_\mu \p_t^{j_2}R)+\big(\p_t^{j_1}(R'R)\big)(\p_t^{j_2+1}v^\alpha)].\label{time derivatives Euler eq}
\end{align}
When $k=3$, multiplying $(R')^{\frac{3}{2}}$ and then taking $L^2$ norm on both sides of \eqref{time derivatives density eq}, we get
\begin{align}
||(R')^{\frac{3}{2}}\p^\alpha \p_t^3 v_\alpha||_0 \leq \epsilon||(R')^{\frac{3}{2}}\p_t^3 v_\alpha||_1\no\\
+C{\sum}_{\substack{j_1+j_2=3\\j_1\geq 1}}||(R')^{\frac{3}{2}}\p_t^{j_1}(Ra^{\mu\alpha})(\p_\mu \p_t^{j_2} v_\alpha)||_0+ C||(R')^{\frac{3}{2}}R_{tttt}||_0,
\end{align}
where we have used Lemma \ref{prelim lemma a}(9)(10). The term 
$$
{\sum}_{\substack{j_1+j_2=3\\j_1\geq 1}}||(R')^{\frac{3}{2}}\p_t^{j_1}(Ra^{\mu\alpha})(\p_\mu \p_t^{j_2} v_\alpha)||_0
$$ 
is of lower order and can be controlled appropriately.  
Squaring and using Theorem \ref{energy estimate E}, we have
\begin{align}
||(\rr)^{\frac{3}{2}}\di v_{ttt}||_0^2\lesssim ||(R')^{\frac{3}{2}}\di v_{ttt}||_0^2 \lesssim \epsilon P(\mathcal{N})+\PP_0+\PP\int_0^t\PP. 
\end{align}
Now, in view of \eqref{v_ttt_1}, invoking \eqref{curl}, \eqref{v_ttt^3} and Theorem \ref{energy estimate E} gives
\begin{align}
||(\rr)^{\frac{3}{2}}v_{ttt}||_1^2 \lesssim \epsilon P(\mathcal{N})+\PP_0+\PP\int_0^t\PP. \label{control of v_ttt H^1}
\end{align}

We now move to estimate $||\rr R_{ttt}||_1$. Invoking \eqref{time derivatives Euler eq} for $k=3$, multiplying $R'$ on both sides and taking $L^2$ norm, we have:
\begin{align}
||R'R_{ttt}||_1 \lesssim \epsilon ||R'R_{ttt}||_1+||(R')^2v_{tttt}||_0+\epsilon \mathcal{N}+\PP_0+\PP\int_0^t \PP.
\end{align}
Here, $\epsilon \mathcal{N}$ appears when controlling the error term of \eqref{time derivatives Euler eq}\footnotemark.  \footnotetext{Specifically, $\epsilon \mathcal{N}$ is required to control $||R'(\p_t^3 a^{\nu\alpha})(\p_\nu R)||_0$. This term involves $||R'(\p v_{tt})(\p R)||_0$ at the top order, which is bounded by $||\rr \p v_{tt}||_1^{\frac{1}{2}}||\rr \p v_{tt}||_0^{\frac{1}{2}}||\p R||_1^{\frac{1}{2}}||\p R||_0^{\frac{1}{2}}\leq \epsilon(\sqrt{\mathcal{N}}+\mathcal{N})+\PP_0+\PP\int_0^t\PP$.  }
Squaring this provides:
\begin{align}
||\rr R_{ttt}||_1^2\lesssim||R'R_{ttt}||_1^2\lesssim \epsilon P(\mathcal{N})+ \PP_0+\PP\int_0^t \PP,\label{control of R'R_ttt H^1}
\end{align}
where Theorem \ref{energy estimate E} is also used.

Next, we estimate $||\mathfrak{R}_\kk \di v_{tt}||_1$. Invoking \eqref{time derivatives density eq} with $k=2$, multiplying $R'$ and then applying $H^1$ norm on both sides, we get
\begin{align}
||R' \p^\alpha \p_t^2 v_\alpha||_1 \leq \epsilon ||R'v_{tt}||_2+C {\sum}_{\substack{j_1+j_2=2\\j_1\geq 1}}||R'\p_t^{j_1}(Ra^{\mu\alpha})(\p_\mu \p_t^{j_2} v_\alpha)||_1+C||R'R_{ttt}||_1.
\end{align}
Using \eqref{control of R'R_ttt H^1}, squaring the above estimate leads to
\begin{align}
||\mathfrak{R}_\kk \di v_{tt}||_1^2\lesssim ||R'\di v_{tt}||_1^2 \lesssim  \epsilon P(\mathcal{N}) + \PP_0+\PP\int_0^t \PP.
\end{align}
In light of \eqref{v_tt_2}, the above bound for $||\rr \di v_{tt}||_1^2$,  together with \eqref{curl}, \eqref{v_tt^3} and Theorem \ref{energy estimate E} give
\begin{align}
||\rr v_{tt}||_2^2 \lesssim \epsilon P(\mathcal{N}) + \PP_0+\PP\int_0^t \PP.\label{control of R'v_tt H^2}
\end{align}

Furthermore, invoking \eqref{time derivatives Euler eq} for $k=2$, multiplying $\sqrt{R'}$ and taking $H^1$ norm and squaring, we get:
\begin{align}
||\sqrt{R'}R_{tt}||_2^2 \lesssim \epsilon ||\sqrt{R'}R_{tt}||_2^2+||(R')^{\frac{3}{2}} v_{ttt}||_1^2+\epsilon\mathcal{N}+\PP_0+\PP\int_0^t \PP,
\end{align}
which implies, after invoking \eqref{control of v_ttt H^1}, that
\begin{align}
||\sqrt{\rr}R_{tt}||_2^2 \lesssim \epsilon P(\mathcal{N})+\PP_0+\PP\int_0^t \PP \label{control of sqrt(R') R_tt H^2}. 
\end{align}

In addition, this allow us to continue this procedure to get an estimate for $R'\di v_t$; let $X=R'\p_t v$ and $s=3$ in \eqref{time derivatives density eq}, we get:
\begin{equation}
||R'\di v_t||_2 \lesssim \epsilon||R' v_t||_3+ ||R' R_{tt}||_2+\PP_0+\PP\int_0^t\PP,
\end{equation}
squaring, and invoking \eqref{curl}, \eqref{v_t^3} and \eqref{control of sqrt(R') R_tt H^2} gives
\begin{align}
||\rr v_{t}||_3^2 \lesssim \epsilon P(\mathcal{N}) + \PP_0+\PP\int_0^t \PP.\label{control of R'v_t H^3}
\end{align}
Now, invoking \eqref{time derivatives Euler eq} for $k=1$, squaring and taking $H^2$ norm yields
\begin{align}
||R_t||_3^2 \lesssim \epsilon ||R_t||_3^2 +||R' v_{tt}||_2+\epsilon\mathcal{N}+\PP_0+\PP\int_0^t \PP\lesssim \epsilon P(\mathcal{N})+ \PP_0+\PP\int_0^t\PP,\label{control of R_t H^3}
\end{align}
as a consequence of \eqref{control of R'v_tt H^2}. 

Finally, the above procedure yields
\begin{align}
||\di v||_3 \lesssim \epsilon ||v||_4 + ||R_t||_3,
\end{align}
and hence
\begin{align}
||v||_4^2 \lesssim \epsilon P(\mathcal{N})+ \PP_0+\PP\int_0^t\PP,
\end{align}
via \eqref{curl}, \eqref{v^3} and \eqref{control of R_t H^3}. Moreover, we have:
\begin{align}
||R||_4^2 \lesssim \epsilon ||R||_4^2 +||R' v_{t}||_3+\PP_0+\PP\int_0^t \PP\lesssim \epsilon P(\mathcal{N})+ \PP_0+\PP\int_0^t\PP,\label{control of R H^4}
\end{align}
via \eqref{control of R'v_t H^3}.
\subsection{The continuity argument, proof of Theorem \ref{main theorem 2} }
\paragraph*{Recovering the a priori assumptions:} We need to control the left hand side of \eqref{apriori v eta}-\eqref{a priori g gamma} by $\epsilon P(\mathcal{N})+\PP_0+\PP\int_0^t\PP$. The control for \eqref{apriori v eta} is a direct consequence of the Sobolev embedding , i.e., 
\begin{align}
||\p \eta||_{L^\infty} +||\p^2 \eta||_{L^\infty} \lesssim ||\eta||_4  \leq \PP_0+\PP\int_0^t\PP.
\end{align}
This also controls the left hand side of \eqref{a priori g gamma} by the definition of $g^{ij}$ and $\Gamma_{ij}^k$.

\paragraph*{Estimates at $t=0$:}
As we have seen that $\PP$ involves quantities involving time derivatives, and so one needs to show that these quantities can be controlled by $\PP_0$. More precisely, we show:
\begin{align}
||\rr v_t(0)||_3+||\rr v_{tt}(0)||_2+||(\rr)^{\frac{3}{2}}v_{ttt}(0)||_1+||R_t(0)||_3+||\sqrt{\rr}R_{tt}(0)||_2+||\rr R_{ttt}(0)||_1\no\\
||v_t(0)||_2+||\sqrt{\rr} v_{tt}(0)||_2+ ||\rr v_{ttt}(0)||_0 + ||R_{tt}(0)||_1+ ||\sqrt{\rr} R_{ttt}(0)||_0 \leq \PP_0.
\label{int norm at t=0}
\end{align}
This estimate is straightforward, i.e., we use \eqref{E mod}  to obtain $||\rr v_t(0)||_{3}\leq ||\rho_0^{-1} \p q(0)||_3 \lesssim \PP_0$. Moreover, we use \eqref{time derivatives density eq} with $k=0$ at $t=0$ to obtain $||R_t(0)||_3\leq ||\rho_0^{-1} \di v(0)||_3 \leq \PP_0$. The other quantities in \eqref{int norm at t=0} can be controlled similarly. In addition, we also need 
\begin{equation}
||\rr v_t(0)||_{3,\Gamma}+||\rr  v_{tt}(0)||_{2,\Gamma}+||(\rr)^{\frac{3}{2}} v_{ttt}(0)||_{1,\Gamma}\leq \PP_0.
\label{bdy norm at t=0}
\end{equation} 
To control $||\rr v_t(0)||_{3,\Gamma}$, we use \eqref{time derivatives Euler eq} to obtain $R' v_t^i(0) = -\delta^{ij}\p_j R(0)$, which implies $||\rr v_t^i(0)||_{3,\Gamma}\leq  ||R(0)||_{4,\Gamma}\leq \PP_0$. On the other hand, we control the normal component $v_t^3(0)$ using the elliptic estimate. Time differentiating \eqref{lap g v^3 = qt} and then restricting at $t=0$ yields:
\begin{align}
 \overline{\Delta} v_t^3(0)=- \sigma^{-1} q_{tt}(0)+ F,
 \end{align}
 where $F$ satisfies $||\rr F||_{1,\Gamma}\leq \PP_0$. From the elliptic theory, the control of $||\rr v_t^3(0)||_{3,\Gamma}$ requires the control of  $||\rr q_{tt}(0)||_{1,\Gamma}$ and hence $||R' q_{tt}(0)||_{1,\Gamma}$. Invoking the wave equation \eqref{wave equation}, this is bounded by $||\lap q(0)||_{1,\Gamma}+||\mathcal{F}_1||_{1,\Gamma}$. There is no problem to control $||\mathcal{F}_1||_{1,\Gamma}$ by $\PP_0$ in light of \eqref{F_r}. Furthermore, invoking the compatibility condition in section \ref{section 5}, i.e., $q_0=\sigma\overline{\Delta} \eta_0^3$, one  controls $||\lap q_0||_{1,\Gamma}$ by $||\eta_0^3||_{5.5}$. 
 
 The estimates for $||\rr  v_{tt}(0)||_{2,\Gamma},||(\rr)^{\frac{3}{2}} v_{ttt}(0)||_{1,\Gamma}$ are treated in a similar way, upon time differentiating more times and proceeding as above. We omit the details, but explain the estimates up to the highest order in an expository way.  First, to control the tangential component, we use \eqref{time derivatives Euler eq} and \eqref{time derivatives density eq} to get 
\begin{align}
R'v_{tt}^i(0) \sim \delta^{ij}\cp_j R_t(0) \sim \delta^{ij} \cp_j \p_\alpha v^\alpha(0),\\
(R')^{\frac{3}{2}} v_{ttt}^i(0) \sim \sqrt{R'} \delta^{ij} \cp_j R_{tt}(0)\sim \sqrt{R'}\delta^{ij} \cp_j \lap q_0\sim \sqrt{R'}\delta^{ij} \cp_j \lap \overline{\lap} \eta_0^3,
\end{align}
where $\sim$ means up to controllable terms. This yields that $$||\rr v_{tt}^i(0)||_{2,\Gamma}, \quad ||(\rr)^{\frac{3}{2}} v_{ttt}^i(0)||_{1,\Gamma}$$ are controlled by $||\di v_0||_{3, \Gamma}$ and $||\eta_0^3||_{6.5}$, respectively. 
 Second, to control the normal component, time-differentiating \eqref{lap g v^3 = qt} two times and restricting at $t=0$ yields $\overline{\lap} v_{tt}^3(0)\sim q_{ttt}(0)$. Therefore, from the elliptic theory, the control of $||\rr  v_{tt}^3(0)||_{2,\Gamma}$ requires that of  $||\rr q_{ttt}(0)||_{0,\Gamma}$ and hence  $||\lap q_t(0)||_{0,\Gamma}$, in light of the wave equation. Invoking the compatibility condition $q_t(0) \sim \overline{\lap} v_0^3$, $||\lap q_t(0)||_{0,\Gamma}$ is controlled by $||\lap \overline{\lap} v_0^3||_{0,\Gamma}$. On the other hand, time-differentiating \eqref{lap g v^3 = qt} three times and restricting at $t=0$ yields $\overline{\lap} v_{ttt}^3(0)\sim q_{tttt}(0)$. Therefore,  from the elliptic theory, the control of $||(\rr)^{\frac{3}{2}}  v_{ttt}^3(0)||_{1,\Gamma}$ requires that of  $||(\rr)^{\frac{3}{2}} q_{tttt}(0)||_{-1,\Gamma}$ and hence $||\sqrt{\rr}\lap q_{tt}(0)||_{-1,\Gamma}$.  Invoking the compatibility conditions $q_{tt}(0)\sim \overline{\lap}\p_3 q(0)$ and $q(0) \sim \overline{\lap} \eta_0^3$, we have that $||\sqrt{R'}\lap q_{tt}(0)||_{-1,\Gamma}$ is bounded by $||\eta_0^3||_{6.5}$.

Hence, Theorem \ref{close esi thm} implies 
\begin{equation}
\mathcal{N}(t)\lesssim \epsilon P(\mathcal{N}(t))+ P(\mathcal{N}(0))+P(\mathcal{N}(t))\int_0^t P(\mathcal{N}(s))\,ds.
\end{equation}
Invoking the continuity-boostrap argument in \cite{tao2006nonlinear}, this implies that there exists $\mathfrak{M}>0$ such that 
\begin{equation}
\mathcal{N}(t) \leq \mathfrak{M}, \q \text{whenever}\,\,t\in[0,T],
\label{N(t) bound}
\end{equation}
for some $T>0$. 
\subsection{Passing to the incompressible limit, proof of Theorem \ref{main theorem 1}}
\paragraph*{Proof for statement 1:}This is standard since we have an uniform a priori estimate.
\paragraph*{Proof for statement 2:}
The bound \eqref{N(t) bound} implies that $||v_\kk||_4+||R_\kk||_4\leq \sqrt{\mathfrak{M}}$ uniformly as $\kk\to\infty$. Therefore, by the Sobolev embedding, we have:
\begin{equation}
{\sum}_{\ell\leq 2}\Big( ||\p^\ell v_\kk||_{L^\infty(\Omega)}+ ||\p^\ell R_\kk||_{L^\infty(\Omega)} \Big)\leq \sqrt{\mathfrak{M}}.
\end{equation}
This yields that for each fixed $t\in [0,T]$, $v_\kk$ and $R_\kk$ are uniformly bounded and equicontinuous in $C^2(\Omega)$, which implies the convergence of $v_\kk$ and $R_\kk$ in $C^2(\Omega)$. Moreover, $v_\kk\rightarrow \mathfrak{v}$ since $a^{\mu\alpha}\p_\mu (v_{\kk})_\alpha\rightarrow 0$ in $L^\infty(\Omega)$, which is a consequence of $||\p_t q_{\kk}||_2$ being bounded independent of $\kk$ and $R'_\kk\to 0$ as $\kk\to \infty$.   


\section{The initial data} \label{section 5}
\subsection{The compatibility conditions}
\label{section 5.1}
The compatibility conditions for the initial data are necessary for construction of solutions, as well as for passing the solution to the incompressible limit. We recall that since 
\begin{align}
q = \sigma g^{ij}\n_\mu \cp_{ij}^2 \eta^\mu,\q \text{on}\,\, \Gamma,
\end{align}
we have:
\begin{align}
q|_{t=0} =\Big(\sigma g^{ij}\n_\mu \cp_{ij}^2 \eta^\mu\Big)\Big|_{t=0}:=H_0,\q \text{on}\,\, \Gamma, \label{cpt 0}
\end{align}
which is the zero-th order compatibility condition. In addition, for each $j\geq 1$, the $j$-th order compatibility reads 
\begin{align}
\p_t^j q|_{t=0} =\p_t^j \Big(\sigma g^{ij}\n_\mu \cp_{ij}^2 \eta^\mu\Big)\Big|_{t=0}:=H_j,\q \text{on}\,\, \Gamma. 
\label{cpt j}
\end{align}
Our goal is to construct $(\vv, \qq)$ that verifies the compatibility condition \eqref{cpt j} for $j=0,1,2,3$. We shall focus on the case when $\Omega=\mathbb{T}^2\times (0,1)$ , whose boundary $\Gamma$ is flat. Our method can easily be generalized to more general domains. 

\subsection{Formal construction}
We shall describe our method formally which serves as a good guideline for readers. Since 
$$
q \sim \overline{\lap}\eta^3, \q \text{on}\,\,\Gamma,
$$
we get
\begin{align*}
q_t \sim \overline{\lap} v^3, \q q_{tt} \sim \overline{\lap} v_t^3, \q
q_{ttt} \sim \overline{\lap} v_{tt}^3,\q \text{on}\,\,\Gamma,
\end{align*}
after taking time derivatives. Moreover, since the Euler equations imply
\begin{align*}
v_t \sim \p q, \q q_t \sim \kk \di v, 
\end{align*}
we have
\begin{align*}
q_{tt} \sim \overline{\lap} \p_3 q,\q q_{ttt}\sim \overline {\lap} \p_3 q_t\sim \kk\overline{\lap}\p_3 \di v, \q \text{on}\,\,\Gamma.
\end{align*}
For each $\ell=0,1,2,3$, we obtain the $\ell$-th order compatibility condition after restricting the above expression at $t=0$, i.e., 
\begin{align*}
q|_{t=0} \sim \overline{\lap}\eta_0^3, \q \text{on}\,\,\Gamma,\\
q_t|_{t=0} \sim \overline{\lap} v_0^3, \q \text{on}\,\,\Gamma,\\
q_{tt}|_{t=0} \sim \overline{\lap} \p_3 q_0, \q \text{on}\,\,\Gamma, \\
q_{ttt}|_{t=0} \sim \kk\overline{\lap}\p_3 \di v, \q \text{on}\,\,\Gamma.
\end{align*} 
On the other hand, since 
\begin{align}
q_t \sim \kk \di v, \q q_{tt}\sim \kk \di v_t\sim \kk \lap q, \q q_{ttt} \sim \kk \lap q_t \sim \kk^2 \lap \di v,
\end{align}
then
\begin{align}
q_0\sim \overline{\lap} \eta_0^3,\q \text{on}\,\,\Gamma,\label{cpt 0'}\\
\di v_0 \sim \kk^{-1} \overline{\lap} v_0^3,\q \text{on}\,\,\Gamma,\label{cpt 1'}\\
\lap q_0 \sim \kk^{-1}\overline{\lap} \p_3 q_0,\q \text{on}\,\,\Gamma,\label{cpt 2'}\\
\lap \di v_0 \sim \kk^{-1}\overline{\lap} \p_3 \di v_0, \q \text{on}\,\,\Gamma\label{cpt 3'}.
\end{align}
In other words, the first order compatibility condition (i.e., \eqref{cpt j} when $j=1$), is expressed in $v_0$, and the second order compatibility condition is expressed in $q_0$, and finally the third order compatibility condition is expressed in $v_0$ again.  

To construct initial data that satisfies the compatibility conditions up to order $3$, our first step is to obtain $(\uu, \pp)$ that satisfies the \eqref{cpt 0'}. This is easy, since we can simply let $\uu$ to be velocity for the incompressible case, i.e., $\uu = \mathfrak{u}_0$, and $\pp$ 
\begin{align}
-\lap\pp = (\p_\mu \uu^\nu)(\p_\nu \uu^\mu),\q\text{in}\,\,\Omega,\\
 \pp=\overline{\lap}\eta_0^3,\q\text{on}\,\,\Gamma.
\label{p_0 construction}
\end{align}  
Our next step is to construct a velocity vector field $\ww=(\ww^1, \ww^2, \ww^3)$ that satisfies \eqref{cpt 1'}. To achieve this, we set $\ww^1=\uu^1$ and $\ww^2=\uu^2$, while we define $\ww^3$ via solving
\begin{align}
\lap^2 \ww^3=\lap^2 \uu^3,\q \text{in}\,\,\Omega,\\
 \ww^3=\uu^3 ,\q \p_3 \ww^3 \sim \kk^{-1}\overline{\lap}\uu^3-\p_1\uu^1-\p_2\uu^2,\q\text{on}\,\,\Gamma.
 \label{w_0 construction}
\end{align}
We now construct $\qq$ that satisfies \eqref{cpt 2'}. We define $\qq$ by the solution of
\begin{align}
\lap^3 \qq = 0,\q \text{in}\,\,\Omega, \no\\
\qq =\pp,\q \p_3 \qq=\p_3\pp, \q \lap \qq \sim \kk^{-1}\overline{\lap}\p_3 \pp,\q \text{on}\,\,\Gamma.
\label{q_0 construction}
\end{align}
Finally, we need to construct $\vv$ using \eqref{cpt 3'}. To achieve this, we set $\vv^1=\uu^1$, $\vv^2=\uu^2$, and  we define $\vv^3$ by solving 
\begin{align}
\lap^4 \vv^3 = \lap^4 \ww^3,\q \text{in}\,\,\Omega,\no\\
\vv^3=\ww^3,\q \p_3\vv^3\sim \kk^{-1}\overline{\lap}\ww^3-\p_1\ww^1-\p_2\ww^2, \q \text{on}\,\,\Gamma,\no\\
 \p_3^2 \vv^3 \sim \kk^{-1}\p_3\overline{\lap}\ww^3-\p_3\p_1\ww^1-\p_3\p_2\ww^2,\q \text{on}\,\,\Gamma,\no\\
\lap \p_3 \vv^3\sim \kk^{-1}\overline{\lap}\p_3 \di \ww-\lap \p_1 \ww^1-\lap \p_2\ww^2,\q\text{on}\,\,\Gamma.
\label{v_0 construction}
\end{align}
\rmk In fact, $\overline{\lap}\eta_0^3=0$ on the boundary of the reference domain $\mathbb{T}^2\times(0,1)$.  But that we do not use this condition exactly because we want to keep the regularity of each argument as it should hold for the general domain.
\thm  \label{data construction}
Let $\mathfrak{u}_0\in H^{6.5}(\Omega)$ be a divergence free vector field in $\Omega$ and $\pp$ be the associated pressure. Then there exists initial data $(\vv, \qq)=(\vv^\kk, \qq^\kk)$ satisfying the compatibility conditions up to order 3, i.e.,  \eqref{cpt 0'}-\eqref{cpt 3'}, such that $\vv^\kk\rightarrow \uu$ in $C^2(\Omega)$ and $\di \vv^\kk \rightarrow 0$ in $C^1(\Omega)$ as $\kk\rightarrow \infty$, and $\PP_0$ is uniformly bounded for all $\kk$.

\begin{proof}
$(\vv, \qq)$ verifies \eqref{cpt 0'}-\eqref{cpt 3'} follows automatically from our construction. Since $\pp$ satisfies the elliptic equation \eqref{p_0 construction}, for $s\geq 4$, we have:
\begin{align}
||\pp||_s \lesssim ||\lap \pp||_{s-2}+||\pp||_{s-0.5, \Gamma},
\label{est p}
\end{align}
which requires $||\uu||_{s-1}$ and $||\eta_0||_{s+2}$ to control. 
Moreover, by the poly-harmonic estimate applied to \eqref{q_0 construction} we have:
\begin{align}
||\qq||_s \lesssim ||\lap \qq||_{s-2.5, \Gamma}+||\p_3 \qq||_{s-1.5,\Gamma}+||\qq||_{s-0.5,\Gamma}\\
\lesssim \kk^{-1}||\overline{\lap}\p_3 \pp||_{s-2}+||\p_3 \pp||_{s-1}+||\pp||_s.
\end{align}
Invoking \eqref{est p}, this requires  $||\uu||_{s}$ and $||\eta_0||_{s+3}$ to control. On the other hand, invoking \eqref{w_0 construction} and the poly-harmonic estimate, we get:
\begin{align}
||\ww^3||_s \lesssim ||\lap^2\uu^3||_{s-4}+||\p_3 \ww^3||_{s-1.5,\Gamma}+||\ww^3||_{s-0.5,\Gamma}\\
\lesssim ||\lap^2\uu^3||_{s-4}+\kk^{-1}||\overline{\lap}\uu^3||_{s-1}+||\p_1\ww^1||_{s-1}+||\p_2\ww^2||_{s-1}+||\uu^3||_s,
\end{align}
which needs $||\uu^3||_{s+1}$ to control. In addition, since  $\ww^{i}=\uu^i$, one controls $||\ww||_s$ via $||\uu||_{s+1}$. Moreover, invoking \eqref{v_0 construction} and the poly-harmonic estimate, we get:
\begin{align}
||\vv^{3}||_s \lesssim ||\lap^4 \uu^3||_{s-8}+ ||\lap \p_3 \vv^{3}||_{s-3.5,\Gamma}+||\p_3^2 \vv^{3}||_{s-2.5,\Gamma}+||\p_3 \vv^{3}||_{s-1.5,\Gamma}+||\vv^{3}||_{s-0.5,\Gamma}\\
\lesssim||\lap^4 \uu^3||_{s-8}+\kk^{-1}||\overline{\lap}\p_3\di \ww||_{s-3}+||\lap\p_1\ww^1||_{s-3}+||\lap\p_2\ww^2||_{s-3}\\
+\kk^{-1}||\p_3 \overline{\lap}\ww^3||_{s-2}+||\p_3\p_1 \ww^1||_{s-2}+||\p_3\p_2 \ww^2||_{s-2}\\
+\kk^{-1}||\overline{\lap}\ww^3||_{s-1}+||\p_1\ww^1||_{s-1}+||\p_2\ww^2||_{s-1}+||\ww^3||_s ,
\end{align}
which requires $||\ww^3||_{s+1}$ and hence $||\uu^3||_{s+2}$ to control. Once again, since  $\vv^{i}=\uu^i$, one controls $||\vv||_s$ through $||\uu||_{s+2}$. 

 Next, since \eqref{w_0 construction} implies 
\begin{align}
\lap^2 (\ww^3-\uu^3)=0,\q \text{in}\,\,\Omega,\\
 \ww^3-\uu^3=0 ,\q \p_3 (\ww^3-\uu^3) \sim \kk^{-1}\overline{\lap}\uu^3,\q\text{on}\,\,\Gamma,
\end{align}  
we have that $||\ww^3-\uu^3||_s\rightarrow 0$ as $\kk\to \infty$, and hence $\ww\to \uu$ in $H^s(\Omega)$ as $\kk\to \infty$. Similarly, \eqref{v_0 construction} implies $\vv\to \ww$ in $H^s(\Omega)$ as $\kk\to\infty$, and so we conclude that $\vv\to\uu$ in $H^s(\Omega)$ as $\kk\to\infty$.  Furthermore, because $s\geq 4$ and $\vv$ is uniformly bounded in $H^s$,  we have that $\vv\to \uu$ in $C^2(\Omega)$ thanks to Arzel\`a-Ascoli and $\di \vv \rightarrow \di \uu =0$ in $C^1(\Omega)$.

Finally,  we recall that $\PP_0$ consists $$||\vv||_4, ||\vv||_{4,\Gamma}, ||\qq||_4, ||\qq||_{4,\Gamma}, ||\di \vv|_{\Gamma}||_{3,\Gamma}, ||\lap \vv|_\Gamma||_{2,\Gamma}),$$
which can all be controlled by $||\uu||_{s+2}=||\mathfrak{u}_0||_{s+2}$ and $||\eta_0||_{s+3}$ with $s=4.5$. 
\end{proof}

\rmk 
The initial data constructed in Theorem \ref{data construction} is given in terms of the initial pressure $\qq$ instead of the initial density $R_0$.  This is because 
the boundary condition is more easily stated in terms of $q$ and 
we need to make sure that the quantities $||\qq||_4$ and $||\qq||_{4,\Gamma}$ are bounded uniformly in $\kk$. But we can compute $R_0$ through the equation of states $R=R(q)$, i.e.,  $R_{0}=[(c_\gamma \kk)^{-1} \qq+\beta]^{1/\gamma}$.

The rest of this section is devoted to provide detailed construction, and for the sake of simple expositions, we assume the equation of state is taken to be
$$
q(R) = \kk (R-1). 
$$
This allows us to exchange $q$ and $R$ in an explicit way. Also, throughout the rest of this section, we shall use $Q$ to denote a rational function.  
\subsection{Construction for $(\uu, \pp, \Omega)$ that satisfies (\ref{cpt 0}) while $j=0$}
Let $\uu = \mathfrak{v}_0$, where $\mathfrak{v}_0$ is the data for the incompressible Euler equations. 
Since 
$
H_0 = \sigma \overline{\lap}\eta_0^3,
$
 we define $\pp$ by solving
\begin{align}
\begin{cases}
-\lap\pp = (\p_\mu \uu^\nu)(\p_\nu \uu^\mu),\q\text{in}\,\,\Omega,\\
 \pp=H_0,\q\text{on}\,\,\Gamma.
\end{cases}
\end{align}

\subsection{Construction for $\ww$ that satisfies (\ref{cpt j}) while $j=1$ }
We next consider the first order compatibility condition, i.e., $\p_t q|_{t=0} = H_1$. Since 
\begin{align}
\p_t \Big(\sigma g^{ij}\n_\mu \cp_{ij}^2 \eta^\mu\Big)= \sigma g^{ij}\n_\mu \cp_{ij}^2 v^\mu+\sigma Q( \n, \cp \eta, \cp v)\cp^2 \eta,
\label{one time derivative}
\end{align}
 and thus
\begin{align}
H_1 = \sigma \overline{\lap} v_0^3+\sigma Q(\cp \eta_0, \cp v_0)\cp^2\eta_0. 
\label{H1}
\end{align}
On the other hand, since $\p_t q=-R\kk a^{\mu\alpha}\p_\mu v_\alpha$, \eqref{cpt j} with $j=1$ becomes:
\begin{equation}
\di v_0 = \kk^{-1}(\kk^{-1}q_0+1)H_1,\q \text{on}\,\,\Gamma,
\end{equation}
and so
\begin{align*}
\p_3 v_0^3 = \kk^{-1}(\kk^{-1}q_0+1)H_1 - \p_1 v_0^1-\p_2 v_0^2,\q \text{on}\,\,\Gamma.
\end{align*}
   Furthermore, this suggests that $\ww$ should be constructed as follows:
let $\ww = (\textbf{u}_0^1, \textbf{u}_0^2, \textbf{w}_0^3)$, where $\ww^3$ solves
\begin{align}
\begin{cases}
\lap^2 \ww^3 = \lap^2 \uu^3,\q \text{in}\,\,\Omega,\\
\ww^3 = \uu^3,\q \text{on}\,\,\Gamma,\\
\p_3\ww^3 = \kk^{-1}\sigma(\kk^{-1}\pp+1)\overline{\lap}\uu^3-\kk^{-1}\sigma(\kk^{-1}\pp+1) Q(\cp \eta_0, \cp \uu)\cp^2\eta_0-\p_1\uu^1-\p_2\uu^2,\q \text{on}\,\, \Gamma.
\end{cases}
\label{ww eq}
\end{align}

\subsection{Construction for $\qq$ that satisfies (\ref{cpt j}) while $j=2$}
The second order compatibility condition reads  $\p_t^2 q|_{t=0} = H_2$, and we need to express this in terms of $\eta_0, v_0$ and $q_0$, which yields a system satisfied by $\pp$.  Invoking \eqref{one time derivative}, we have
\begin{align}
\p_t^2 \Big(\sigma g^{ij}\n_\mu \cp_{ij}^2 \eta^\mu\Big)= \sigma g^{ij}\n_\mu \cp_{ij}^2 v_t^\mu+\sigma Q(\n, \cp \eta, \cp v)\cp^2 v 
+\sigma Q(\n, \cp \eta, \cp v)\cp^2 \eta(\cp v_t+1)  \label{two time derivatives sub}.
\end{align}
In addition, since $Rv_t^\mu + a^{\nu\mu}\p_\nu q=0$, we get for $s=1,2$ that
\begin{align}
\cp^s(v_t^\mu)= -R^{-1}a^{\nu\mu} \cp^s\p_\nu q- \sum_{1\leq k \leq s}\cp^k (R^{-1}a^{\nu\mu})\cp^{s-k}\p_\nu q.
\end{align}
This, together with \eqref{two time derivatives sub} and the equation of state $R=\kk^{-1}q+1$ imply
\begin{align}
H_2 = H_2(\eta_0, p_0, v_0)= -\sigma (\kk^{-1}q_0+1)^{-1}\overline{\lap}\p_3 q_0+\sigma Q(\cp \eta_0, \cp v_0)\cp^2 v_0\no\\
-\sigma Q\Big((\kk^{-1}q_0+1)^{-1}, \kk^{-1}\cp q_0, \p\eta_0, \cp\p \eta_0\Big)\cp \p q_0\\
-\sigma Q\Big((\kk^{-1}q_0+1)^{-1}, \p\eta_0, \kk^{-1}\p q_0,\kk^{-1} \cp^2 q_0\Big) \p q_0\no\\
+\sigma Q\Big((\kk^{-1}q_0+1)^{-1}, \cp \eta_0, \cp\p\eta_0, \cp v_0,  \p q_0 \Big)(\cp \p q_0+\cp^2 \eta_0).\label{H2}
\end{align}
On the other hand,  the continuity equation implies $R a^{\mu\alpha}\p_\mu v_\alpha = -\kk^{-1} \p_t q$, and hence
\begin{align}
-\kk^{-1} \p_t^2 q = \p_t(R a^{\mu\alpha})\p_\mu v_\alpha+R a^{\mu\alpha}\p_\mu \p_t v_\alpha
= \p_t(R a^{\mu\alpha})\p_\mu v_\alpha- R a^{\mu\alpha} \p_\mu (R^{-1}a^\nu_{\,\,\alpha} \p_\nu q)\\
=-a^{\mu\alpha}a^{\nu}_{\,\,\alpha} \p_\mu \p_\nu q-Ra^{\mu\alpha}\p_\mu (R^{-1}a^\nu_{\,\,\alpha})\p_\nu q+ \p_t(R a^{\mu\alpha})\p_\mu v_\alpha.
\label{q_tt pre}
\end{align}
Restricting the above identity to the boundary $\Gamma$ and then taking $t=0$, we get
\begin{align}
\kk^{-1}\p_t^2 q|_{t=0}=\lap q_0- Q\Big((\kk^{-1}q_0+1)^{-1}, \p \eta_0, \p^2\eta_0, \cp v_0,  \kk^{-1}\p q_0\Big)\p q_0\no\\
+Q(\kk^{-1}q_0,\p \eta_0, \p v_0 )\p v_0.
\label{q_tt}
\end{align}
Invoking \eqref{H2} and \eqref{q_tt}, we are able to rewrite \eqref{cpt j} when $j=2$ as
\begin{align}
\lap q_0 = Q\Big((\kk^{-1}q_0+1)^{-1}, \p \eta_0, \p^2\eta_0, \cp v_0,  \kk^{-1}\p q_0\Big)\p q_0\\
-Q(\kk^{-1}q_0,\p \eta_0, \p v_0 )\p v_0+\kk^{-1}H_2(\eta_0, p_0, v_0).
\label{j=2 explicit}
\end{align} 
This yields that $\qq$ should solve:
\begin{align}
\begin{cases}
\lap^3 \qq = 0,\q \text{in}\,\,\Omega,\\
\qq = \pp, \q\text{on}\,\, \Gamma,\\ 
\frac{\p \qq}{\p N}=\p_3 \qq = \p_3 \pp=\frac{\p \pp}{\p N},\q\text{on}\,\,\Gamma,\\
\lap \qq = \varphi, \q\text{on}\,\, \Gamma.
\end{cases}
\end{align}
Here, 
\begin{align*}
\varphi = Q\Big((\kk^{-1}\pp+1)^{-1}, \p \eta_0, \p^2\eta_0, \cp \ww,  \kk^{-1}\p \pp\Big)\p \pp
-Q(\kk^{-1}\pp,\p \eta_0, \p \ww )\p \ww 
+\kk^{-1}H_2(\eta_0, \pp, \ww),
\end{align*}
which is obtained from \eqref{j=2 explicit}.
\subsection{Construction for $\vv$ that satisfies (\ref{cpt j}) while $j=3$}
Our last step is to construct $\vv$ that satisfies third order compatibility condition, i.e., $\p_t^3 q|_{t=0}=H_3$ on $\Gamma$. Similar to what has been done for the previous cases when $j=0,1,2$, we shall first compute the compatibility condition explicitly. Invoking \eqref{two time derivatives sub}, as well as $v_t^\mu=-R^{-1}a^{\nu\mu } \p_\mu q $ and $\p_t q = -\kk R a^{\mu\alpha}\p_\mu v_\alpha$, we have
\begin{align}
\p_t^3 \Big(\sigma(g^{ij}\n_\mu \cp_{ij}^2 \eta^\mu\Big)= -\p_t\Big(\sigma(g^{ij}\n_\mu \cp_{ij}^2(R^{-1}a^{\nu\mu} \p_\nu q)\Big)\no\\
+\sigma Q( \n, \cp \eta, \cp v)\cp^2 v_t+\sigma Q(g, \n, \cp \eta, \cp v, \cp v_t)\cp^2 v
+\sigma Q(\n, \cp \eta, \cp v, \cp v_t)(\cp v_{tt}+\cp^2 \eta)\no\\
= -\sigma g^{ij}\n_\mu \cp_{ij}^2\p_t(R^{-1}a^{\nu\mu} \p_\nu q)+\sigma Q( \n, \cp \eta, \cp v)\cp_{ij}^2 (R^{-1}a^{\nu\mu} \p_\nu q)\no\\
+\sigma Q(\n, R^{-1}, \cp R^{-1}, \cp^2 R^{-1}, \cp v, \p \eta, \cp\p\eta, \cp^2\p\eta, \p q, \cp\p q)\cp^2\p q\no\\
+\sigma Q(\n, R^{-1}, \cp R^{-1},\cp v, \p \eta, a, \cp \p \eta,\p q, \cp \p q)\cp^2 v\no\\
+\sigma \kk Q(\n, R^{-1}, \p R^{-1},\p v, \p^2 v, \p \eta,  \cp\p\eta,\p q, \cp \p q)(a^{\mu\alpha}\p_\mu \cp \p v_\alpha+\cp^2 \eta),
\label{H3 1}
\end{align}
where
\begin{align}
\sigma g^{ij}\n_\mu \cp_{ij}^2\p_t(R^{-1}a^{\nu\mu} \p_\nu q)= \sigma g^{ij}\n_\mu R^{-1}a^{\nu\mu} \p_\nu  \cp_{ij}^2 q_t\no\\
+\sigma Q(\n , R, \cp R, \cp^2 R, \p v, \p \eta, \cp\p\eta, \cp^2\p\eta, \p q, \cp \p q, \cp^2 \p q)(\p q_t+ \cp \p q_t)\no\\
= -\kk \sigma g^{ij}\n_\mu a^{\nu\mu} \p_\nu  \cp_{ij}^2 (a^{\alpha\beta} \p_\alpha v_\beta)
+ \kk \sigma {\sum}_{k=1,2,3}Q(\n ,\p^k R)\Big(\p^{3-k} (a^{\alpha\beta}\p_\alpha v_\beta)\Big)\no\\
+\kk \sigma Q(\n , R, \p R, \p^2 R,\p v, \p^2 v, \p^3 v, \p \eta, \cp\p\eta, \cp^2\p\eta, \p q, \cp \p q, \cp^2 \p q)\Big( a^{\alpha\beta}\p_\alpha (\cp \p v_\beta+\p v_\beta)\Big).
\label{H3 2}
\end{align}
Restricting \eqref{H3 1} and \eqref{H3 2} at $t=0$, we get
\begin{align}
H_3 = H_3(\eta_0, q_0, v_0) = -\kk \sigma \p_3 \overline{\lap} \di v_0-\sigma \sum_{\ell=1,2,3}(\p^\ell q_0 )(\p^{3-\ell} \di v_0 )\no\\
+\kk\sigma Q\Big((\kk^{-1}q_0+1)^{-1}, \p v_0, \p^2 v_0, \p^3 v_0, \p \eta_0, \p^2 \eta_0, \p^3\eta_0, \p q_0, \p^2 q_0, \cp^2 \p q_0  \Big)\sum_{\ell=1,2}\p^\ell \di v_0.
\label{H3}
\end{align}
Next, invoking \eqref{q_tt pre}, we obtain 
\begin{align}
\kk^{-1} q_{ttt} 
=\p_t\Big(a^{\mu\alpha}a^{\nu}_{\,\,\alpha} \p_\mu \p_\nu q+Ra^{\mu\alpha}\p_\mu (R^{-1}a^\nu_{\,\,\alpha})\p_\nu q- \p_t(R a^{\mu\alpha})\p_\mu v_\alpha\Big)\no\\
=a^{\mu\alpha}a^{\nu}_{\,\,\alpha} \p_\mu \p_\nu q_t+Ra^{\mu\alpha}\p_\mu (R^{-1}a^\nu_{\,\,\alpha})\p_\nu q_t
+Q(R, R^{-1}, \p R^{-1}, \p \eta, \cp\p\eta, v, \p v)\p^2 q \\
=-R\kk a^{\mu\alpha}a^{\nu}_{\,\,\alpha} \p_\mu \p_\nu (a^{\beta\gamma}\p_\beta v_\gamma)
-2\kk a^{\mu\alpha}a^{\nu}_{\,\,\alpha} (\p_\mu R) \p_\nu (a^{\beta\gamma}\p_\beta v_\gamma)
\\-\kk a^{\mu\alpha}a^{\nu}_{\,\,\alpha} (\p_\mu \p_\nu R) (a^{\beta\gamma}\p_\beta v_\gamma)
+Q(R, R^{-1}, \p R^{-1},  \p \eta, \cp\p\eta)\p (a^{\beta\gamma}\p_\beta v_\gamma)\\
+Q(R,\p R, R^{-1}, \p R^{-1},  \p \eta, \cp\p\eta)a^{\beta\gamma}\p_\beta v_\gamma
+Q(R, R^{-1}, \p R^{-1}, \p \eta, \cp\p\eta, v, \p v)\p^2 q.
\label{q_ttt}
\end{align}
Restricting \eqref{q_ttt} to the boundary $\Gamma$ and then taking $t=0$, we have
\begin{align}
\kk^{-1} q_{ttt}|_{t=0}= -\kk R_0 \lap \di v_0-{\sum}_{\ell=1,2}2(\p^\ell q_0)(\p^{2-\ell} \di v_0)\no\\
+Q\Big((\kk^{-1} q_0+1)^{-1}, \kk^{-1}q_0, \p v_0, \p \eta_0, \p^2 \eta_0\Big){\sum}_{\ell=0,1}\p^\ell \di v_0\no\\
+Q\Big((\kk^{-1} q_0+1)^{-1}, \kk^{-1}q_0, v_0, \p v_0,  \p \eta_0, \p^2 \eta_0\Big)\p^2 q_0.
\end{align}
Invoking \eqref{H3}, the compatibility condition $q_{ttt}|_{t=0}=H_3$ can then be re-expressed as
\begin{align}
 \lap \di v_0= \psi(\eta_0, q_0, v_0)
\end{align}
where
\begin{align}
\psi (\eta_0, q_0, v_0) = -\kk^{-1}(\kk^{-1}q_0+1){\sum}_{\ell=1,2}2(\p^\ell q_0)(\p^{2-\ell} \di v_0)\no\\
+\kk^{-1}Q\Big((\kk^{-1} q_0+1)^{-1}, \kk^{-1}q_0, \p v_0, \p \eta_0, \p^2 \eta_0\Big){\sum}_{\ell=0,1}\p^\ell \di v_0\no\\
+\kk^{-1}Q\Big((\kk^{-1} q_0+1)^{-1}, \kk^{-1}q_0, v_0, \p v_0,  \p \eta_0, \p^2 \eta_0\Big)\p^2 q_0-\kk^{-2}(\kk^{-1} q_0+1)^{-1}H_3(\eta_0, q_0, v_0).
\end{align}
This implies that $\vv = (\vv^1, \vv^2, \vv^3)$ should be constructed such that $\vv^1=\uu^1$ and $\vv^2=\uu^2$, whereas $\vv^3$ solves
\begin{align}
\begin{cases}
\lap^4 \vv^3 = \lap^4\ww^3,\q \text{in}\,\,\Omega,\\
\vv^3 = \ww^3,\q \text{on}\,\,\Gamma,\\
\p_3 \vv^3 =  \kk^{-1}\sigma(\kk^{-1}\qq+1)\overline{\lap}\ww^3-\kk^{-1}\sigma(\kk^{-1}\qq+1) Q(\cp \eta_0, \cp \ww)\cp^2\eta_0-\p_1\ww^1-\p_2\ww^2,\q \text{on}\,\,\Gamma,\\
\p_3^2 \vv^3 =\p_3\Big(\kk^{-1}\sigma(\kk^{-1}\qq+1)\overline{\lap}\ww^3-\kk^{-1}\sigma(\kk^{-1}\qq+1) Q(\cp \eta_0, \cp \ww)\cp^2\eta_0-\p_1\ww^1-\p_2\ww^2\Big),\q \text{on}\,\,\Gamma,\\
\lap \p_3 \vv^3=\psi(\eta_0, \qq, \uu)-\lap\p_1\ww^1-\lap\p_2\ww^2,\q \text{on}\,\,\Gamma.
\end{cases}
\end{align}

\begin{appendix}
\section*{Appendix}
\section{Basic estimates}
\thm (Standard div-curl estimates )
Let $X$ be a vector field on $\Omega$ with sufficiently regular boundary $\Gamma$. Define $\di X=\p_j X^j$ and $(\curl X)_{ij}=\p_i X_j-\p_j X_i$, then for $1\leq s\leq 4$, we have
\begin{align}
||X||_s \lesssim ||\di X||_{s-1}+||\curl X||_{s-1}+||X\cdot N||_{s-0.5,\Gamma}+||X||_0,\label{div curl normal}\\
||X||_s \lesssim ||\di X||_{s-1}+||\curl X||_{s-1}+||X\cdot \mathcal{T}||_{s-0.5,\Gamma}+||X||_0,\label{div curl tang}
\end{align}
where $N$ is the outward unit normal to $\Gamma$, whereas $\mathcal{T}$ is the unit vector which is tangent to $\Gamma$. 

\begin{proof}
We refer \cite{lindblad2009priori} for the detailed proof. 
\end{proof}

\section{The energy identity for the wave equations of order $3$}
\label{section B'}
 We recall that for $r=1,2,3$, the wave equation reads:
\begin{align}
JR'\p_t^{r+1} q - a^{\nu\alpha}A^\mu_{\,\,\alpha}\p_\nu\p_\mu \p_t^{r-1} q = \mathcal{G}_r+\mathcal{S}_r, 
\end{align}
where
\begin{align}
\mathcal{G}_r=
-{\sum}_{\substack{j_1+j_2=r\\j_1\geq 1}}\big(\p_t^{j_1}(JR')\big)(\p_t^{j_2+1}q)+a^{\nu\alpha}(\p_\nu \rho_0)\p_t^r v_\alpha\no\\
+{\sum}_{\substack{j_1+j_2=r-1\\j_1\geq 1}}a^{\nu\alpha}\p_\nu(\p_t^{j_1}A^{\mu}_{\,\,\alpha}\cdot\p_\mu \p_t^{j_2}q)
-\rho_0{\sum}_{j_1+j_2=r-1}(\p_t^{j_1+1} a^{\nu\alpha})(\p_t^{j_2}\p_\nu v_\alpha). 
\end{align}
and
\begin{align}
\mathcal{S}_r=a^{\nu\alpha}(\p_\nu A^\mu_{\,\,\alpha})\p_\mu \p_t^{r-1}q.
\end{align}
\thm  For $r=1,2,3$, let
\label{est W thm 1}
\begin{align}
W_r^2= \frac{1}{2} \int_\Omega \rho_0^{-1}(JR'\p_t^r q)^2\,dy+\frac{1}{2}\int_\Omega \rho_0^{-1}R'(A^{\nu\alpha}\p_\nu \p_t^{r-1}q)(A^\mu_{\,\,\alpha}\p_\mu \p_t^{r-1}q)\,dy\no\\
+\frac{\sigma}{2}\int_{\Gamma} \rr \sqrt{g}g^{ij}\Pi_\mu^\alpha(\cp_i\p_t^r \eta^\mu)(\cp_j \p_t^r \eta_\alpha)\,dS. \label{W 1}
\end{align}
Then, 
\begin{equation}
{\sum}_{r\leq 3}W_r^2 \leq  \epsilon P(\mathcal{N}) +\epsilon (||q||_2^2+||q_t||_2^2) +\PP_0 +\PP\int_0^t\PP, \q t\in[0,T], \label{est W 1}
\end{equation}
where $T>0$ is sufficiently small. 

\paragraph*{Proof of Theorem \ref{est W thm 1}}
It suffices to consider the case when $r=3$. Invoking \eqref{p_t J} and \eqref{R_kk assumption}, we have:
\begin{align}
\frac{d}{dt}\frac{1}{2}\int_\Omega \rho_0^{-1}(JR'\p_t^3 q)^2 \,dy = \int_\Omega \rho_0^{-1}(JR'\p_t^3 q)(a^{\nu\alpha}A^\mu_{\,\,\alpha}\p_\nu\p_\mu \p_t^{2} q) \,dy\no\\
+\int_\Omega \rho_0^{-1}(JR'\p_t^3 q)(\mathcal{G}_3+\mathcal{S}_3)\,dy+\mathcal{R},\label{3.13}
\end{align} 
where $\mathcal{R}$ consists of error terms that are generated when $\p_t$ falls on either $J$ or $R'$, which we have no problem to control. In addition, 
\begin{align}
\int_\Omega \rho_0^{-1}(JR'\p_t^3 q)(a^{\nu\alpha}A^\mu_{\,\,\alpha}\p_\nu\p_\mu \p_t^{2} q) \,dy \no\\
= \int_\Omega \rho_0^{-1}(R'\p_t^3 q)(A^{\nu\alpha}\p_\nu)(A^\mu_{\,\,\alpha}\p_\mu \p_t^{2}q)\,dy
-\int_\Omega \rho_0^{-1}(JR'\p_t^3 q)\mathcal{S}_3. \label{3.14}
\end{align}
The last term in \eqref{3.14} cancels with the corresponding term in \eqref{3.13}, which is essential since $||\mathcal{S}_3||_0$ cannot be controlled uniformly when $R'\to 0$. Moreover, the first term on the right hand side of \eqref{3.14} is treated as:
\begin{align}
\int_\Omega \rho_0^{-1}(R'\p_t^3 q)(A^{\nu\alpha}\p_\nu)(A^\mu_{\,\,\alpha}\p_\mu \p_t^{2}q)\,dy= -\int_\Omega \rho_0^{-1} R'(A^{\nu\alpha}\p_\nu \p_t^3 q)(A^\mu_{\,\,\alpha}\p_\mu\p_t^2 q)\,dy\no\\
+ \int_\Gamma\rho_0^{-1} R'(A^{\nu\alpha}N_\nu \p_t^3 q)(A^\mu_{\,\,\alpha}\p_\mu\p_t^2 q)\,dS+\mathcal{R}.\label{3.16}
\end{align}
The first term on the right hand side of \eqref{3.16} is equal to $$-\frac{d}{dt}\frac{1}{2}\int_\Omega \rho_0^{-1} R'(A^{\nu\alpha}\p_\nu \p_t^{2}q)(A^\mu_{\,\,\alpha}\p_\mu \p_t^{2}q)\,dy+\mathcal{R}$$ and hence moved to the left. In addition, 
\begin{align}
\int_\Gamma\rho_0^{-1} R'(A^{\nu\alpha}N_\nu \p_t^3 q)(A^\mu_{\,\,\alpha}\p_\mu\p_t^{2}q)\,dS= \int_\Gamma\rho_0^{-1} R' \p_t^3(A^{\nu\alpha}N_\nu q)\p_t^{2}(A^\mu_{\,\,\alpha}\p_\mu q)\,dS\no\\
\underbrace{-\int_{\Gamma}\rho_0^{-1}R'\p_t^3(A^{\nu\alpha}N_\nu q)(\p_t A^\mu_{\,\,\alpha})(\p_\mu\p_t q)}_{\mathcal{WB}_1}\underbrace{-{\sum}_{\substack{j_1+j_2=3\\j_1\geq 1}}\int_{\Gamma}\rho_0^{-1} R'\p_t^2(A^\mu_{\,\,\alpha}\p_\mu q)(\p_t^{j_1}A^{\nu\alpha})(N_\nu \p_t^{j_2}q)}_{\mathcal{WB}_2}\nonumber\\
+\underbrace{{\sum}_{\substack{j_1+j_2=3\\j_1\geq 1}}\int_{\Gamma}\rho_0^{-1}R'(\p_t A^\mu_{\,\,\alpha})(\p_\mu \p_t q)(\p_t^{j_1}A^{\nu\alpha})(N_\nu \p_t^{j_2}q)}_{\mathcal{WB}_3},
\end{align}
which is due to
\begin{align}
A^{\nu\alpha}N_\nu\p_t^3 q = \p_t^3(A^{\nu\alpha}N_\nu q)-{\sum}_{\substack{j_1+j_2=3\\j_1\geq 1}}(\p_t^{j_1}A^{\nu\alpha})N_\nu \p_t^{j_2}q,\\
A^\mu_{\,\,\alpha}\p_\mu \p_t^{2}q = \p_t^{2}(A^\mu_{\,\,\alpha}\p_\mu q)-(\p_t A^\mu_{\,\,\alpha})\p_\mu \p_t q.
\end{align}
Next, invoking \eqref{E mod}, \eqref{equivlent R_kk} and \eqref{boundary cond}, the main boundary term is equal to
\begin{align}
 \sigma\int_\Gamma \rr \sqrt{g}g^{ij}\Pi_\mu^\alpha (\p_t^3\cp^2_{ij}\eta^\mu)(\p_t^4 \eta_\alpha)+\underbrace{\sigma{\sum}_{\substack{j_1+j_2=3\\j_1\geq 1}}\int_\Gamma \rr (\p_t^{j_1}\sqrt{g}g^{ij}\Pi_\mu^\alpha)(\cp_{ij}^2\p_t^{j_2} \eta^\mu)(\p_t^3 v_\alpha)}_{\mathcal{WB}_4}\no\\
 = -\sigma \int_\Gamma \rr\sqrt{g}g^{ij}\Pi_\mu^\alpha (\p_t^3\cp_i\eta^\mu)(\cp_j\p_t^4 \eta_\alpha)+\mathcal{WB}_4-\underbrace{\sigma\int_\Gamma \rr \cp_j(\sqrt{g}g^{ij}\Pi_\mu^\alpha)(\p_t^3\cp_i\eta^\mu)(\p_t^3 v_\alpha)}_{\mathcal{WB}_5}.
\end{align}
The first term on the last line is equal to 
$$
-\frac{d}{dt}\frac{\sigma}{2} \int_\Gamma \rr \sqrt{g}g^{ij}\Pi_\mu^\alpha (\p_t^3\cp_i\eta^\mu)(\cp_j\p_t^3 \eta_\alpha)+\underbrace{\frac{\sigma}{2}\int_\Gamma \rr\p_t(\sqrt{g}g^{ij}\Pi_\mu^\alpha) (\p_t^3\cp_i\eta^\mu)(\cp_j\p_t^3 \eta_\alpha)}_{\mathcal{WB}_6},
$$
where the main term is moved to the left, and this completes the construction for \eqref{W}. 

The proof of Theorem \ref{est W thm 1} requires the bound for $\int_0^t||\mathcal{G}_3||_0$ and ${\sum}_{1\leq j\leq 6}\int_0^t\mathcal{WB}_j$.
There is no problem to control $\int_0^t ||\mathcal{G}_3||_0$. In addition, using the duality, we have:
\begin{equation}
\mathcal{WB}_1\lesssim P(||v||_3, ||\eta||_3)||R' \p_t^3(A^{3\alpha}q)||_0||\p_t q||_2,
\end{equation}
and
\begin{align}
\mathcal{WB}_2 \lesssim P(||v||_3, ||\eta||_3)\Big(||(\sqrt{R'}\p_t^2(A^\mu_{\,\,\alpha}\p_\mu q)||_0(\sqrt{R'}\p \p_t^2 v)||_0||q||_2\no\\
+||(\sqrt{R'}\p_t^2(A^\mu_{\,\,\alpha}\p_\mu q)||_0(\sqrt{R'}\p \p_t v)||_1||\p_t q||_2
+||\sqrt{R'}\p_t^2(A^\mu_{\,\,\alpha}\p_\mu q)||_0||\sqrt{R'}q_{tt}||_1\Big).
\end{align}
Therefore, $\int_0^t \mathcal{WB}_1+\mathcal{WB}_2$ can be controlled appropriately. 
Moreover, $\int_0^t \mathcal{WB}_3+\mathcal{WB}_6$ is controlled in a routine fashion. On the other hand, $\int_0^t \mathcal{WB}_4+\mathcal{WB}_5$ is treated in \cite{disconzi2017prioriC}, where the $\rr$-weights are incorporated so that the estimates in \cite{disconzi2017prioriC} can go through.  

\section{The energy identity for $\rr$-weighted wave equations}
\label{C}
We recall that the $\rr$-weighted wave equation reads:
\begin{align}
\rr^\ell R' JD^3\p_t^{2} q - \rr^\ell a^{\nu\alpha}A^\mu_{\,\,\alpha}\p_\nu\p_\mu D^3 q = \widetilde{\mathcal{G}}_4+\widetilde{\mathcal{S}}_4, 
\end{align}
where
\begin{align}
\widetilde{\mathcal{G}}_4= -\rr^{\ell}[D^3\p_t, JR']\p_t q+\rr^\ell [D^3, \rho_0]\p_t (R^{-1}R'\p_t q)
+\rr^{\ell}a^{\nu\alpha}(\p_\nu \rho_0)D^3\p_t v_\alpha\no\\
+\rr^{\ell}a^{\nu\alpha}\p_\nu\big([D^3, A^{\mu}_{\,\,\alpha}]\p_\mu q\big)
+\rr^{\ell}a^{\nu\alpha}\p_\nu\big([D^3,\rho_0]\p_t v_\alpha\big)
-\rr^{\ell}\rho_0[D^3\p_t, a^{\nu\alpha}]\p_\nu v_\alpha,
\end{align}
and
\begin{align}
\widetilde{\mathcal{S}}_4=\rr^{\ell}a^{\nu\alpha}(\p_\nu A^\mu_{\,\,\alpha})\p_\mu D^3 q.
\end{align}
Here, $\ell=1$ when $D^3=\p_t^3$, $\ell=\frac{1}{2}$ when $D^3=\p_t^2 \cp$ and $\ell=0$ when $D^3= \p_t\cp^2$. 
\thm  Let
\label{est W4 thm 1}
\begin{align}
W_4^2= \frac{1}{2} \int_\Omega \rho_0^{-1}\rr^{2\ell}(JR'D^3\p_t q)^2\,dy+\frac{1}{2}\int_\Omega \rho_0^{-1}\rr^{2\ell} R'(A^{\nu\alpha}\p_\nu D^3 q)(A^\mu_{\,\,\alpha}\p_\mu D^3 q)\,dy\no\\
+\frac{\sigma}{2}\int_{\Gamma} \rr^{2\ell+1} \sqrt{g}g^{ij}\Pi_\mu^\alpha(\cp_iD^3\p_t \eta^\mu)(\cp_j D^3 \p_t \eta_\alpha)\,dS. \label{W4 1}
\end{align}
Then, 
\begin{equation}
W_4^2 \leq  \epsilon P(\mathcal{N}) +\PP_0 +\PP\int_0^t\PP, \q t\in[0,T], \label{est W4 1}
\end{equation}
where $T>0$ is sufficiently small.

\paragraph*{Proof of Theorem \ref{est W4 thm 1}}
Invoking \eqref{p_t J} and \eqref{R_kk assumption}, we have:
\begin{align}
\frac{d}{dt}\frac{1}{2}\int_\Omega \rho_0^{-1}\rr^{2\ell}(JR'D^3\p_t q)^2 \,dy = \int_\Omega \rho_0^{-1}\rr^{2\ell}(JR'D^3\p_t q)(a^{\nu\alpha}A^\mu_{\,\,\alpha}\p_\nu\p_\mu D^3 q) \,dy\no\\
+\int_\Omega \rho_0^{-1}\rr^{2\ell}(JR'D^3\p_t q)(\widetilde{\mathcal{G}}_4+\widetilde{\mathcal{S}}_4)\,dy+\mathcal{R},\label{3.13'}
\end{align} 
where $\mathcal{R}$ consists error terms that are generated when $\p_t$ falls on either $J$ or $R'$, which we have no problem to control. In addition, 
\begin{align}
\int_\Omega \rho_0^{-1}\rr^{2\ell}(JR'D^3\p_t q)(a^{\nu\alpha}A^\mu_{\,\,\alpha}\p_\nu\p_\mu D^3 q) \,dy \no\\
= \int_\Omega \rho_0^{-1}\rr^{2\ell}(R'D^3\p_t q)(A^{\nu\alpha}\p_\nu)(A^\mu_{\,\,\alpha}\p_\mu D^3 q)\,dy
-\int_\Omega \rho_0^{-1}\rr^{2\ell}(JR'D^3 \p_t q)\widetilde{\mathcal{S}}_4. \label{3.14'}
\end{align}
The last term in \eqref{3.14'} cancels with the corresponding term in \eqref{3.13'}, which is essential since $||\widetilde{\mathcal{S}}_3||_0$ cannot be controlled uniformly when $R'\to 0$. Moreover, the first term on the right hand side of \eqref{3.14'} is treated as:
\begin{align}
\int_\Omega \rho_0^{-1}\rr^{2\ell}(R'D^3\p_t q)(A^{\nu\alpha}\p_\nu)(A^\mu_{\,\,\alpha}\p_\mu D^3 q)\,dy= -\int_\Omega \rho_0^{-1}\rr^{2\ell} R'(A^{\nu\alpha}\p_\nu D^3\p_t q)(A^\mu_{\,\,\alpha}\p_\mu D^3 q)\,dy\no\\
+ \int_\Gamma\rho_0^{-1}\rr^{2\ell} R'(A^{\nu\alpha}N_\nu D^3\p_t q)(A^\mu_{\,\,\alpha}\p_\mu D^3 q)\,dS+\mathcal{R}.\label{3.16'}
\end{align}
The first term on the right hand side of \eqref{3.16'} is equal to $$-\frac{d}{dt}\frac{1}{2}\int_\Omega \rho_0^{-1}\rr^{2\ell} R'(A^{\nu\alpha}\p_\nu D^3q)(A^\mu_{\,\,\alpha}\p_\mu D^3q)\,dy+\mathcal{R}$$ and hence moved to the left. In addition, 
\begin{align}
 \int_\Gamma\rho_0^{-1}\rr^{2\ell} R'(A^{\nu\alpha}N_\nu D^3\p_t q)(A^\mu_{\,\,\alpha}\p_\mu D^3 q)\,dS= \int_\Gamma\rho_0^{-1}\rr^{2\ell} R' D^3\p_t(A^{\nu\alpha}N_\nu q)D^3(A^\mu_{\,\,\alpha}\p_\mu q)\,dS\no\\
\underbrace{-\int_{\Gamma}\rho_0^{-1}\rr^{2\ell}R'D^3 \p_t(A^{\nu\alpha}N_\nu q)\big([D^3, A^\mu_{\,\,\alpha}]\p_\mu q\big)}_{\widetilde{\mathcal{WB}_1}}\underbrace{-\int_{\Gamma}\rho_0^{-1} R'\rr^{2\ell}D^3(A^\mu_{\,\,\alpha}\p_\mu q)\big([D^3\p_t, A^{\nu\alpha}]N_\nu q\big)}_{\widetilde{\mathcal{WB}_2}}\nonumber\\
+\underbrace{\int_{\Gamma}\rho_0^{-1}R'\rr^{2\ell}\big( [D^3\p_t, A^{\nu\alpha}]N_\nu q\big)\big([D^3, A^\mu_{\,\,\alpha}]\p_\mu q \big)}_{\widetilde{\mathcal{WB}}_3},
\end{align}
which is due to
\begin{align}
A^{\nu\alpha}N_\nu D^3\p_t q = D^3 \p_t(A^{\nu\alpha}N_\nu q)-[D^3\p_t, A^{\nu\alpha}]N_\nu q,\\
A^\mu_{\,\,\alpha}\p_\mu D^3 q = D^3 (A^\mu_{\,\,\alpha}\p_\mu q)-[D^3, A^\mu_{\,\,\alpha}]\p_\mu q.
\end{align}
Next, invoking \eqref{E mod}, \eqref{equivlent R_kk} and \eqref{boundary cond}, the main boundary term is equal to
\begin{align}
 \sigma\int_\Gamma \rr^{2\ell+1} \sqrt{g}g^{ij}\Pi_\mu^\alpha (D^3\p_t\cp^2_{ij}\eta^\mu)(D^3\p_t^2 \eta_\alpha)+\underbrace{\int_\Gamma \rr^{2\ell+1} [D^3\p_t, \sqrt{g}g^{ij}\Pi_\mu^\alpha](\cp_{ij}^2 \eta^\mu)(D^3\p_t v_\alpha)}_{\widetilde{\mathcal{WB}}_4}+\mathcal{R}\no\\
 = -\sigma\int_\Gamma \rr^{2\ell+1} \sqrt{g}g^{ij}\Pi_\mu^\alpha (D^3\p_t\cp_i\eta^\mu)(\cp_jD^3\p_t^2 \eta_\alpha)+\widetilde{\mathcal{WB}}_4-\underbrace{\sigma\int_\Gamma \rr^{2\ell+1} \cp_j(\sqrt{g}g^{ij}\Pi_\mu^\alpha)(D^3\p_t\cp_i\eta^\mu)(D^3\p_t v_\alpha)}_{\widetilde{\mathcal{WB}}_5}+\mathcal{R}.
\end{align}
The first term on the last line is equal to 
$$
-\frac{d}{dt}\frac{\sigma}{2} \int_\Gamma \rr^{2\ell+1} \sqrt{g}g^{ij}\Pi_\mu^\alpha (D^3\p_t\cp_i\eta^\mu)(\cp_jD^3\p_t \eta_\alpha)+\underbrace{\frac{\sigma}{2}\int_\Gamma \rr^{2\ell+1}\p_t(\sqrt{g}g^{ij}\Pi_\mu^\alpha) (D^3\p_t\cp_i\eta^\mu)(\cp_jD^3\p_t \eta_\alpha)}_{\widetilde{\mathcal{WB}}_6},
$$
where the main term is moved to the left, and this completes the construction for \eqref{W}. 

The proof of Theorem \ref{est W4 thm 1} requires the bound for $\int_0^t||\widetilde{\mathcal{G}}_4||_0$ and ${\sum}_{1\leq j\leq 6}\int_0^t\widetilde{\mathcal{WB}}_j$. First, $\int_0^t \widetilde{\mathcal{WB}}_1+\widetilde{\mathcal{WB}}_2$ can be controlled similar to $\int_0^t \mathcal{WB}_1+\mathcal{WB}_2$ in the previous section, after distributing correct $\rr$-weights. Second, the control of $\int_0^t ||\widetilde{\mathcal{G}}_4||_0$ and  $\int_0^t\widetilde{\mathcal{WB}}_3$ can be done in a routine fashion. Finally,  $\int_0^t\widetilde{\mathcal{WB}}_4+\widetilde{\mathcal{WB}}_5+\widetilde{\mathcal{WB}}_6$ is treated similar to $\int_0^t \mathcal{B}$ in section \ref{section B}.
\end{appendix}

\end{document}